\documentclass[12pt]{amsart}
\usepackage{epsf,amscd,amssymb}
 
\newtheorem{thm}{Theorem}[section]
\newtheorem{lem}[thm]{Lemma}

\newtheorem{prop}[thm]{Proposition}

\newtheorem{thmA}{Theorem A}

\newtheorem{thmB}{Theorem B}

\newtheorem{corA}{Corollary A}

\newtheorem{corB}{Corollary B}

\newtheorem{corC}{Corollary C}

\newtheorem{corD}{Corollary D}

\theoremstyle{remark} 
\newtheorem{rem}[thm]{Remark}

\theoremstyle{definition}
\newtheorem{defn}[thm]{Definition}

\newcommand\bR{{\mathbb{R}}}

\newcommand\bZ{{\bf Z}}

\newcommand\GL{{\rm GL}}
\newcommand\SL{{\rm SL}}
\newcommand\PGL{{\rm PGL}}

\newcommand\PO{{\rm PO}}
\newcommand\PSL{{\rm PSL}}

\newcommand\Hom{{\rm Hom}}

\newcommand\dev{{\bf dev}}
\newcommand\SI{{\bf S}}

\newcommand\clo{{\rm Cl}}

\newcommand\ra{\rightarrow}

\newcommand\emp{\emptyset}
\newcommand\eps{\epsilon}

\newcommand\vth{\vartheta}
\newcommand\vpi{\varphi}

\newcommand\ovl{\overline}
\newcommand\Aut{{\rm Aut}}

\newcommand\Out{{\rm Out}}
\newcommand\Inn{{\rm Inn}}

\newcommand\smc{\hbox{\rm ; }}
\newcommand\tri{\triangle}
\newcommand\rpt{\mathbb{RP}^2}
\newcommand\rpts{\rpt{}^*}
\newcommand\rpo{\mathbb{RP}^1}
\newcommand\idnt{\rm{Id}}

\begin{document}
\title[Convex real projective $2$-orbifolds]
{The deformation spaces of convex 
$\rpt$-structures on $2$-orbifolds}
\author{Suhyoung Choi}
\address{Department of Mathematics \\
College of Natural Sciences \\ 
Seoul National University \\ 
151--742 Seoul, Korea}
\email{shchoi@math.snu.ac.kr}
\author{William M. Goldman}
\address{Department of Mathematics\\
University of Maryland \\
College Park, MD, 20742 USA}
\email{wmg@math.umd.edu}
\date{November 15, 2002} 
\thanks{The first author gratefully acknowledges  
support from Korea Research Foundation Grant (KRF-2002-070-C00010).
The second author gratefully acknowledges partial support from
National Science Foundation grant DMS-0103889.}

\begin{abstract}
We determine that the deformation space of convex real 
projective structures, 
that is, projectively flat torsion-free connections
with the geodesic convexity property
on a compact $2$-orbifold of negative Euler characteristic is homeomorphic
to a cell of certain dimension. 
The basic techniques are
from Thurston's lecture notes on hyperbolic $2$-orbifolds,
and the previous work of Goldman on convex real 
projective structures on surfaces. 
\end{abstract}  
\maketitle

%\tableofcontents

\section*{Introduction}
An {\em orbifold\/} is a structure on a topological space 
with neighborhoods modeled on quotient spaces of open
subsets of $\bR^n$ by finite group actions.
A {\em geometric structure\/} modeled on $(X, G)$, 
where $X$ is a manifold and $G$ is a Lie group acting on $X$,
on an orbifold is given by such identifications where the open 
sets are assumed to be in $X$ and the finite groups subgroups of $G$
and the transition functions in $G$. 

An orbifold $Y$ covers another orbifold $Z$ if each point of $Z$
has a neighborhood which is a quotient of a neighborhood of $Y$.

In this paper, we consider only $1$- or $2$-dimensional orbifolds.
They were classified by Thurston along with their hyperbolic,
Euclidean, spherical structures, that is, determined their Teichm\"uller
spaces.  We consider only {\em good} orbifolds those that have
orbifold-covering by surfaces.  We also know that good geometric
orbifolds are often finitely and regularly covered by compact surfaces
by Selberg's lemma: More precisely, if a discrete group $\Gamma$ act
on a surface $S$ properly, then $S/\Gamma$ carries a natural orbifold
structure induced from charts of $S$, where $S/\Gamma$ is said to be a
quotient orbifold.  An orbifold in this paper normally is of such a
form $S/\Gamma$ for a finite group $\Gamma$; such an orbifold
is said to be {\em very good}.

Every good $2$-orbifold is covered by a simply connected 
surface, that is, a disk or a sphere. The cover is said to
be the {\em universal cover} and the group of automorphisms 
the {\em deck transformation group}. The universal cover 
is unique up to orbifold covering automorphisms (see \cite{dgorb} 
and Thurston \cite{Thnote}.)
 
Let $X$ be the real projective plane $\rpt$ and $G$ the group
$\PGL(3, \bR)$ of {\em collineations,\/} (projective automorphisms of
$\rpt)$.  An {\em $\rpt$-structure} or {\em real projective structure}
on a $2$-orbifold $\Sigma$ is 
an $(\rpt, \PGL(3, \bR))$-structure on $\Sigma$.  Two
$\rpt$-structures on $\Sigma$ are equivalent if a so-called isotopy of
$\Sigma$ induces one from the other.  The deformation space
$\rpt(\Sigma)$ of $\rpt$-structures on $\Sigma$ is the space of
equivalence classes of $\rpt$-structures with appropriate topology.

An orbifold $\Sigma$ is {\em closed\/} if it is compact and 
$\partial\Sigma = \emptyset$.

A closed $2$-orbifold $\Sigma$ 
with an $\rpt$-structure is called {\em convex} if it is
projectively diffeomorphic to the quotient of a convex domain in an
affine patch by a properly discontinuous action of a group of
collineations.  (When $\partial\Sigma\neq\emptyset$, 
boundaries are required to be principal geodesic (see \S 3).)
%sc. section number
In particular 
the {\em hyperbolic} $\rpt$-structures arise as quotients of
domains bounded by conics.
They correspond naturally to ones with hyperbolic metrics via the
Klein model of hyperbolic geometry.  % The group of isometries is the
% collineation group of the disk, isomorphic to $\PO(1,2)
% \subset \PGL(3, \bR)$.
% i don't think this is necessary to include here
%sc. OK delete

Let $\Sigma$ be a compact $2$-orbifold with $\chi(\Sigma) < 0$. (Here,
$\chi(\Sigma)$ is the orbifold Euler characteristic.)  The subspace of
the deformation space $\rpt(\Sigma)$ of $\rpt$-structures on $\Sigma$
corresponding to convex ones is denoted by $\mathcal{C}(\Sigma)$ and
the subspace corresponding to hyperbolic ones is denoted by
$\mathcal{T}(\Sigma)$, identified as the Teichm\"uller space of
$\Sigma$ as defined by Thurston \cite{Thnote}.  Then we see that
$\mathcal{T}(\Sigma)$ is a subspace of $\mathcal{C}(\Sigma)$, and
$\mathcal{C}(\Sigma)$ is an open subset of $\rpt(\Sigma)$.

Recall that the orbifold Euler characteristic of orbifolds,
a signed sum of the number of cells with weights given by
$1$ divided by the orders of groups associated to cells.

\begin{thmA}
Let $\Sigma$ be a compact $2$-orbifold with $\chi(\Sigma) < 0$ and
$\partial \Sigma = \emp$.  Then the deformation space
$\mathcal{C}(\Sigma)$ of convex $\rpt$-structures on $\Sigma$ is
homeomorphic to a cell of dimension
\begin{equation*}
-8\chi(X_\Sigma) + (6k_c -2 b_c) + (3 k_r - b_r) 
\end{equation*}
where $X_\Sigma$ is the underlying space of $\Sigma$,
$k_c$ is the number of cone-points, 
$k_r$ the number of corner-reflectors, 
$b_c$ the number of cone-points of order two, and 
$b_r$ the number of corner-reflectors of order two. 
\end{thmA}
The terms should be familiar to geometric topologists but
they will be explained in \S \ref{sec:orbifolds}.
Loosely, a cone-point has a neighborhood which is 
an orbit space of a cyclic action 
on an open disk with a fixed point, with the order being
the order of the group and a corner-reflector has a neighborhood
which is an orbit space of a dihedral group acting on an open disk.  
The order is the half of the group order.

Recall that in the earlier work in Chapter 5: ``Orbifolds and Seifert
fibered spaces'' by W. Thurston \cite{Thnote} proved that the
deformation space of hyperbolic structures on $2$-orbifolds, that is, the
Teichm\"uller space, is homeomorphic to a cell of dimension
\begin{equation*}
-3\chi(X_\Sigma) + 2k + l 
\end{equation*}
where $k$ is the number of cone-points and
$l$ is the number of corner-reflectors.  (We give a detailed proof of
this, Theorem \ref{thm:Thur}, while a sketchy proof is in some
versions of his notes.)  
%Later Ohshika filled in the details
%\cite{Oh:87} but for cases without corner-reflectors. 
Kulkarni-Lee-Raymond \cite{KLR:87} have worked out the Teichm\"uller
spaces also.  Their techniques are not used here but one can possibly try
to apply their method to $\rpt$-structures as well.

The Teichm\"uller space of a $2$-orbifold with negative Euler
characteristic is a real analytic subspace of the deformation space of
convex $\rpt$-structures on the orbifold. This can be seen by a Klein
model of hyperbolic geometry; a hyperbolic structure naturally induces
a canonical convex $\rpt$-structure on the orbifold.

We remark that the deformation spaces of convex $\rpt$-structures on
closed surfaces of negative Euler characteristic are very interesting
spaces admitting complex structures, which were proved using
M\"onge-Ampere equations by Labourie \cite{Lab:97} and Loftin
\cite{Loft:01}, and they might have K\"ahler structures as
conjectured by one of us (Goldman) and even stronger structures as
Labourie suggests.

Loftin \cite{Loft:01} has estimated the dimension of subspaces of
fixed points of automorphism groups acting on the deformations spaces
of convex $\rpt$-surfaces using techniques from differential geometry
and the Riemann-Roch theorem.  We can compute the dimension using
orbifold techniques based on the knowledge of the quotient orbifold.

Assigning a geometric structure to an orbifold $\Sigma$ 
is equivalent to giving its universal cover $\tilde \Sigma$ 
an immersion $\dev:\Sigma \ra \rpt$ equivariant with
respect to the homomorphism $h$ from the (orbifold) fundamental group 
$\pi_1(\Sigma)$ to $\PGL(3, \bR)$. 
The pair $(\dev, h)$ is called a {\em development
pair}, $\dev$ a {\em developing map}, and $h$ the associated {\em holonomy
homomorphism}.  The pair $(\dev, h)$ is only defined up to the action
of $g \in \PGL(3, \bR)$ so that $g(\dev, h(\cdot)) = (g\circ \dev, g
\circ h(\cdot) \circ g^{-1})$.

The map assigning a geometric structure to 
the conjugacy class of the associated holonomy homomorphism 
induces the following map,
so-called holonomy map,
\[\mathcal{H}: \rpt(\Sigma) \ra 
\Hom(\pi_1(\Sigma), \PGL(3, \bR))/\PGL(3, \bR)\]
where $\PGL(3,\bR)$ acts on \[\Hom(\pi_1(\Sigma), \PGL(3, \bR))\] 
by conjugation. 

Let us denote by $C_{\mathcal{T}}(\Sigma)$ the unique component of 
\[\Hom(\pi_1(\Sigma), \PGL(3, \bR))\]
containing the holonomy homomorphisms of hyperbolic 
$\rpt$-structures on $\Sigma$. 
Then $C_{\mathcal{T}}(\Sigma)$ is also a component of 
the part
\[\Hom(\pi_1(\Sigma), \PGL(3, \bR))^{st}\]
of \[\Hom(\pi_1(\Sigma), \PGL(3, \bR))\] where 
$\PGL(3, \bR)$ acts properly. 
$C_{\mathcal{T}}/\PGL(3,\bR)$ is said to be 
a {\em Hitchin-Teichm\"uller component} (see \cite{Hit:92}).
We prove:

\begin{thmB}
Let $\Sigma$ be a closed $2$-orbifold with negative Euler
characteristic. Then 
\[\mathcal{H}:\mathcal{C}(\Sigma) \ra C_{\mathcal{T}}(\Sigma)/\PGL(3, \bR)\]
is a diffeomorphism, and $C_{\mathcal{T}}(\Sigma)$ consists of 
discrete faithful representations of $\pi_1(\Sigma)$.
%In this case, $C_{\mathcal{T}}(\Sigma)$ is a component of 
%$\Hom(\pi_1(\Sigma), \PGL(3, \bR))$.
%sc. delete
\end{thmB}

Hence, by Theorem A,
$\mathcal{C}(\Sigma)$ is homeomorphic to a cell of dimension
$-8\chi(X_\Sigma) + 6 k_c -2 k_b + 3 l_c - l_b$ where
$X_\Sigma$ is the underlying space of $\Sigma$, 
$k_c$ is the number of cone-points, $l_c$ the number
of corner-reflectors, $k_b$ the number of cone-points of
order two, and $l_b$ the number of corner-reflectors of
order two. 

\begin{corA} 
The Hitchin-Teichm\"uller component
$C_{\mathcal{T}}(\Sigma)/\PGL(3, \bR)$ 
is homeomorphic to a cell of the dimension as above.
\end{corA}
This gives us a partial classification of discrete representations 
of fundamental groups using topological ideas:
Benoist \cite{Ben:00} characterized the group of 
projective transformations acting on a convex domain
for general dimensions. (We only consider 
the three-dimensional cases here.)
An element of $\GL(3, \bR)$ is {\em proximal}
if it has an attracting fixed point in $\rpt$
for its standard action. 
An element is {\em positive proximal} if the eigenvalue
corresponding to the fixed point is positive. 
A subgroup $\Gamma$ of $\GL(3, \bR)$ is {\em positive proximal}
if all proximal elements of $\Gamma$ are positive 
proximal. Proposition 1.1 of \cite{Ben:00} shows that 
if $\Gamma$ is an irreducible subgroup of $\GL(3, \bR)$,
then $\Gamma$ preserves a properly convex cone in $\bR^3$ 
if and only if $\Gamma$ is positive proximal. 
Such a subgroup $\Gamma$ of $\GL(3, \bR)$, if discrete, acts on
a convex domain $\Omega$ in an affine patch
so that $\Omega/\Gamma$ is a $2$-orbifold. Suppose that
$\Omega/\Gamma$ is compact, and $\Gamma$ 
contains a free subgroup of two generators, then 
$\Omega/\Gamma$ is an orbifold of negative Euler-characteristic.
By Theorem A, such groups are parameterized by cells.

We discuss the rigidity of $2$-orbifolds with hyperbolic or convex
$\rpt$-structures in a way related to the interesting recent work of
Dunfield and Thurston \cite{DT}.  The $2$-orbifolds with Teichm\"uller
spaces single points must be orbifolds with empty boundary and contains
no $1$-dimensional suborbifold cutting them into smaller orbifolds
with negative Euler characteristics.  From the classifications of such
orbifolds, we see that such an orbifold is a sphere with three
cone-points, a disk with one cone-point and one corner-reflector, and
a disk with three corner-reflectors. The proofs of the following 
corollaries are omitted as they follow from 
Corollary A and Theorem \ref{thm:Thur}.

\begin{corB} 
The sphere $\Sigma$ with cone-points of order $p, q, r$ 
satisfying $p\leq q \leq r, 1/p + 1/q + 1/r < 1$ has as its Teichm\"uller 
space a single point. 
If $p =2$, then so is $\mathcal{C}(\Sigma)$. 
If $p > 2$, then $\mathcal{C}(\Sigma)$ is homeomorphic to 
$\bR^2$.
\end{corB}

\begin{corC} 
Let $\Sigma$ be a $2$-orbifold whose underlying space is 
a disk and with one cone point of order $p$ and 
a corner-reflector of order $q$ so that $1/p + 1/2q < 1/2$
has as its Teichm\"uller space a single point. 
If $q = 2$, then so is $\mathcal{C}(\Sigma)$. 
If $q > 2$, then $\mathcal{C}(\Sigma)$ is homeomorphic to 
$\bR$. 
\end{corC}

\begin{corD} 
Let $\Sigma$ be a $2$-orbifold whose underlying space is 
a disk and with three corner-reflectors of order $p\leq q \leq r$, 
$1/p + 1/q + 1/r < 1/2$. Then $\mathcal{T}(\Sigma)$ is 
a single point. If $p=2$, then so is $\mathcal{C}(\Sigma)$.
If $p >2$, then $\mathcal{C}(\Sigma)$ is homeomorphic to 
$\bR$. 
\end{corD}

The idea of the proof of Theorem A follows \cite{Gconv:90} closely: 
given a compact $2$-orbifold of negative Euler characteristic, 
we find ``essential'' $1$-orbifolds decomposing it into 
twelve types of ``elementary'' $2$-orbifolds of 
negative Euler characteristic 
which can no longer be decomposed. Given one of these elementary ones, 
we determine the deformation space with the projective structures on
the boundary $1$-orbifolds fixed. We realize the deformation space 
as a fibration over the deformation space of the union of the boundary 
$1$-orbifolds. Next, we rebuild the original orbifold 
by various geometric constructions. In the constructions, 
we build the deformation space again using the fibration 
property.

To prove Theorem B, we follow the article \cite{CG:93}: 
Given an $\rpt$-structure, we associate a homomorphism 
to $\PGL(3, \bR)$ from the fundamental group of an orbifold. 
Since the homomorphism is determined only up to 
conjugation, we show that the deformation space 
is locally homeomorphic to the $\PGL(3,\bR)$-quotient of the
space of representations of the fundamental group to $\PGL(3, \bR)$.
We will show that 
the image of the deformation space of convex $\rpt$-structures is 
an open and closed subset. By Theorem~A, this space is 
connected, implying Theorem B.

\S 1 contains preliminary materials on topological 
$2$-orbifolds.  We discuss the inverse topological processes of splitting
and sewing  
$2$-orbifolds along $1$-orbifolds. 
% \marginpar{sc: sewing, sewing on... correct usage?}

\S 2 concerns $2$-dimensional orbifolds with projectively flat
structures, henceforth called $\rpt$-orbifolds. 
We discuss the properties of the deformation space and
the character variety of the fundamental group of 
a compact $2$-orbifold $\Sigma$. 
We define 
{\em the Hitchin-Teichm\"uller component  \/}
of $\Hom(\pi_1(\Sigma), G)/G$ where $G = \PGL(3, \bR))$.
The deformation space of convex $\rpt$-structures on
$\Sigma$ identifies with an open and closed subset of
$\Hom(\pi_1(\Sigma), G)/G$.

\S 3 details the geometric processes of splitting and sewing
$\rpt$-oribifolds.  We discuss how these geometric processes affect the
deformation spaces of $\rpt$-structures.

\S 4 describes the decomposition of convex $\rpt$-orbifolds of negative
Euler characteristic into elementary ones of twelve types.  They are
{\em elementary \/} in that they cannot be further split along 
$1$-orbifolds into ones of negative Euler characteristic.  
This follows the work 
of Thurston \cite{Thnote} for hyperbolic $2$-orbifolds.

\S 5 determines the Teichm\"uller space of elementary
$2$-orbifolds, that is, the deformation spaces of hyperbolic structures.
This is used later to show the existence of convex structures on 
elementary orbifolds.

\S 6 computes the deformation space of convex
$\rpt$-structures on elementary $2$-orbifolds. Some of the elementary
$2$-orbifolds are classified by decomposing them into unions of
triangles, and reducing these to configurations of triangles and the
corresponding easily solvable algebraic relations, following \cite{Gconv:90}.
For others, we will identify the deformation spaces with
other types of configuration spaces. 

The authors benefited much from conversations with 
Thierry Barbot,
Yves Benoist,
Nathan Dunfield, Vladimir Fock, David Fried, 
Hyuk Kim, Inkang Kim,
Fran{\c c}ois Labourie, John Loftin,
John Millson, William Thurston, and S.-T. Yau. 
The final writing of this paper was done while one 
of us (Choi) was visiting Departments of 
Mathematics at Boston University and Stanford University. 
We appreciate very much the hospitality of the both departments.
The figures were drawn using xfig and Maple.
%%%!
\section{Preliminaries on orbifolds}
\label{sec:orbifolds}

We define orbifolds, covering maps, orbifold maps, suborbifolds, and
Euler characteristics of orbifolds.  We discuss how to obtain
orbifolds by cutting along $1$-orbifolds and how to sew along
$1$-orbifolds to obtain bigger orbifolds.  Euler characteristic zero
$2$-orbifolds and regular neighborhoods of $1$-orbifolds are defined
and classified.  The topological operations will be given three
different interpretations.

The material here can 
be found principally in Chapter 5 of 
various versions of Thurston's note \cite{Thnote} and an
expository paper on $2$-dimensional orbifolds by Scott
\cite{Scott:83}.  See also Ratcliffe \cite{RT:94}, 
Bridson-Haefliger \cite{BH:99}, and Kapovich \cite{Kap:01}. 
Most of the technical 
details for this paper can be found in another paper
\cite{dgorb} on geometric structures on orbifolds, 
which should be read ahead of this paper.

%%% I changed several ``homeomorphisms'' to ``diffeomorphisms''
%%% we claim to work with differentiable objects
%%% I think it is important that we be consistent.

In this paper, we will only work with differentiable objects although
we won't require differentiability for the spaces of such objects.
Moreover, we assume that group actions are %%%%! locally faithful.
strongly effective.  That is, if an element $g$ is so that $g$ equals
identity on an open subset, then $g$ is the identity element.  
%(We
%believe that there is some topological, PL, smooth structure theories
%on orbifolds, and think that there are no essential differences
%between these theories at least in dimension two and three.)

\subsection{Definition of orbifolds}
Let $Q$ be a Hausdorff, second countable space. 
An {\em orbifold atlas } is a open covering
$\{U_i\}_{i\in I}$ that for each $U_i$,
there is an open subset $\tilde U_i$ of $\bR^n$ and a finite group
$\Gamma_i$ of diffeomorphisms acting on $\tilde U_i$ with a
homeomorphism $\phi_i: \tilde U_i/\Gamma_i \ra U_i$. Given an
inclusion map $U_i \ra U_j$, there is an injective homomorphism
$f_{ij}: \Gamma_i \ra \Gamma_j$ and an embedding $\tilde \phi_{ij}:
\tilde U_i \ra \tilde U_j$ equivariant with respect to $f_{ij}$ (that is,
$\tilde \phi(g x) = f_{ij}(g)\circ \tilde \phi_{ij}(x)$ for all $g \in
\Gamma_i, x \in \tilde U_i$) so that that the following diagram is
commutative:
\begin{eqnarray}
\tilde U_i  &\stackrel{\tilde \phi_{ij}}{\longrightarrow} 
& \tilde U_j \nonumber\\
\downarrow &   &\downarrow  \nonumber \\
\tilde U_{i}/\Gamma_i &\stackrel{\phi_{ij}}{\longrightarrow}&
\tilde U_j/f_{ij}\Gamma_i \nonumber \\ 
\phi_i \downarrow & & {\begin{array}{c} \downarrow \nonumber \\ 
\tilde U_j/\Gamma_j \\ \downarrow \phi_j \end{array}}\\
U_i &\subset & U_j 
\end{eqnarray}
where $\phi_{ij}$ is induced from $\tilde \phi_{ij}$.
Actually, $(\tilde \phi_{ij}, f_{ij})$ is determined 
up to the action given by
\[g(\tilde \phi_{ij}, f_{ij}(\cdot)) = (g \circ \tilde \phi_{ij}, 
g f_{ij}(\cdot) g^{-1}) \hbox{ for } g \in \Gamma_j.\]
That is, the equivalence class of the pair is given in
the information about the orbifold structure.
% (We work in the differentiable category.)
%sc. delete
An {\em orbifold structure} is a maximal 
family of coverings satisfying the above conditions. 
The space $Y$ with the structure is an {\em $n$-dimensional 
orbifold} or {\em $n$-orbifold}. 
Given an orbifold $Q$ the underlying space is 
denoted $X_Q$. Clearly, a smooth structure on a manifold 
is an orbifold structure. 

An {\em orbifold\/} is a topological space with a maximal orbifold atlas.

Given a Lie group $G$ acting on a space $X$, we define 
an {\em $(X, G)$-structure} on an orbifold $\Sigma$ to
be a maximal collection 
\[\{(U_i, \tilde U_i, \Gamma_i, \phi_i), (f_{ij}, \tilde \phi_{ij})\}\]
in the orbifold structure
where $\tilde U_i$ is identified with an open subset of $X$ and 
$\Gamma_i$ and $f_{ij}$ are 
required to be restrictions of elements of $G$
and $\tilde \phi_{ij}$ conjugation homomorphisms by $f_{ij}$s.

Here is the basic example. Suppose $X$ is a homogeneous Riemannian manifold 
with isometry group $G$ and let $\Gamma\subset G$ be a discrete subgroup
acting  on an open subset $\Omega/\Gamma$.
Then the quotient orbifold $\Omega/\Gamma$ carries an $(X, G)$-structure. 

A point of $X_Q$ is said to be {\em regular} if it has a neighborhood 
with trivial associated group. Otherwise it is {\em singular}.
Let $p$ be a singular point. Then for every chart 
$(\tilde U, \Gamma)$ in the orbifold atlas about $p$, the point
$p$ corresponds to a fixed point of a nontrivial element of $\Gamma$.

An {\em orbifold with boundary} has neighborhoods modeled
on open subsets of the closed upper half space $\bR^{n, +}$.
A {\em suborbifold} $Q'$ on a subspace 
$X_{Q'}\subset X_Q$ is the subspace so that each point of $X_{Q'}$ has
a neighborhood in $X_Q$ modeled on an open subset $U$ of 
$\bR^n$ with a finite group $\Gamma$ preserving $U \cap \bR^d$ where 
$\bR^d \subset \bR^n$ is a proper subspace, so that 
$(U\cap \bR^d, \Gamma')$ is in the orbifold structure of $Q'$.
Here $\Gamma'$ denotes the restricted group of $\Gamma$ 
to $U \cap \bR^d$, which is in general a quotient group.
%\marginpar{I need to check, M. Davis or Kato?,}

The {\em interior} of $Q$ is defined as the set of
points with neighborhoods modeled on open subsets of $\bR^n$.
The {\em boundary} of $Q$ is the complement of
the interior. The boundary is denoted by $\partial Q$. 
The boundary of an $n$-orbifold is clearly
an $(n-1)$-suborbifold without boundary.
(It is a subset of the boundary of the underlying
space $X_Q$ but is not necessarily all of it.) 

\begin{defn} A map $f: Q \ra Q'$ between two orbifolds $Q$ and $Q'$ 
is an {\em orbifold map} or simply a {\em map}, 
if it induces a continuous function
$X_Q \ra X_{Q'}$, and for each point of form $f(x)$ for $x \in X_Q$ with 
$(U, \Gamma)$ with a homeomorphism
$\phi_U$ from $U/\Gamma$ to a neighborhood of $f(x)$
there is a pair $(V, \Gamma')$ with homeomorphism $\phi_V$ from $V/\Gamma'$
to a neighborhood of $x$ so that there exists a differentiable 
map $\tilde f: V \ra U$ 
equivariant with respect to a homomorphism $\psi: \Gamma' \ra \Gamma$ 
so that the following diagram is commutative:
\begin{eqnarray}
V  &\stackrel{\tilde f}{\longrightarrow} & U  \nonumber \\
\downarrow & & \downarrow  \nonumber \\
V/\Gamma'  & \stackrel{f}{\longrightarrow} & U/\Gamma.
\end{eqnarray}
That is, we need to record $(\tilde f, \psi)$ but $\tilde f$ 
is determined only up to the actions of $\Gamma$ and $\Gamma'$
and $\psi$ changed correspondingly.
\end{defn}

%%4/24 5:20

Let $I$ be the unit interval with an obvious smooth structure 
seen as an orbifold structure. 
Given an orbifold $Q$, $X_Q \times I$ has an obvious orbifold structure 
with boundary equal to 
the union of two orbifolds $Q\times \{0\}$ and $Q\times \{1\}$,
orbifold-diffeomorphic to $Q$ itself. 
Let $Q \times I$ denote the orbifold.
A {\em homotopy} between two orbifold maps from $Q$ to another orbifold $Q'$ 
is an orbifold map $Q \times I$ to $Q'$ which restricts 
to the two maps at $0$ and $1$. 

An {\em isotopy} of an orbifold $\Sigma$ is 
a self-diffeomorphism $f$ that so that there exists an orbifold map 
$F:\Sigma \times I \ra \Sigma$ so that 
$F_t:\Sigma \ra \Sigma$ given by $F_t(x) = F(x, t)$ 
is a diffeomorphism for each $t$ and 
$F_0$ is the identity and $F_1 = f$. 

\begin{defn}\label{defn:cov} 
A {\em covering orbifold}\/ of an orbifold $Q$ is an orbifold
$\tilde Q$ with a surjection $p:X_{\tilde Q} \ra X_Q$ such that 
each point $x \in X_Q$ has a neighborhood $U$ with
a homeomorphism $\phi:\tilde U/\Gamma \ra U$ for
an open subset of $\tilde U$ in $\bR^n$ or $\bR^{n, +}$ with
a group $\Gamma$ acting on it  so that 
each component $V_i$ of $p^{-1}(U)$ has a diffeomorphism
$\tilde \phi_i:  \tilde U/\Gamma_i \ra V_i$ (in the orbifold
structure) where $\Gamma_i$ is a subgroup of $\Gamma$.
We require the quotient map $\tilde U \ra V_i$ induced by $\tilde \phi_i$ 
composed with $p$ is the quotient map $\tilde U \ra V$ induced by $\phi$.
\end{defn}

Clearly, if $f:Q' \ra Q$ is an orbifold covering and 
$Q$ has an $(X, G)$-structure, then so does $Q'$.

Given a smooth manifold $M$ and a group $\Gamma$ acting
properly discontinuously, $M/\Gamma$ has a unique orbifold
structure for which the quotient projection $M\longrightarrow M/\Gamma$
is an orbifold covering map.

\begin{defn} A {\em good} orbifold is an orbifold which
has a covering that is a manifold.
\end{defn}

The inverse image of a suborbifold under the covering
map of an orbifold is a suborbifold again. If the covering 
orbifold is a manifold, then the inverse image is a submanifold.

We remark that an orbifold always has 
a so-called universal covering orbifold:
\begin{prop} 
An orbifold $Q$ has a covering orbifold $p:\tilde Q \ra Q$
with the following property. If $x$ is a nonsingular point,  
$p(\tilde x) = x$ for $\tilde x \in \tilde Q$, and $p':Q'\ra Q$ 
is a covering map with $p'(x') = x$, then there is a lifting 
orbifold map $q:\tilde Q \ra Q'$ with $q(\tilde x) = x'$. 
\end{prop}
\begin{proof} See \cite{dgorb} or \S 5 of Thurston \cite{Thnote} or
Chapter 13 of Ratcliffe \cite{RT:94}.
\end{proof}

A universal covering orbifold is unique up to isomorphisms
of covering spaces; that is, given two universal coverings
$p_1: \tilde Q_1 \ra Q$ and $p_2: \tilde Q_2 \ra Q$, 
there is a diffeomorphism $f: \tilde Q_1 \ra \tilde Q_2$ 
so that $p_1 \circ f = p_2$. 

For good orbifolds, the universal covering orbifolds are simply
connected manifolds. For two-dimensional orbifolds, they are
diffeomorphic to either a disk or a sphere.
 
Let $p:Q' \ra Q$ be an orbifold covering map.  A {\em deck
transformation} of a covering orbifold $Q'$ of $Q$ is an orbifold
self-diffeomorphism of $Q'$ which composed with $p$ is equal to
$p$. When $Q'$ is the universal cover of $Q$, then the group of deck
transformations are said to be the ({\em orbifold}) 
{\em fundamental group} of $Q$ and denoted by $\pi_1(Q)$.
(See Chapter 13 of Ratcliffe \cite{RT:94}.)  Thus, a good
orbifold is a quotient orbifold of a simply connected manifold 
by the fundamental group. We denote by $\pi_1(Q)$ the group
of deck transformations.

Let $M$ and $M'$ be two orbifolds, and let 
$\tilde M$ and $\tilde M'$ be their universal covers 
with deck transformations $\pi_1(M)$ and $\pi_1(M')$ respectively.
Given a map $f: \tilde M \ra \tilde M'$ lifting 
a diffeomorphism $f':M \ra M'$, define a homomorphism
\begin{align*}
\tilde f_*: \pi_1(M) & \ra \pi_1(M') \\ 
\gamma & \mapsto \tilde f \circ \gamma \circ \tilde f^{-1}.
\end{align*}

A singular point of $1$-orbifold has always a group $\bZ_2$ associated
with it which acts as a reflection.  We call this singular point a
{\em mirror point}.

We can easily classify $1$-orbifolds with compact connected underlying
spaces. Each of them is diffeomorphic to
a circle, a segment with both endpoints a mirror point,
a segment with two endpoints one of which is a mirror point,
or a segment without singular points. The second one is said to be 
a {\em full $1$-orbifold}, the third {\em half} $1$-orbifold,
and the fourth a {\em segment}.

\begin{defn} 
The singular points and the boundary points
of $1$-orbifolds are said to be {\em endpoints}.
\end{defn}

Note that all $1$-orbifolds are very good:

The singular points of a two-dimensional orbifold fall into three types:
\begin{itemize}
\item[(i)] The mirror point: $\bR^2/\bZ_2$ where $\bZ_2$ acts by
reflections on the $y$-axis. 
\item[(ii)] The cone-points of order $n$: $\bR^2/\bZ_n$ where 
$\bZ_n$ acting by rotations by angles $2\pi m/n$ for integers $m$.
\item[(iii)] The corner-reflector of order $n$: $\bR^2/D_n$ where 
$D_n$ is the dihedral group generated by reflections about
two lines meeting at an angle $\pi/n$. 
\end{itemize}
(The actions here are isometries on $\bR^2$.)

\begin{figure}[ht]
\centerline{\epsfxsize=2.5in \epsfbox{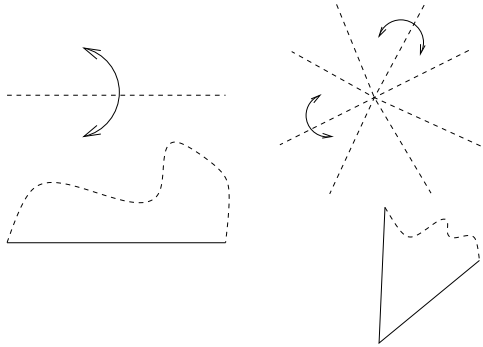}}
\caption{\label{fig:singular} 
The singular points in two-dimensional orbifolds.
Mirror points and corner-reflectors are drawn.}
\end{figure}
\typeout{<>}

\subsection{The Euler characteristics of orbifolds}

In dimension $1$ or $2$, the underlying space $X_Q$ of an orbifold $Q$ has 
a cellular decomposition such that each point of an open cell has 
the same model open set and the same finite group action.
We define the Euler characteristic to be 
\[\chi(Q) = \sum_{c_i} (-1)^{\dim(c_i)} (1/|\Gamma(c_i)|),\]
where $c_i$ ranges over the open cells and $|\Gamma(c_i)|$ is
the order of the group $\Gamma_i$ associated with $c_i$.

For example, a full $1$-orbifold has Euler characteristic zero.

We recall that the cardinality of inverse image under a covering map
$p:Q' \ra Q$ is constant over nonsingular points.
If $p:Q' \ra Q$ is $k$-sheeted, then $\chi(Q') = k\chi(Q)$.
(This follows since over the regular points, the map is an ordinary
covering map.)

Suppose that a $2$-orbifold $\Sigma$ without boundary 
has the underlying space 
$X_\Sigma$ and $m$ cone-points of order $q_i$ and $n$ corner-reflectors
of order $r_j$. Then
the following generalized Riemann-Hurwitz formula is very useful also:
\begin{equation}\label{eqn:riehur}
\chi(\Sigma) = \chi(X_\Sigma) - \sum_{i=1}^m \left( 1- \frac{1}{q_i} \right) 
- \frac{1}{2}\sum_{j=1}^n\left( 1 - \frac{1}{r_j} \right).
\end{equation}
If $\Sigma$ has nonempty boundary, then it is easy to show 
that the boundary consists of circles and full $1$-orbifolds,
which are mutually disjoint suborbifolds.
Let $n_\Sigma$ denote
the number of boundary full $1$-orbifolds of $\Sigma$, and 
we obtain 
\begin{equation}\label{eqn:riehur2}
\chi(\Sigma) = \chi(X_\Sigma) - \sum_{i=1}^m \left( 1- \frac{1}{q_i} \right) 
- \frac{1}{2}\sum_{j=1}^n\left( 1 - \frac{1}{r_j} \right)
- \frac{1}{2} n_\Sigma.
\end{equation}
Again, this formula is proved by a doubling argument. 
(See Thurston \cite{Thnote} or Scott \cite{Scott:83} for
details.)

For $2$-orbifolds $\Sigma_1, \Sigma_2$ meeting in 
a compact $1$-orbifold 
$Y$ forming a $2$-orbifold $\Sigma$ as a union, 
we have the following additivity formula: 
\begin{equation}\label{eqn:additivity} 
\chi(\Sigma) = \chi(\Sigma_1) + \chi(\Sigma_2) - \chi(Y),
\end{equation}
to be verified by counting cells with weights
since the orders of singular points in the boundary orbifold 
equal the ambient orders.

Thurston showed that compact $2$-orbifolds of nonpositive Euler
characteristic with or without boundary admit euclidean or hyperbolic
structures (with geodesic boundary when there is a nonempty boundary).
See Theorem \ref{thm:Thur}. They are good orbifolds; that is, they are
covered by surfaces, since orbifolds admitting geometric structures
are good again by Thurston \cite{Thnote} (see \cite{dgorb} also).  We
also claim that these orbifolds admit finite regular-covers by
surfaces; that is, they are very good.  In other words, there is a finite
group $F$ acting on a surface $S$ such that our orbifold is of form
$S/F$, a quotient orbifold:  Since such an orbifold is of form the
hyperbolic plane $H^2$ quotient by an infinite discrete group
$\Gamma$. There is a finite-index torsion-free normal subgroup
$\Gamma'$ by Selberg's lemma.  Thus $H^2/\Gamma'$ is a closed surface
and our orbifold $H^2/\Gamma$ is a finite quotient orbifold of it.

Since $\pi_1(S) = \Gamma$ is finitely presented, we obtain:
\begin{thm}\label{thm:finitegen} 
Let $\Sigma$ be a compact orbifold of negative Euler characteristic. 
Then the group $\pi_1(\Sigma)$ of deck-transformations 
is finitely presented. 
\end{thm}

\subsection{Splitting and sewing on topological 
$2$-orbifolds}

We now describe the process of splitting and sewing: Let $S$ be a very
good orbifold so that its underlying space $X_S$ is a pre-compact open
surface with a path-metric admitting a compactification to a surface
with boundary, which amounts to attaching boundary components, but not
points.  (For example, $S$ maybe a suborbifold obtained by removing
from an orbifold $S'$ boundary components or embedded $1$-orbifolds or
circles. The completion is in general different from $X_{S'}$.)

Let $\hat S$ be a very good cover, that is, a finite regular cover, of
$S$, so that $S$ is orbifold-diffeomorphic to $\hat S/F$ where $F$ is
a finite group acting on $\hat S$.  Since $X_{\hat S} =\hat S$ is also
pre-compact and has a path-metric, complete it to obtain a compact
surface $X'_{\hat S}$.  Since $F$ acts isometric with respect to some
metric of $\hat S$, it acts on $X'_{\hat S}$. We easily see that $X_S$
can be identified naturally as a subspace of the $2$-nd countable
Hausdorff space $X'_{\hat S}/F$ which has a quotient orbifold
structure induced from $X'_{\hat S}$ by the quotient-orbifold process
described briefly.

\begin{defn}\label{defn:completion}
$X'_{\hat S}/F$ with the quotient orbifold
structure is said to be the {\em orbifold-completion} of $S$.
\end{defn}

Let $S$ be a $2$-orbifold with an embedded circle or a full
$1$-orbifold $l$ in the interior of $S$.  The completion $S'$ of $S-l$
is said to be obtained from {\em splitting} $S$ along $l$.  Since
$S-l$ has an embedded copy in $S'$, we see that there exists a map $S'
\ra S$ sending the copy to $S-l$.  Let $l'$ denote the boundary
component of $S$ corresponding to $l$ under the map.  Conversely, $S$
is said to be obtained from {\em sewing} $S'$ {\em along} $l'$.

%Note that splitting can be done for $2$-orbifolds of nonpositive Euler
%characteristic if there exists an embedded full $1$-orbifold.  Sewing
%can be done for any bounded $2$-orbifold of nonpositive Euler
%characteristic.
%sc. delete

If the interior of the underlying space of $l$ lies in the interior of
the underlying space of $S$, then the components of $S'$ are said to
be {\em decomposed components of $S$ along $l$}, and we also say that
$S$ {\em decomposes} into $S'$ along $l$.  Of course, if $l$ is a
union of disjoint embedded circles or full $1$-orbifolds, the same
definition holds.

%\subsection{Silvering and Clarifying}

%\begin{rem}\label{rem:freeorb}
There are two distinguished classes of splitting and sewing operations:
A simple closed curve boundary component can be made into a set of
mirror points and conversely in a unique manner: a boundary point has
a neighborhood which is realized as a quotient of an open ball by a
$\bZ_2$-action generated by a reflection about a line. Such a system
of model neighborhoods can be chosen consistently to produce an
orbifold structure.  A boundary full $1$-orbifold can be made into a
$1$-orbifold of mirror points and two corner-reflectors of order two
and conversely in a unique manner: the interior points of the
$1$-orbifold have neighborhoods as above, and a boundary point has a
neighborhood which is a quotient space of a dihedral group of order
four acting on the open ball generated by two reflections.  Again,
such a system produces an orbifold structure.  The forward process is
called {\em silvering} and the reverse process {\em clarifying}. 
%In this paper, we never clarify the original mirror points as we
%decompose (unless we silvered them in artificially for a purpose). See
%the last part of this section for more details.

%Such a process is always possible topologically, but not always
%geometrically (see \S3).
%\end{rem}
%C. delete

\subsection{Euler-characteristic-zero $2$-orbifolds}

An {\em edge} is a segment in the singular locus of a $2$-orbifold
which ends in corner-reflectors or in the boundary. An edge is a
$1$-orbifold only if its endpoints are corner-reflectors of order two
or boundary points. 
%An edge can be clarified if and only if it is a
%$1$-orbifold.
%C. delete

Let $A$ be a compact annulus with boundary.  The quotient orbifold of
an annulus has Euler characteristic zero.  From equation
\eqref{eqn:riehur2}, we can determine all of the Euler characteristic
zero $2$-orbifolds with nonempty boundary. We call them the {\em
annular} orbifolds.  They are quotients of annuli.  Each of
them is diffeomorphic to one of the following orbifolds:
\begin{itemize}
\item[(1)] an annulus,
\item[(2)] a M\"obius band,
\item[(3)] an annulus with one boundary component silvered
(a {\em silvered annulus}),
\item[(4)] a disk with two cone-points of order two
with no mirror points ( a {\em $(;2,2)$-disk} from Thurston's notation),
\item[(5)] a disk with two boundary $1$-orbifolds, 
two edges ({\em a silvered strip}), 
\item[(6)] a disk with one cone-point and one boundary 
full $1$-orbifold  (a {\em bigon with a cone-point of order two}),
that is, it has only one edge, and
\item[(7)] a disk with two corner-reflectors of order two
and one boundary full $1$-orbifold (a {\em half-square}). (It has three 
edges.)
\end{itemize} 
To prove this simply notice that the underlying space must
have a nonnegative Euler characteristic.
When the Euler characteristic of the space is zero,
then there are no cone-points, corner-reflectors, and 
the boundary full $1$-orbifolds, since they will make 
the Euler characteristic negative.
Hence the first three occur.
Suppose that the underlying space is a disk. 
If there are no singular points in the boundary,
then we obtain the fourth case as there has to
be exactly two cone-points of order two for
the Euler characteristic to be zero. 
If there are two boundary full $1$-orbifolds, then
there are no singular points in the interior and 
no corner-reflector can exist; thus, we have 
the fifth case. Assume that there exists exactly
one boundary full $1$-orbifold.
If there is a cone-point, then it has to be 
a unique one and of order two. Thus, we have the
sixth case. If there are no cone-points, but
corner-reflectors, then there are exactly two
corner-reflectors of order two and no more.
We have the seventh case. 

\begin{figure}[ht]
\centerline{\epsfxsize=3.2in \epsfbox{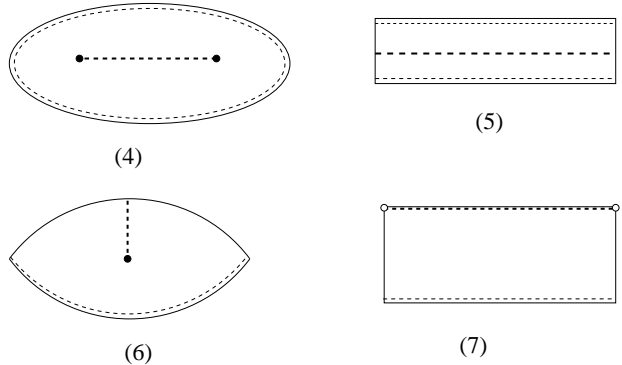}}
\caption{\label{fig:eulerzero} 
Orbifolds of zero-Euler-characteristic in (4)-(7).
A thin dashed arc indicates boundary and a thick dashed arc 
a suborbifold $l$. 
A black dot indicates a cone-point of
order two and a white dot a corner-reflector of order two.}
\end{figure}
\typeout{<>}

\subsection{Regular neighborhoods of $1$-orbifolds}
%510 0:01
Suppose that there exists a circle or a $1$-orbifold $l$
embedded in the interior of a $2$-orbifold $S$,
and assumed not to be homotopic to points.
$l$ has Euler characteristic zero. Thus in 
a good cover $\hat S$ of $S$, the inverse image
of certain of its neighborhoods
is a disjoint union of annuli or M\"obius bands. Thus, $l$ has 
a neighborhood of zero Euler characteristic:
Since the inverse image
of $l$ consists of closed curves which represent
generators of the fundamental group of the neighborhoods,
it follows that in the first two cases (1) and (2), $l$ 
is the closed curve representing the generator of
the fundamental group; in case (3), $l$ is 
the mirror set that is a boundary component; in case (4),
$l$ is the arc connecting the two cone-points unique up to
homotopy; in case (5), $l$ is an arc connecting
two interior points of two edges respectively; 
in case (6), $l$ is an arc connecting 
an interior point of an edge and 
the cone-point of order two; and in the final case (7), 
the edge in the topological boundary connecting the 
two corner-reflectors of order two. 

%%4/25 3:00
\begin{defn}\label{defn:regularn}
Given a $1$-orbifold $l$ and a neighborhood $N$ of it 
in some ambient $2$-orbifold, $N$ is said to be 
a {\em regular neighborhood} if the pair $(N, l)$ is 
diffeomorphic to one of the above. 
\end{defn}

\begin{prop}\label{prop:regularn}
A $1$-orbifold in a good $2$-orbifold has a regular 
neighborhood which is unique up to isotopy.
\end{prop}
\begin{proof}
The existence is proved above. The uniqueness up to isotopy is 
proved as follows:
Each regular neighborhood fibers over a $1$-orbifold 
with fibers connected $1$-orbifolds in the orbifold sense. 
A regular neighborhood can be isotoped into 
any other regular neighborhood by contracting 
in the fiber directions. We can modify the proof of 
Theorem 5.3 of Chapter 4 of \cite{Hirsch} to be adopted 
to an annulus with a finite group acting on it and an imbedded 
circle. 
\end{proof}

\subsection{Splitting and sewing on $2$-orbifolds reinterpreted}
An orbifold is said to be an {\em open orbifold} if the underlying
space is noncompact and the boundary is empty.  If one removes $l$
from the interior of these orbifolds, we obtain either a union of one
or two open annuli, or a union of one or two open silvered strip.  In
(2)-(4), an open annulus results.  For (1), a union of two open annuli
results.  For (6)-(7), an open silvered strip results.  For (5), we
obtain a union of two open silvered strips.  These can be easily
completed to be a union of one or two compact annuli or a union of one
or two silvered strips respectively.  To see this simply identify them
to be dense open suborbifolds of the unions of compact annuli or
silvered strips.

Let $l$ be a $1$-orbifold embedded in the interior of an orbifold $S$.
We can complete $S -l$ in this manner: We take a closed regular
neighborhood $N$ of $l$ in $S$.  We remove $N-l$ to obtain the above
types and complete it and re-identify with $S-l$ to obtain a
compactified orbifold. This process is the splitting of $S$ along $l$.

Conversely, we can describe sewing: Take an open annular
$2$-orbifold $N$ which is a regular neighborhood of a $1$-orbifold
$l$.  Suppose that $l$ is a circle.  We obtain $U = N-l$ which is a
union of one or two annuli.  Take an orbifold $S'$ with a union $l'$
of one (resp. two) boundary components which are circles.  Take an
open regular neighborhood of $l'$ and remove $l'$ to obtain $V$. $U$
and $V$ are the same orbifold. We identify $S'-l'$ and $N -l$ along
$U$ and $V$. This gives us an orbifold $S$, and it is easy to see that
$S$ is obtained from $S'$ by sewing along $l'$.  $l$ corresponds to a
$1$-orbifold $l''$ in $S$ in a one-to-one manner.  We can obtain
(1),(2),(3)-type neighborhoods of $l''$ in this way.
The operation in case (1) is said to be {\em pasting}, 
in case (2) {\em cross-capping}, and in case (3) {\em silvering}
along simple closed curves.

Suppose that $l$ is a full $1$-orbifold.  $U= N -l$ is either an open
annulus or a union of one (resp. two) silvered strips.  The former
happens if $N$ is of type (4) and the latter if $N$ is of type
(5)-(7).  In case (4), take an orbifold $S'$ with a boundary component
$l'$ a circle. Then we can identify $U$ with a regular neighborhood of
$l'$ removed with $l'$ to obtain an orbifold $S$.  Then $l$
corresponds a full $1$-orbifold $l''$ in $S$ in a one-to-one
manner. $l''$ has a type-(4) regular neighborhood.
The operation is said to be {\em folding} along a simple closed 
curve.

In the remaining
cases, take an orbifold $S'$ with a union $l'$ of one (resp. two)
boundary full $1$-orbifolds. Take a regular neighborhood $N$ of $l'$
and remove them to obtain $V$. Identify $U$ with $V$ for $S'-l'$ and
$N-l$ to obtain $S$.  Then $S$ is obtained from $S'$ by sewing along
$l'$.  Again $l$ corresponds to a full $1$-orbifold $l''$ in $S$ in a
one-to-one manner.  We obtain (5),(6), and 
(7)-type neighborhoods of $l''$ in this way, where
the operations are said to be {\em pasting}, {\em folding}, 
and {\em silvering} along full $1$-orbifolds respectively.

In other words, silvering is the operation of removing a regular 
neighborhood and replacing by a silvered annulus or 
a half square. Clarifying is an operation of 
removing the regular neighborhood and replacing 
an annulus or a silvered strip.

\begin{prop}\label{prop:eulerspli}
The Euler characteristic of an orbifold 
before and after splitting or sewing remains unchanged.
\end{prop}
\begin{proof}
Form regular neighborhoods of the involved 
boundary components of the split orbifold and those of 
the original orbifold. They have zero Euler characteristic. 
Since their boundary $1$-orbifolds have zero Euler characteristic,
the lemma follows by the additivity formula \eqref{eqn:additivity}.
\end{proof}

\subsection{Identification interpretations of 
splitting and sewing}
\label{subsec:topop}

In the following we describe the topological identification process of
the underlying space involved in these six types of sewings. The
orbifold structure on the sewed orbifold should be clear.

Let an orbifold $\Sigma$ have a boundary component $b$.  ($\Sigma$ is
not necessarily connected.)  $b$ is either a simple closed curve or a
full $1$-orbifold.  We will now find a $2$-orbifold $\Sigma''$
constructed from $\Sigma$ sewed along $b$ and/or some other component
of $\Sigma$.

We need to look at the cases when $b$ corresponds to a simple closed
curve or $1$-orbifold in $\Sigma''$.  The underlying space
$X_{\Sigma''}$ of $\Sigma''$ is a surface with or without
boundary. First, we need to introduce some identification of $b$ so
that $X_{\Sigma''}$ is a surface. Second, we need to find a suitable
orbifold structure on $X_{\Sigma''}$ so that the splitting will give
us back $\Sigma$:

%\subsection{A}
(A) Suppose that $b$ is diffeomorphic to a circle; that is,
$b$ is a closed curve.
Let $\Sigma'$ be a component of the $2$-orbifold $\Sigma$ 
with boundary component $b'$. 
Note that $\Sigma'$ may be the same component as the component
where $b$ lies in, and $b'$ may equal $b$. 
Suppose that there is a diffeomorphism $f: b \ra b'$. Then 
we obtain a bigger orbifold $\Sigma''$ glued along $b$ and $b'$
topologically. 
\begin{itemize}
\item[(I)] The construction so that $\Sigma''$ does
not create any more singular point results in an orbifold $\Sigma''$ so
that \[\Sigma'' - (\Sigma - b \cup b')\]
is a circle with neighborhood either diffeomorphic
to an annulus or a M\"obius band. 
\begin{itemize}
\item[(1)] In the first case, $b \ne b'$ (pasting).
\item[(2)] In the second case, $b = b'$ and 
$\langle f\rangle $ is of order two without fixed points 
(cross-capping).
\end{itemize}
\item[(II)] When $b = b'$,
the construction so that $\Sigma''$ does 
introduce more singular points to occur
in an orbifold $\Sigma''$ so that 
\[\Sigma'' - (\Sigma - b)\] 
is a circle of mirror points or is a full $1$-orbifold
with endpoints in cone-points of order two
depending on whether $f:b \ra b$ is 
\begin{itemize}
\item[(1)] the identity map (silvering), or
\item[(2)] is of order two and 
has exactly two fixed points (folding).
\end{itemize}
\end{itemize}
%\subsection{B}
(B) Consider when $b$ is a full $1$-orbifold
with endpoints mirror points.
\begin{itemize}
\item[(I)] Let $\Sigma'$ be a component orbifold (possibly the same 
as one containing $b$) with boundary full $1$-orbifold $b'$
with endpoints mirror points where $b \ne b'$.
We obtain a bigger orbifold $\Sigma''$ by gluing 
$b$ and $b'$ by a diffeomorphism $f: b \ra b'$. 
This does not create new singular points (pasting). 
\item[(II)] Suppose that $b = b'$.
Let $f: b\ra b$ be the attaching map. Then 
\begin{itemize}
\item[(1)] if $f$ is the identity, then we silvered $b$
and morphed the end points into corner-reflectors of order two
(silvering).
%(See Remark \ref{rem:freeorb}.) 
%C. delete
\item[(2)] If $f$ is of order two, then $\Sigma''$ has 
a new cone-point of order two and have one-boundary
orbifold reduced away. $b$ corresponds to 
a mixed type $1$-orbifold in $\Sigma'$ (folding).
\end{itemize}
\end{itemize}
It is obvious how to put the orbifold structure 
on $\Sigma''$ using the previous descriptions using
regular neighborhoods above. 

We shall use these labellings to describe the topological
operations.

\section{Projective orbifolds and the Hitchin-Teichm\"uller components}
\label{sec:projorb}

%In this paper, we will let $X= \rpt$ and $G =\PGL(3, \bR)$. 
%C. delete
By an {\em $\rpt$-structure} or {\em projectively flat structure}
on a $2$-orbifold $\Sigma$ we mean an $(\rpt, \PGL(3, \bR))$-structure
on $\Sigma$. From now on, we look at two-dimensional 
$\rpt$-orbifolds, that is, orbifolds with $\rpt$-structures.

In this section, we will define the deformation spaces of
$\rpt$-structures on orbifolds, describe local properties, and define
convex $\rpt$-structures (when the orbifolds are boundaryless).

We will discuss the relationship between the $\rpt$-structures and
holonomy representations.  First, we deduce that the deformation
space is locally Hausdorff from the corresponding property
of the holonomy representation variety.
Next, we discuss convex $\rpt$-structures.  We show that the
deformation space of convex $\rpt$-structures on an orbifold is an
open subset of the full deformation space.  We identify the
deformation space of convex $\rpt$-structures on orbifolds with 
a subset of the space of conjugacy classes of representations of its
fundamental group using the above relationship.

% We will prove facts needed to prove Theorem B in this section 
% but give a proof of Theorem B in  Section \ref{sec:thmB} since 
% we will need Theorem A 
% to complete the proof. However, the facts we present 
% in this section is needed for Theorem A. Therefore,
% we are following logical progression.
%%% UNTANGLE THIS!
%C. Maybe say it better like:
%The results of this section is needed to prove Theorem A
%and Theorem B for which we need Theorem A.

\subsection{Developing orbifolds}

Thurston shows that all orbifolds admitting $(X, G)$-structures are
good \cite{Thnote}.  It also follows from his work that for a
$2$-orbifold $\Sigma$, the existence of $(\rpt, \PGL(3,
\bR))$-structure is equivalent to giving a developing map $\dev:\tilde
\Sigma\ra \rpt$ from the universal covering space $\tilde \Sigma$
of $\Sigma$ equivariant with respect to the holonomy homomorphism
$h:\pi(\Sigma) \ra G$ where $\pi(\Sigma)$ is the deck transformation
group of $\tilde \Sigma$.  In other words $h$ satisfies:
\begin{equation}\label{eqn:equiva}
h(\vth)\circ \dev = \dev \circ \vth, \vth \in \pi(\Sigma).
\end{equation}
(See Proposition 5.4.2 of Thurston \cite{Thnote}.)

Two pairs $(\dev, h)$ and $(\dev', h')$ are {\em equivalent\/} if and
only if there exists $\psi \in G$ satisfying
\begin{equation}\label{eqn:conj}
\dev' = \psi\circ \dev, h'(\cdot) = \psi \circ h(\cdot) \circ \psi^{-1} 
\mbox{ for some } \psi \in G.
\end{equation}
The development pair is determined by the structure and a germ of a
structure at the basepoint; changing the germ of the structure
changes the pair in its equivalence class. Thus the {\em equivalence class\/}
of the development pair is uniquely determined by the structure.
%%%%!
%\marginpar{** }
%C. don't understand.

The image of $\dev$ is a {\em developing image}
and that of $h$ is a {\em holonomy group}.

%\subsection{Reflections and rotations}

\subsection{Types of Singularities}
An automorphism of $\rpt$ is said to be a {\em reflection} if
its matrix is conjugate to 
\[\begin{bmatrix} 1 & 0 & 0 \\ 0 & 1 & 0 \\ 0 & 0 & -1 \end{bmatrix}.\]
A reflection has a line of fixed points and an isolated fixed point,
which is said to be the {\em reflection point}. 
An automorphism of $\rpt$ is said to be a {\em rotation of order $n$}, 
$n = 2, 3, \dots $, if its matrix is conjugate to 
\[\begin{bmatrix} \cos 2\pi/n & -\sin 2\pi/n & 0 \\
\sin 2\pi/n & \cos 2\pi/n & 0 \\
0 & 0 & 1 \end{bmatrix}. \]

A rotation has a unique isolated fixed point, called a {\em rotation point}, 
and an invariant line. %%%%! where the rotation acts like a ``rotation''. 
A one-parameter family of invariant ellipses fills the complement
in $\rpt$ of the rotation point and the invariant line.
%%%%! maybe should define elliptic element
%C. elliptic? How to define?
A rotation of order two is a reflection also
and conversely.

For $\rpt$-orbifolds, the singular points have 
the neighborhoods with model open sets and finite group
actions corresponding to one of the following:
\begin{itemize}
\item[(i)] A mirror point: An open disk in $\rpt$ meeting a line of
fixed points of a reflection.
\item[(ii)] A cone-point of order $n$: An open disk in 
$\rpt$ containing a rotation point of the rotation of 
order $n$. 
\item[(iii)] A corner-reflector of order $n$: An open disk 
in $\rpt$ containing the intersection point of the lines of
fixed points of two reflections $g_1$ and $g_2$ so that 
$g_1\circ g_2$ is a rotation of order $n$.
\end{itemize}
Actually, these models are all projectively isomorphic up to
choices of the open disks.

To give some examples, consider the quotient of $\bR^2$ by reflections
about horizontal and vertical lines through integer points. This gives
us a square with boundary mirror points and corner-reflectors of order two.

\subsection{Example:  Elementary annuli}
Let $\vth$ 
be a collineation represented by a diagonal matrix with
positive eigenvalues.  Then it has three fixed points in $\rpt$: an
attracting fixed point of the action of $\langle \vth \rangle$, a
repelling fixed point, and a saddle type fixed point. Three lines
passing through two of them are $\vth$-invariant, as are four open
triangles bounded by them. Choosing two open sides of an open triangle
ending at an attracting fixed point or a repelling fixed point
simultaneously, their union is acted properly and freely upon by
$\langle \vth \rangle$.  The quotient space is diffeomorphic to an
annulus.  The $\rpt$-surface projectively diffeomorphic to the
quotient space is said to be an {\em elementary annulus}.

\subsection{Example: $\pi$-Annuli}

Take two adjacent triangles, and three open sides of them all ending
in an attracting fixed point or a repelling fixed point. Then the
quotient of the union by $\langle \vth \rangle$ is diffeomorphic to an
annulus. The projectively diffeomorphic surfaces are said to be {\em
$\pi$-annuli}.  (See \cite{cdcr1} and \cite{cdcr2} for more details.)

A reflection sending one triangle to the other induces an order-two
group.  The quotient map is an orbifold map, and the quotient
space carries an orbifold structure so that one boundary component is
made of mirror points.  (See Page \pageref{page:silvered} for more details).

\begin{figure}[ht]
\centerline{\epsfxsize=3in \epsfbox{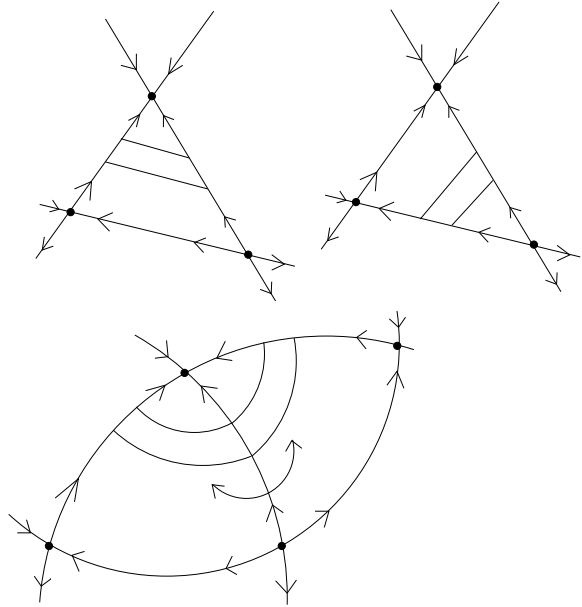}}
\caption{\label{fig:annuli} Elementary annuli and 
a $\pi$-annulus and an action on it.} 
\end{figure}
\typeout{<>}

\subsection{The deformation spaces and holonomy}
We now define the deformation space $\rpt(\Sigma)$ of
$\rpt$-structures on a connected $2$-orbifold $\Sigma$ as follows
(assuming $\Sigma$ is connected and has empty boundary): Give
the $C^1$-topology to the set ${\mathcal{S}}(\Sigma)$ of all pairs
$(\dev, h)$ satisfying equation \eqref{eqn:equiva} on $\tilde \Sigma$.
Two pairs $(\dev, h)$ and $(\dev', h')$ are {\em equivalent under
isotopy} if there exists a self-diffeomorphism $f$ of the universal
cover $\tilde \Sigma$ of $\Sigma$ commuting with the deck
transformations so that $\dev' =\dev \circ f$ and $h' = h$.
We denote by ${\rpts}(\Sigma)$ the space of
equivalence classes with the quotient topology.

The pairs $(\dev, h)$ and $(\dev', h')$ are equivalent 
under $\PGL(3, \bR)$-action, if there exists an element 
$g$ of $\PGL(3, \bR)$ so that $\dev' = g \circ \dev$ and 
$h'(\cdot) = g h(\cdot) g^{-1}$.
The quotient space of ${\rpts}(\Sigma)$ under the 
$\PGL(3, \bR)$-equivalence relation
is denote by $\rpt(\Sigma)$. 

%An {\em isotopy} of  $2$-orbifold $\Sigma$ is a self-diffeomorphism $f$ 
%that so that there exists an orbifold map 
%$F:\Sigma \times I \ra \Sigma$ so that 
%$F_t:\Sigma \ra \Sigma$ given by $F_t(x) = F(x, t)$ 
%is a diffeomorphism for each $t$ and 
%$F_0$ is the identity and $F_1 = f$. 

Another interpretation of the deformation space is 
to consider all $\rpt$-structures on $\Sigma$ and 
quotient by the isotopies. 
One can easily obtain that there is a one-to-one 
correspondence between the above two spaces. 

If two $\rpt$-structures are distinct up to isotopy, we say that
they are {\em isotopically distinct}.  
Isotopically distinct $\rpt$-structures 
represent different points in the deformation spaces.

Two $\rpt$-orbifolds with non-conjugate holonomy 
homomorphisms cannot be isotopic (see \cite{dgorb} for details).

By forgetting $\dev$ from the pair $(\dev, h)$, we obtain an induced map
\[\mathcal{H}': {\rpts}(\Sigma) \ra 
   \mbox{Hom}(\pi_1(\Sigma), \PGL(3, \bR))\]
to the space of homomorphisms of $\pi_1(\Sigma)$
since the isotopy does not change the holonomy homomorphism.

By Theorem \ref{thm:finitegen}, we see that $\pi_1(\Sigma)$
is a finitely presented group.
From now on, we denote
\[H(\Sigma) = \mbox{Hom}(\pi_1(\Sigma), \PGL(3, \bR))\]
for the $\bR$-algebraic subset
of $\PGL(3, \bR)^n$ where $n$ is the number of the generators of 
$\pi_1(\Sigma)$. 
% (To see this, represent the fundamental group
% by group presentations, and generators correspond to each factor of
% $\PGL(3, \bR)$, and the group relations give us the defining 
% equations). 
% C. I am not so sure to delete.

The main result of \cite{dgorb} is that 
the map $\mathcal{H}'$ is 
a local homeomorphism since $\pi_1(\Sigma)$ is finitely presented.
The proof of this is not much different
from the manifold case.
The idea of proof is based on \cite{Grep:88} and 
Morgan and Lok \cite{Lok:84} (from lectures of Morgan)
generalizing Weil's work \cite{We:60}.

Let $U^n$ denote the open subset of $\PGL(3, \bR)^n$ consisting
of $(X_1, \dots, X_n)$ such that no line in $\bR^3$ is simultaneously
invariant under $X_1, \dots, X_n$ where the 
$\PGL(3, \bR)$-action is proper and free (see \cite{Gconv:90}).
Let $U(\Sigma)$ denote $H(\Sigma) \cap U^n$. 

\begin{thm}\label{thm:Haus} 
Let $\Sigma$ be a connected closed $2$-orbifold with 
$\chi(\Sigma)<0$.
Then $\rpt(\Sigma)$ has the structure of
Hausdorff real analytic variety modeled on $U(\Sigma)/\PGL(3, \bR)$,
and the induced map 
\[\mathcal{H}: \rpt(\Sigma) \ra U(\Sigma)/\PGL(3, \bR)\] 
is a local homeomorphism.
\end{thm}
\begin{proof}
First, we show that the image of $\mathcal{H}'$ is in $U(\Sigma)$.
By Lemma 2.5 of \cite{Gconv:90},the holonomy group 
a $2$-orbifold of negative Euler characteristic fixes no line,
since a finite-index subgroup is a fundamental
group of a closed surface of negative Euler characteristic.
Since $\pi_1(\Sigma)$ is finitely presented by 
Theorem \ref{thm:finitegen}, the holonomy map is 
a local homeomorphism follows from \cite{dgorb}. 
\end{proof}

Suppose now that $\Sigma$ has more than one components.
Let $\Sigma_1, \dots, \Sigma_n$
denote the connected components of $\Sigma$. 
Define $\rpt(\Sigma)$ to be the product
\begin{equation*}
\rpt(\Sigma_1) \times \dots \times \rpt(\Sigma_n).
\end{equation*}
Similarly define ${\rpts}(\Sigma)$ and $\mathcal{S}(\Sigma)$ 
and $H(\Sigma)$. 
If the Euler characteristic of each component of $\Sigma$ is negative, 
the product map
\[\mathcal{H}': {\rpts}(\Sigma) \ra 
H(\Sigma) := \prod_{i=1}^n H(\Sigma_i) \] 
is a local homeomorphism, and so
is the product map
\[\mathcal{H}: \rpt(\Sigma) \ra 
U(\Sigma)/\PGL(3, \bR) :=
\prod_{i=1}^n (U(\Sigma_i)/\PGL(3, \bR)).\]

\subsection{Hitchin-Teichm\"uller components}
Let $\Sigma$ be a closed connected $2$-orbifold with $\chi(\Sigma)<0$.
An $\rpt$-structure on $\Sigma$ is {\em convex} 
if $\dev$ is a diffeomorphism to 
a convex subset of an affine patch in $\rpt$.
Then $\Sigma$ is projectively diffeomorphic to
$\Omega/\Gamma$ where $\Gamma$ is a discrete subgroup of 
$\PGL(3, \bR)$ acting on a convex domain $\Omega$. 
By Theorem 3.2 of \cite{Gconv:90},
$\Omega$ is strictly convex in an affine patch of
$\rpt$, which is also precompact. $\partial \Omega$ is $C^1$
and contains no line segment. 

A conic in $\rpt$ is given by a zero set of 
a quadratic form of signature $(2, 1)$.
in $\rpt$ with homogeneous coordinates $x_0, x_1, x_2$. 
The group acting on a conic is conjugate to the projectivized
linear Lorentz group $\PO(1, 2)$.
By applying a collineation we may assume that the conic is defined
by the standard diagonal quadratic form	of signature $(2,1)$ and
that the holonomy homomorphism maps to $\PO(1, 2)$.

Hyperbolic structures on $\Sigma$ form a distinguished class of convex
$\rpt$-structures:  The Klein model of hyperbolic geometry identifies
a hyperbolic plane and the group of isometry with the convex domain
$\Omega$ bounded by a conic and its group of collineations
$\PO(1, 2)$.  A hyperbolic structure on a $2$-orbifold has a chart
into $\Omega$ with transition functions in $\PO(1, 2)$.
Since $\Omega$ is a subset of $\rpt$ and $\PO(1, 2)$ is a subgroup
of $\PGL(3, \bR)$, such a $2$-orbifold has a $\rpt$-structure, which
is said to be a {\em hyperbolic $\rpt$-structure}.  
If a $2$-orbifold has a hyperbolic structure with geodesic boundary, 
then it has an $\rpt$-structure with geodesic boundary. The
orbifold is said to have a {\em hyperbolic $\rpt$-structure}.

The subspace of $\rpt(\Sigma)$ of elements represented by convex
$\rpt$-structures will be denoted by $\mathcal{C}(\Sigma)$.  The
subspace of $\mathcal{C}(\Sigma)$ corresponding to
hyperbolic $\rpt$-structures is denoted by
$\mathcal{T}(\Sigma)$ and identifies with 
the Teichm\"uller space of hyperbolic structures on $\Sigma$ 
as determined by Thurston \cite{Thnote}. 
This follows since a projective diffeomorphism of two
hyperbolic $\rpt$-orbifolds is obviously an isometry of them.
This is also a topology preserving identification 
since the both topologies of the deformation spaces and 
the Teichm\"uller spaces are defined by $C^1$-topology 
of developing maps.

For later purposes, we define $\mathcal{C}'(\Sigma)$
be the subset of ${\rpts}(\Sigma)$ consisting of  isotopy
classes of convex structures on $\Sigma$.  

The {\em pre-Hitchin-Teichm\"uller component} of $H(\Sigma)$ is a component
$C_{\mathcal{T}}$ of it which contains representations
\[\Hom(\pi_1(\Sigma), \PO(1, 2)\] corresponding
to holonomy homomorphisms of hyperbolic structures on $\Sigma$
(discrete embeddings $\pi_1(\Sigma)\longrightarrow\PO(1,2)$).

The group $\PGL(3, \bR)$ acts on $H(\Sigma)$ by conjugation, that is,
\begin{equation*}
h(\cdot) \mapsto \vth h(\cdot) \vth^{-1}, \vth \in \PGL(3, \bR)  
\end{equation*}
for $h\in H(\Sigma)$.  Let $H(\Sigma)^{st}$ be the subspace
of representations $r$ acting freely on $\rpt$. By Lemma
1.12 of \cite{Gconv:90}, $\PGL(3, \bR)$ acts properly on this subset.

%\begin{prop}\label{prop:component}
%The pre-Hitchin-Teichm\"uller component $C_{\mathcal{T}}$ 
%is an open and closed subset of the stable set above, 
%and $\PGL(3,\bR)$ acts properly on $C_{\mathcal{T}}$ properly. 
%\end{prop}
%The proof is the first part of the proof of Theorem B in subsection
%\ref{subsec:thmB}.

%%4/25 9:30
We call a component $C_{\mathcal{T}}/\PGL(3, \bR)$ of
\[H(\Sigma)^{st}/\PGL(3, \bR)\]
a {\em Hitchin-Teichm\"uller component} of $\Sigma$ following 
Hitchin~\cite{Hit:92}. 
In the proof of the following theorem, we show that 
$C_{\mathcal{T}}$ may be defined to be a component 
of $H(\Sigma)$ as we defined in the introduction.
Theorem B states that this component is identical 
with the deformation space $\mathcal{C}(\Sigma)$ of
convex $\rpt$-structures on $\Sigma$. 

\subsection{Openness of convex $\rpt$-structures}
We will need Propositions \ref{prop:openness}
and \ref{prop:closedness} to prove Theorem B. 
They will be modified to Propositions \ref{prop:openness2} 
and \ref{prop:closedness2} in the next section.

\begin{prop}\label{prop:openness}
Let $S$ be a closed orbifold with $\chi(S)<0$.
Then $\mathcal{C}(S)$ is an open subset of $\rpt(S)$.
{\rm (}So is $\mathcal{C}'(S)$ of $\rpts(S)$.{\rm )} 
\end{prop}
\begin{proof}
It is shown in \cite{Gconv:90} using J. L. Koszul's result \cite{Kz:68}
that $\mathcal{C}(S')$ is an open subset of $\rpt(S')$
if $S'$ is a surfaces. This is done by taking product
$S' \times \SI^1$ and finding an affine structure on it
corresponding to an $\rpt$-structure on $S'$.
An affine structure on $S' \times \SI^1$ is convex 
if and only if the $\rpt$-structure on $S'$ is convex. 
Since all discussions in \cite{Kz:68} apply to differentiable 
orbifolds as well, we can replace $S'$ by $S$. 
\end{proof}

\subsection{Closedness of $\rpt$-structures}

The following proposition shows that $2$-orbifolds
of negative Euler characteristic with isomorphic fundamental
groups are diffeomorphic. The harmonic map theory
of Schoen-Yau~\cite{SY:78} and the Nielsen realization theorem proved by 
Kerckhoff~\cite{Kerck2} are essential for the proof:

\begin{prop}\label{prop:SY}
Let $\Sigma_1$ and $\Sigma_2$ be closed $2$-orbifolds 
of negative Euler characteristic 
where $\tilde \Sigma_1$ and $\tilde \Sigma_2$ are homeomorphic to disks. 
If $k: \pi_1(\Sigma_1) \ra \pi_1(\Sigma_2)$ is an isomorphism,
then there is an orbifold diffeomorphism 
$f:\Sigma_1 \ra \Sigma_2$ so that $\tilde f_* = k$
for a lift $\tilde f:\tilde \Sigma_1 \ra \tilde \Sigma_2$. 
\end{prop} 
\begin{proof} 
Since $\Sigma_1$ admits a hyperbolic structure, 
$\pi_1(\Sigma_1)$ is isomorphic to a discrete cocompact subgroup
of $\PSL(2, \bR)$. There is a torsion-free finite-index normal subgroup 
$\Gamma$ of $\pi_1(\Sigma)$ by Selberg's lemma~\cite{RT:94}. 
Let $\Gamma'$ be $k(\Gamma)$ in $\pi_1(\Sigma_2)$. 
There is a finite covering surface $\Sigma'_1$ of $\Sigma_1$
corresponding to $\Gamma$ and $\Sigma'_2$ of $\Sigma_2$
corresponding to $\Gamma'$. 

The finite group $G_1 =\pi_1(\Sigma_1)/\Gamma$ 
maps injectively into $\Out(\Gamma)= \Aut(\Gamma)/\Inn(\Gamma)$.
Similarly, $G_2 = \pi_1(\Sigma_2)/\Gamma'$ maps 
into $\Out(\Gamma')/\Inn(\Gamma')$.
Thus, the commutative diagram 
\begin{eqnarray}
G_1  & \stackrel{k}{\longrightarrow}  & G_2  \nonumber \\
\downarrow & & \quad \downarrow \nonumber \\
\Out(\Gamma)/\Inn(\Gamma) 
&\stackrel{k_*}{\longrightarrow} & \Out(\Gamma')/\Inn(\Gamma')
\end{eqnarray}
holds where $k_*$ is an induced isomorphism. 

By the Nielsen realization theorem of Kerckhoff \cite{Kerck2},
$G_1$ acts on $\Sigma_1$ and $G_2$ on $\Sigma_2$.
A homeomorphism $f':\Sigma'_1 \ra \Sigma'_2$ realizes $k|\Gamma$ 
and for each $g \in G_1$, $f'\circ g$  
is homotopic to $k(g)\circ f'$ for $k(g) \in G_2$.

Give $\Sigma'_1$ and $\Sigma'_2$ arbitrary hyperbolic
metrics which are $G_1$- and $G_2$-invariant respectively
(that is, using Thurston's orbifold hyperbolization of 
$\Sigma_1$ and $\Sigma_2$).
Then choose a unique harmonic diffeomorphism
$\hat f:\Sigma'_1\ra \Sigma'_2$ in the homotopy class of $f'$
as obtained by Schoen-Yau \cite{SY:78}.
Since $\hat f\circ g = k(g)\circ \hat f$ by uniqueness,
$\hat f$ induces a desired orbifold diffeomorphism 
$f:\Sigma_1 \ra \Sigma_2$.
\end{proof}

\begin{prop}\label{prop:closedness}
Let $S$ be a closed orbifold with $\chi(S)<0$.  
The image of
\[\mathcal{H}':{\mathcal{C}'}(S) \ra H(S)\] 
is closed. 
\end{prop}
\begin{proof}
Choose a sequence of representations in $\mathcal{H}'$
\[h_i:\pi_1(S) \ra \PGL(3, \bR)\]
so that 
\[(h_i(g_1), \dots, h_i(g_m)) \ra (h(g_1), \dots, h(g_m))\]
for a representation $h:\pi_1(S) \ra \PGL(3, \bR)$;
that is, $h_i$ converges to $h$ algebraically.
We show that $h$ is in the image proving the closedness. 

Let $S_i$ be the orbifold $S$ with an $\rpt$-structure 
corresponding to $h_i$, and let $\tilde S$ be the universal cover of $S$.
Let $\tilde S_i$ be the universal cover $\tilde S$ with 
the induced $\rpt$-structure from $S_i$, and 
$D_i$ a developing map of $S_i$ associated with $h_i$.
Then $D_i:\tilde S_i \ra \rpt$ maps 
onto a strictly convex domain $\Omega_i$
in an affine patch of $\rpt$. 
(This follows since 
$\tilde S_i$ is a universal cover of a convex 
$\rpt$-surface finitely covering $S_i$.
See \cite{Gconv:90}.)
$D_i:\tilde S_i \ra \Omega_i$
induces an $\rpt$-diffeomorphism
$S_i \ra \Omega_i/h_i(\pi_1(S))$.

There is a unique $\rpt$-structure on the sphere $\SI^2$
such that the covering projection 
$p_{\rpt}:\SI^2 \ra \rpt$ is a projective map.
The nontrivial deck transformation is represented by the antipodal
map of the sphere.
Its collineation group $\Aut(\SI^2)$ 
is isomorphic to $\SL_\pm(3, \bR)$, generated by $\PGL(3,\bR)$ and
the antipodal map,

We can show that $D_i:\tilde S_i \ra \rpt$ always
lifts to an embedding $D'_i:\tilde S_i \ra \SI^2$ 
and $h_i$ lifts to a homomorphism 
$h'_i:\pi_1(S) \ra \Aut(\SI^2)$
(see \cite{cdcr1}): we can lift first,
and for a deck-transformation $\vth$ of
$\pi_1(S)$, $D'_i\circ \vth$ is another developing map,
and hence it must equal $\vpi\circ D'_i$ for $\vpi \in \Aut(\SI^2)$. 
Defining $h'_i(\vth)= \vpi$, we see that $h'_i$ is a lift of $h_i$.

The image $\Omega'_i$ of $D'_i$ is a convex open subset of an open
hemisphere in $\SI^2$ with a standard geodesic structure.

By choosing a subsequence if necessary, the sequence of the closures
$\clo(\Omega'_i)$ converges to a compact convex subset of $\SI^2$ in a
closed hemisphere (Choi-Goldman \cite{CG:93}). We claim that the limit
$\Omega'_\infty$ is neither a point, a line segment, a lune, nor a closed
hemisphere.  Otherwise, (taking a finite subcover $S'$ of $S$ if necessary), 
all $D'_i(\tilde S_i)$ are images of a sequence of developing
images of convex $\rpt$-structures on a closed surface $S'$. We
showed that such a degeneration cannot happen in \cite{CG:93}.  (See
also \cite{mar:96}.)  Thus, $\Omega'_\infty$ is a compact convex
subset of an open hemisphere in $\SI^2$.

By choosing a subsequence, $h'_i$ converges to a
representation $h':\pi_1(S) \ra \Aut(\SI^2)$ lifting $h$. 
As in \cite{CG:97}, $h'(\pi_1(S))$ acts on $\Omega'_\infty$. 
Since $h'$ is a map to $\SL_\pm(3, \bR)$, $h'$ is discrete and
faithful by Lemma 1.1 of Goldman-Millson \cite{GM:88}. ($\pi_1(S)$ has
a finite index subgroup which is torsion-free.  Apply Lemma 1.1 of
\cite{GM:88} here and the finite index extension argument is trivial.)

Therefore, $h'(\pi_1(S))$ acts on an open disk 
$\Omega^{\prime o}_\infty$ with quotient orbifold $S''$. 
By Proposition \ref{prop:SY}, there is a
diffeomorphism $S''\longrightarrow S$ inducing $h'$. 
Since $p_{\rpt}|\Omega^{\prime o}_\infty$ 
is an embedding onto an 
$h(\pi_1(S))$-invariant convex open domain $\Omega\subset\rpt$,
we see that $S''$ is realized also as 
the quotient space $\Omega/h(\pi_1(S))$.
Thus, $h$ is realized as a holonomy homomorphism
of a convex $\rpt$-structure on $S$ and 
lies in the image of 
$\mathcal{H}'(\mathcal{C}(S))$.
\end{proof}

\subsection{The proof of Theorem B}\label{subsec:thmB}

For the proof, we need Theorem A, which will be proved in \S 6.

%4/26 11:45
\begin{proof}[Proof of Theorem B]
Let $\Sigma$ be a closed $2$-orbifold with $\chi(\Sigma)<0$.
Let $\mathcal{T}'(\Sigma)$ the connected subset of $\mathcal{C}'(\Sigma)$
consisting of hyperbolic $\rpt$-structures on $\Sigma$.
Since by Theorem A, $\mathcal{C}'(\Sigma)$ is connected,
and $\mathcal{H}'$ sends $\mathcal{T}'(\Sigma)$
into $C_{\mathcal{T}}$, it follows that 
$\mathcal{H'}$ sends $\mathcal{C}'(\Sigma)$ into 
$C_{\mathcal{T}}$.

By Proposition \ref{prop:openness}, 
$\mathcal{C}'(\Sigma)$ is an open subset of 
$\rpts(\Sigma)$. Since $\mathcal{H}'$ is an open map, 
$\mathcal{H}'(\mathcal{C}'(\Sigma))$ is an open subset of $C_{\mathcal{T}}$. 
By Proposition \ref{prop:closedness}, the image 
is a closed subset of $C_{\mathcal{T}}$. 
Hence, the image equals $C_{\mathcal{T}}$. 

The holonomy group of a convex $\rpt$-orbifold 
is discrete since it acts on an open domain
discontinuously. Thus, $C_{\mathcal{T}}$ consists of discrete 
embeddings.

Recall also that since $\PGL(3, \bR)$ acts properly on $U(\Sigma)$,
$C_{\mathcal{T}}$ is a subset of $H(\Sigma)^{st}$. 
%proving Proposition \ref{prop:component} as promised.

To complete the proof of Theorem~B, we show that 
$\mathcal{H}'|\mathcal{C'}$ is injective: That is, given
two holonomy representations $h$ and $h'$ for 
$(D, \tilde f:\tilde \Sigma \ra \tilde M)$ and 
$(D', \tilde f':\tilde \Sigma \ra \tilde M')$ 
for convex  $\rpt$-orbifolds $M$ and $M'$.  Then, 
\begin{equation*}
h = h_1\circ \tilde f_*, \qquad h' = h'_1\circ \tilde f'_* 
\end{equation*}
for holonomy
homomorphisms of $h_1$ and $h'_1$ of $M$ and $M'$ respectively. We
show that if $h=h'$, then $(D, \tilde f)$ and $(D', \tilde f')$ are 
isotopic equivariantly with respect to $h' \circ h^{-1}$.
(According to the definition of deformation spaces in \cite{dgorb}, 
this will prove the injectivity of $\mathcal{H}'$.)

Let $\Sigma'$ be a closed surface finitely covering $\Sigma$.  Let
$\Omega$ be the image of $D$ composed with $\tilde f$ and $\Omega'$
that of $D'$ composed with $\tilde f'$.  Since $h$ and $h'$ restricted
to $\pi_1(\Sigma')$ are the same, Proposition 3.4 of \cite{Gconv:90}
shows $\Omega = \Omega'$.  The images of $D\circ \tilde f$ and
$D'\circ\tilde f'$ are the same, and they are both equivariant under
the homomorphism $h=h':\pi_1(\Sigma) \ra \PGL(3, \bR)$.  The map 
$g= D^{\prime,-1}\circ D: \tilde M \ra \tilde M' $ is so that 
$D'\circ g = D$ and is equivariant under the homomorphism
\[g_*=i_{M, M'}=\tilde f'_* \circ \tilde f^{-1}_*: \pi_1(M) \ra \pi_1(M').\]

Let $\Sigma$ have a hyperbolic metric $\mu$. Then 
$\tilde f$ and $\tilde f'$ induce metrics on 
$\tilde M$ and $\tilde M'$ respectively. 
We now show that $g \circ \tilde f$ and $\tilde f'$ 
are $i_{M, M'}$-equivariantly isotopic. 
$\tilde M'$ has induced 
hyperbolic metrics $\mu_0$ and $\mu_1$
induced from $g\circ \tilde f$ and $\tilde f'$ respectively.
There is a path of Riemannian metrics 
\begin{equation*}
\mu_t = t\mu_1 + (1-t)\mu_0  
\end{equation*}
for $t\in [0,1]$ from $\mu_0$ to 
$\mu_1$. By the equivariance, they induce metrics on $M'$ to be denoted
by same letters.
Recall that $\Sigma'$ is the closed surface covering $\Sigma$.
Let $M'_s$ denote the corresponding closed surface 
covering $M'$.
Let $\mu'_t$ denote the Riemannian metrics of $M'_s$ 
corresponding to $\mu_t$. 
By Theorem B.26 of Tromba \cite{Trom:92}, 
there exists a smooth one-parameter family of harmonic diffeomorphisms 
$S'(\mu_t): (M'_s, \mu'_t) \ra (\Sigma', \mu)$.
Since these harmonic diffeomorphisms are unique in
their homotopy classes, they should be equivariant 
under automorphisms of $M'_s$ and $\Sigma'$, and
$S'(\mu_t)$ descend to
a one-parameter family of diffeomorphisms 
$S(\mu_t):  (M', \mu_t) \ra (\Sigma, \mu)$. 
One can lift the inverse map $S(\mu_t)^{-1}$ 
of $S(\mu_t)$ to a smooth one-parameter
family of diffeomorphisms
$\tilde S(\mu_t)^{-1}: \tilde \Sigma \ra \tilde M'$
using analytic continuations.
By uniqueness of harmonic diffeomorphisms,
we see that the inverse map $S(\mu_0)^{-1}$ lifts to
$g\circ \tilde f: \tilde \Sigma_1 \ra \tilde M'$,
and by analytic continuation $S(\mu_1)^{-1}$ lifts to
$\gamma \circ \tilde f':\tilde \Sigma \ra \tilde M'$
for some deck transformation $\gamma$ of $\tilde M'$, where
$g\circ \tilde f$ is isotopic to $\gamma \circ \tilde f$
equivariant with respect to $i_{M, M'}$.
$\tilde S(\mu_1)^{-1}_*(\cdot)$ must equal 
$\gamma \circ \tilde f'_*(\cdot) \circ \gamma$. 
Since $g_* \circ \tilde f_* = \tilde f'_*$,
Proposition 8 of \cite{dgorb} implies that $\gamma$ equals
the identity since the center of $\pi_1(\Sigma)$ is trivial.
Therefore $\tilde S(\mu_1)^{-1}$ 
equals $\tilde f'$. Thus 
$g\circ \tilde f$ and $\tilde f'$ are 
$i_{M, M'}$-equivariantly isotopic.

Applying $D'$ to $g\circ \tilde f$ and $\tilde f'$ 
and $D'_*$ to $i_{M, M'}$, implies that
$(D, \tilde f)$ and $(D',\tilde f')$ are equivalent. 
Therefore, $\mathcal{H}':\mathcal{C}'(\Sigma) \ra C_{\mathcal{T}}$ is 
a homeomorphism, inducing one
\[\mathcal{H}:\mathcal{C}(\Sigma) \ra C_{\mathcal{T}}/\PGL(3, \bR).\]
\end{proof}

% do you think any reader would care about this?
% \begin{rem}\label{rem:openclosed}
% We will need Propositions \ref{prop:openness2} and 
% \ref{prop:closedness2} modified from Propositions 
% \ref{prop:openness} and \ref{prop:closedness}
% to prove Theorem A. Theorem A was needed in the proof of
% Theorem B. 
% \end{rem}
% C. delete

\section{Splitting and sewing $\rpt$-orbifolds}
\label{sec:splitting}

We describe geometric operations on 
convex $\rpt$-orbifolds corresponding to the topological
operations in \S\ref{sec:orbifolds}.
(We recommend \cite{Gconv:90}, \cite{cdcr1},
\cite{cdcr2}, and \cite{cdcr3} for background knowledge of 
$\rpt$-structures on surfaces.)  
The boundary invariants of a convex $\rpt$-orbifold 
are used to build bigger
convex orbifolds from by constructions
along $1$-dimensional suborbifolds.
Finally, we discuss how the deformation space 
of $\rpt$-orbifolds with boundary relates to
the space of conjugacy classes of representations, following 
\S\ref{sec:projorb}.  We prove the openness 
and closedness of convex $\rpt$-structures, i.e., 
Propositions \ref{prop:openness2}
and \ref{prop:closedness2}, which we need later, 
generalizing Propositions \ref{prop:openness} and 
\ref{prop:closedness}. The geometric operations discussed here 
induce fibrations of orbifold deformation spaces.

% We remark that for the following discussions, the proofs when
% $\Sigma$ is a surface works just as well. 
% C. delete

\subsection{The deformation spaces of boundary closed curves}
Let $\Sigma$ be a compact
convex $\rpt$-orbifold with nonempty boundary.  Let $(\dev, h)$ be its
development pair and $\tilde \Sigma$ the universal cover.  Let $b$ be
a closed curve in $\Sigma$ and $\tilde b$ a lift to $\tilde \Sigma$,
which is an embedded arc. Let $\gamma$ be the corresponding deck
transformation. Then $h(\gamma)$ must be conjugate to one of the
following matrices: % since $\gamma$ may be considered 
% a root of an element of the
% fundamental group of a finite covering $\rpt$-surface:
% C. delete
\begin{equation}
\begin{bmatrix} 
\lambda_1 & 0 & 0\\
0 &\lambda_2 & 0\\
0 & 0 &\lambda_3 
\end{bmatrix} 
\end{equation}
where $|\lambda_1| \geq |\lambda_2| \geq |\lambda_3|$, and
$\lambda_1\lambda_2\lambda_3 = 1$
or 
\begin{equation}
\begin{bmatrix}
\lambda_1 & 1 & 0 \\
0 &\lambda_1 & 0 \\
0 & 0 & \lambda_2
\end{bmatrix}
\end{equation}
where $\lambda_1^2\lambda_1 = 1$. 
In the former case, if all eigenvalues are positive,
$b$ or the holonomy of $b$ is said
to be {\em hyperbolic}, and in the second case, 
if all eigenvalues are positive again
{\em quasi-hyperbolic}.
A curve with quasi-hyperbolic holonomy
is homotopic to the boundary. (See Theorem 3.2 of
\cite{Gconv:90} and Proposition 4.5 in \cite{cdcr2}.)

\begin{rem}\label{rem:qh} 
Let $\Sigma$ be an orbifold with $\chi(\Sigma)<0$ and a closed
geodesic boundary component $\gamma$. Then 
the holonomy of $\gamma$ is either 
hyperbolic or quasi-hyperbolic.
% A matrix is 
% {\em quasi-hyperbolic} if it is conjugate 
% to a nondiagonalizable matrix with two positive 
% eigenvalues. 
%%%% didn't we just say this a paragraph ago?
%C. delete
This follows from the analogous property for
closed surfaces; see \cite{cdcr1,cdcr2}.
\end{rem}

In the hyperbolic case, if an eigenvalue is negative,
we say that $b$ or the holonomy of $b$ is said to
be a {\em hyperbolic slide-reflection}.

\begin{rem} 
The conjugacy classes
of hyperbolic automorphisms are classified 
by two real numbers $\lambda = \lambda_3$,
and $\tau = \lambda_1 + \lambda_2$. They are in the space
\[{\mathcal R} =\{(\lambda, \tau)| 0< \lambda < 1, 
2/\sqrt{\lambda} < \tau < \lambda + \lambda^{-2}\}. \]
Also, for $A \in \SL(3, \bR)$, $A$ is hyperbolic
if and only if $A \in (\lambda, \tau)^{-1}({\mathcal{R}})$.
If $\tau = 1 + \lambda^{-1},$ then $\lambda_2 = 1$, 
and the hyperbolic element is called {\em purely hyperbolic}.
Since hyperbolic elements of $\PO(1,2)$ are purely hyperbolic,
the holonomy of a closed essential curve in a hyperbolic $\rpt$-orbifold 
is purely hyperbolic.
\end{rem}

A hyperbolic collineation of $\rpt$ has three or two fixed points, and
$\dev\circ \tilde b$ is a line connecting two of the fixed points. In
the first case, there are three fixed points, which are attracting,
repelling, or saddle type ones in the dynamics of infinite cyclic action of
the powers of the automorphism. If $\dev\circ \tilde b$ connects the
attractor to the repelling fixed point, then $b$ is said to be {\em
principal}.

If $\Sigma$ has no boundary and is convex, 
then any closed geodesic $b$ is principal
(see \cite{Gconv:90} and \cite{cdcr2}).
When $\Sigma$ has boundary diffeomorphic to a circle, 
we will require it to be a principal closed geodesic 
in this paper for convenience, and in this case, 
all closed curves are homotopic to a unique principal
closed geodesics (see \cite{cdcr2}).

For a principal closed geodesic $b$, we define 
the {\em space of invariants}
as the subspace: 
\[{\mathcal{R}(b)} =\{(\lambda, \tau)| 0< \lambda < 1, 
2/\sqrt{\lambda} < \tau < \lambda + \lambda^{-2}\} \subset \bR^2 \]
which is diffeomorphic to $\bR^2$.

An {\em open} $2$-orbifold is an orbifold with empty boundary whose 
underlying space is noncompact. An $\rpt$-orbifold
with principal closed geodesic boundary component $b'$
is contained in an ambient $2$-orbifold so that $b'$ has 
an annulus neighborhood.

We could interpret the space of invariants as 
a deformation space of germs of convex 
$\rpt$-structures on $b$ as the point of this space determines 
an $\rpt$-structure on a thin neighborhood:
That is, a principal simple closed geodesic is characterized 
by its holonomy along it; two principal simple closed 
geodesics have projectively isomorphic neighborhoods
(in some open ambient $2$-orbifolds)
if and only if they have conjugate holonomy.
We may easily see this since there exists a neighborhood of
a principal simple closed geodesic which is projectively diffeomorphic 
to a quotient under $\langle h(\vth) \rangle$ 
of a domain which is a sufficiently 
thin-neighborhood of the line connecting
the attracting fixed point and the repelling one. 
Such a diffeomorphism is induced by $\dev$.

%\subsection{The deformation space of boundary $1$-orbifolds}
\subsection{A classification of geodesic $1$-orbifolds}
A {\em geodesic} $1$-suborbifold in $\Sigma$ is a suborbifold in
$\Sigma$ so that it is locally modeled on a subspace $\rpo$ in $\rpt$ 
with projective group actions on $\rpt$ preserving $\rpo$.
Since a geodesic full $1$-suborbifold has two points which are mirror points,
the universal cover of a full $1$-orbifold 
is diffeomorphic to an open interval.
(A full $1$-orbifold is two-fold covered by
a circle with one-dimensional projective structure.) 
$\rpt$-structures on a circle are easily classified:
the universal cover of a full $1$-orbifold 
is isomorphic to one of the following $\rpo$-manifolds:
\begin{itemize}
\item $(0,1) $ in $\bR$, considered as an affine patch of $\rpo$.
\item $\bR$ itself.
\item An infinite cyclic cover of $\rpo$.
\end{itemize}
%The universal cover of $1$-orbifold is of the above three types.
In the first case, $l$ is said to be {\em principal}.
A geodesic $1$-dimensional suborbifold of 
a convex two-dimensional orbifold $S$ of negative 
Euler characteristic is 
always principal as a component of the inverse image of it in the universal
cover of $S$ must be isomorphic to the first item.

Let $\Sigma$ be an $\rpt$-orbifold.
A singular point of a principal geodesic full $1$-orbifold in $\Sigma$ 
is modeled on
an open set with an order-two group acting on it fixing a point. The group 
must be generated by an element with a matrix conjugate to:
\[ \begin{bmatrix} -1 & 0 & 0 \\ 0 & -1 & 0 \\ 0 & 0 & 1\end{bmatrix}.\]

Thus, either it is an isolated fixed point of a reflection or 
in a line of fixed point. A principal geodesic full $1$-suborbifold in 
$\Sigma$ 
\begin{description}
\item[cone-type] 
either connects two points which are both cone-points of order two
and lies in the interior of $\Sigma$; 
\item[mirror-type] connects two points which are both mirror points 
and lies in the interior or the boundary of $Q$ entirely and
\begin{description}
\item[boundary-mirror-type] it furthermore is in the boundary,
\item[singular-mirror-type] it furthermore lies in the singular locus 
of $\Sigma$ and connects two corner-reflectors of order two 
(it lies in the interior);
\end{description}
\item[mixed-type] connects a cone-point of order two with 
a mirror point and lies in the interior.
\end{description}

%A geodesic half $1$-orbifold in $\Sigma$ is of:
%\begin{description}
%\item[cone-type] has a singular point in a cone-point of order two, 
%\item[mirror-type] has a singular point in a mirror point
%\begin{description}
%\item[boundary-mirror type] and furthermore lies in the boundary
%\item[singular-mirror type] 
%and furthermore lies in the singular locus of $\Sigma$ and has a singular
%point in a corner-reflector of order two.
%\end{description}
%\end{description}

A geodesic segment is of {\em singular-type} 
if it lies in the singular locus of $\Sigma$.
We say that such a segment is a {\em singular segment}.

Notice from these that a boundary component of a compact convex
$2$-orbifold always is a closed curve 
or a full mirror-type $1$-orbifold.

\subsection{Cross ratios}
In order to introduce invariants of the boundary components,
we recall: 
\begin{defn}\label{defn:cross-ratio} 
Let $y, z, u, v$ be four distinct collinear points with 
$u = \lambda_1 y + \lambda_2 z$ 
and $v = \mu_1 y + \mu_2 z$. The cross-ratio $[y, z; u, v]$ is 
defined to be $\lambda_2\mu_1/\lambda_1 \mu_2$. 
\end{defn}
In particular, the cross-ratio $[0,\infty;1,z]$ equals $z$.

The cross-ratio of four concurrent lines is also defined similarly
(see Busemann-Kelly \cite{BK:53}) using the dual projective
plane where they become four collinear points. 
Another convenient formula is 
given by 
\[[y,z; u, v] = \frac{(\bar u - \bar y) (\bar v - \bar z)}
{(\bar u - \bar z)(\bar v - \bar y)}\]
where $\bar x$ is the coordinate of an affine coordinate 
system on the line containing $y, z, u, v$. (See Berger \cite{Berger:87} for
definitions of affine coordinates.)

%% is this really helpful? Where is this used later?

% It might be helpful to recall that the cross ratio becomes its inverse, if 
% one switches the elements in one of pairs $(u,v)$ and $(y,z)$.
% The cross ratio is invariant if one changes the pair
% $(u,v)$ with $(y,z)$. It is invariant if one 
% switches the elements in both pairs $(u,v)$ and $(y,z)$.
% Also, if a cross ratio is $r$, 
% then permutations of $u,v,y,z$ give rise to
% only six values \[r, 1/r, 1-r, 1/(1-r), 1-(1/r), r/(r-1).\]
% In fact the subgroup of the permutation group not changing
% the cross ratio is generated by above two elements. 
% C. delete 

For example, if $y=1, z =0$, and $1 > u > v > 0$, then the cross ratio
$[1,0,u,v]$ equals \[\frac{1-u}{u}\frac{v}{1-v}\]
which is positive and can realize
any values in the open interval $(0,1)$.

Given a notation $[y,z; u, v]$ with four points 
$y, z, u, v$, they are to be on an image of a segment 
under a projective map where $y, z$ the endpoints and 
$y,v$ separates $u$ from $z$. This is the standard 
position of the four points in this paper.

\subsection{The deformation spaces of full $1$-orbifolds}
A principal geodesic full $1$-orbifold is covered by a component $l$
of its inverse image in $\Sigma$.  Since $l$ is projectively
diffeomorphic to $(0,1)$, $\dev|l$ is an embedding onto a line, a
precompact subset, in an affine patch of $\rpt$.  The holonomy
group of the $1$-orbifold is generated by two reflections $r_1$ and $r_2$.

Let us discuss mirror-type $1$-orbifolds first.  The lines of fixed
points of two reflections $r_1$ and $r_2$ meet $\dev(l)$ at two points
$p_1$ and $p_2$ respectively.  Since $\dev(l)$ is invariant under
$r_1$ and $r_2$, the respective isolated fixed points $f_1$ and $f_2$
of $r_1$ and $r_2$ lie in the one-dimensional subspace containing
$\dev(l)$.  The points $f_1$ and $f_2$ may not coincide with any of
$p_1$ or $p_2$, as $\dev(l)/\langle r_1, r_2 \rangle$ is a
$1$-orbifold isomorphic to $l$ itself and so
$\langle r_1, r_2 \rangle$ acts properly discontinuously on $\dev(l)$.
Also $f_1$ may not coincide with $f_2$ since then $\dev(l)$ is
projectively diffeomorphic to an entire affine line $\bR$ (In this
case, $l$ is said to be an affine $1$-orbifold).

Since $l$ is projectively diffeomorphic to $(0,1)$, the four points
$f_1, f_2, p_1$, and $p_2$ are distinct, and these points are located
on a segment in an affine patch with endpoint $f_1$ and $f_2$ so that
$p_1$ separates $f_1$ and $p_2$.  Otherwise, we obtain a noninjective
developing map, a holonomy element with non-real eigenvalues, or
affine $1$-orbifolds, which contradicts principality.  The cross ratio
$[f_1, f_2; p_1, p_2]$, which is in the interval $(0,1)$, of these
points is invariant under choices of $\dev$ or the conjugation of holonomy.

There alway is an affine 
coordinate so that $(f_1, f_2, p_1, p_2) = (1,0,y,x)$
where $0<x<y<1$. The cross ratio equals
\[ \frac{1-y}{y}\frac{x}{1-x} \]
and hence it is positive and may assume any 
value in $(0,1)$.

Conversely, if two $1$-orbifolds $l_1$ and $l_2$ of mirror type 
have the same invariants, 
then there exist isomorphic neighborhoods
(in some ambient $\rpt$-orbifolds).
% (This again can be seen from $\dev|\tilde \Sigma$ and
% the fact that reflections are unique up to conjugations. 
% Thus, a pair of reflections are determined by the cross-ratio
% above.)
% C. delete
If the interiors of the underlying spaces of $l_1$ and $l_2$ lie 
in the interiors of projective $2$-orbifolds, then the neighborhoods 
can be chosen to be open $2$-orbifolds. If $l_1$ and $l_2$ lie in
the boundary, then the $2$-orbifolds can be enlarged so that 
the neighborhoods become open. The same can be said for 
each of the boundary-mirror-orbifold case or 
the singular-mirror-orbifold case.

% the use of the word ``again'' is very repetitive

For orbifolds of cone-type,
the isolated fixed points $f_1$ and $f_2$ 
lie on $\dev(l)$. Let $p_1$ and $p_2$ be the points of
intersection of the lines of fixed points of $r_1$ and $r_2$ 
meet the one-dimensional subspace containing $\dev(l)$. 
The points $f_1, f_2, p_1, p_2$ are distinct, and lie on
a segment with endpoints $p_1$ and $p_2$. We assume that
$f_1$ separates $p_1$ and $f_2$. The cross-ratio
$[f_1, f_2; p_1, p_2] \in (0,1)$ 
is independent of the choice of
of $\dev$. Conversely, if two $1$-orbifolds of cone-type have 
the same invariants, then there exists 
isomorphic open neighborhoods (possibly after extending the collar
neighborhoods).

The mixed-type case is entirely similar
with invariant defined by $[f_1, f_2; p_1,p_2]$ again
for $f_1, f_2, p_1, p_2$ defined as above.
(Here, $f_1$ and $f_2$ are not the endpoints of a segment.)

We remark that a principal $1$-orbifold has 
a double-covering circle $s$, and
the generator of whose fundamental group has holonomy which
is hyperbolic with eigenvalues $\lambda, 1, \lambda^{-1}$; i.e.,
it is purely hyperbolic.
This can be easily seen since the intersection point of 
the lines of fixed points of $r_1$ and $r_2$ correspond to
an eigenvector with eigenvalue $1$ for the generator.
%(In this paper, we will only need to consider type  
%boundary mirror $1$-orbifolds.)

Given a full principal geodesic $1$-orbifold $b$ in $\Sigma$,
the {\em space of invariants} is defined
as  $\mathcal{C}(b) = (0,1)$.

\begin{defn}\label{defn:bddeform}
Given a convex $2$-orbifold $\Sigma$, 
let $\partial \Sigma$ denote the union of boundary
$1$-orbifolds. Let $\mathcal{C}(\partial \Sigma)$ denote
the product of the spaces of invariants of
all components of $\partial \Sigma$.
\end{defn}

As with closed case, $\mathcal{C}(\Sigma)$ is 
a subspace of $\rpt(\Sigma)$ and 
$\mathcal{T}(\Sigma)$ is a subspace of 
$\mathcal{C}(\Sigma)$.

%%4/26 3:30

\subsection{Geometric constructions of $\rpt$-orbifolds}
\label{subsec:geomop}
Now, let $\Sigma'$ be a compact convex $\rpt$-orbifold
with principal boundary. Given that certain boundary
conditions are met, we will describe how
to obtain a convex  $\rpt$-orbifold $\Sigma''$ obtained from $\Sigma'$ by
the above topological operations in \S \ref{subsec:topop}
and construct all convex 
structures on $\Sigma''$ so that $\Sigma'$ with
its original convex structure is 
obtained back when we split. We of course obtain 
principal boundary for the resulting $2$-orbifolds. 

We follow the notation of \S\ref{subsec:topop}.

\subsubsection{Pasting or crosscapping {\rm (A)(I)}}
In the former case, if 
two boundary component curves $b$ and $b'$ 
%parameterized by $f$ 
%C. delete
have a conjugate holonomy, then $\Sigma''$ is also 
a projective $2$-orbifold. This construction is given in 
\S3.6 of \cite{Gconv:90} for surface cases.
The constructions are called {\em pasting} or {\em crosscaping}
depending on whether the curve is two-sided or one-sided respectively.

Here is the construction. We distinguish two cases, depending on whether
$b$ and $b'$ are non-isomorphic or isomorphic respectively.
% (1) 
% why is this numbered?
% C. delete

Let $b$ and $b'$ be distinct with equal invariants.  Suppose also that
the component $\Sigma'$ of $\Sigma$ containing $b'$ is distinct from
the component $S$ of $\Sigma$ containing $b$. Without loss of generality,
assume that $\Sigma'$ has the two components only.
Let $\tilde S$ be the universal cover of $S$ 
and $(\dev_S, h_S)$ the development pair.  Let
$\tilde \Sigma'$ and $(\dev_{\Sigma'}, h_{\Sigma'})$ be the universal
cover and the pair for $\Sigma'$.  Let $l$ and $l'$ denote components
of inverse images of $b$ and $b'$ in $\tilde S$ and $\tilde \Sigma'$
respectively. Let $\vth$ and $\vth'$ denote the deck transformations
corresponding under $f$.  Then $h_S(\vth)$ and $h_{\Sigma'}(\vth')$
act on $\dev_{S}(l)$ and $\dev_{\Sigma'}(l')$ respectively.  Since
$h_S(\vth)$ and $h_{\Sigma'}(\vth')$ are conjugate,
\begin{equation}\label{eqn:conj2} 
f' h_S(\vth) f^{\prime -1} = h_{\Sigma'}(\vth')
\end{equation}
for some collineation $f'$.
Let $\Omega = \dev_S(\tilde \Sigma)$ 
and $\Omega' = \dev_{\Sigma'}(\tilde \Sigma')$. 
Then by post-composing $f'$ with a reflection if necessary, we may 
assume without loss of generality that  
$f'(\Omega)$ and $\Omega'$ meet exactly in $f'(\dev(l)) =\dev(l')$. 

%%% 
%%% this is vague and should be rewritten: ``a supportine line consideration''
%%% 
%C. rewrite

Their union is convex: $f'(\Omega)$ is a subset of
$f' h_S(\vth) f^{\prime -1}$-invariant triangle with 
an open side $f'(\dev(l))$ 
and $\Omega$ is a subset of $h_{\Sigma'}(\vth')$-invariant 
triangle with side $\dev(l')$.
The second triangle must be adjacent to the first one. 
Since a supporting line of $f'(\Omega)$ at a vertex of $\dev(l')$ 
coincide with a side of the first triangle and that of $\Omega'$ 
coincide with the second triangle and the sides extend 
each other being the sides of the invariant triangles, 
it follows that the support lines coincide. The same holds 
at the other vertex of $\dev(l')$. Elementary geometry shows that  
$f'(\Omega) \cup \Omega'$ is convex.  

Let $\Gamma$ be the image of the homomorphism 
$f' h_S(\cdot) f^{\prime -1}$ and $\Gamma'$ that of
$h_{\Sigma'}$. Let $\Gamma''$ be the group generated
by $\Gamma$ and $\Gamma'$, which is isomorphic to an amalgamated product 
of $\Gamma$ and $\Gamma'$ actually.
Let $\Omega''$ be the union of images of 
$\Omega$ and $\Omega'$ under $\Gamma''$.
We claim that $\Omega''$ is a convex domain: 
any two points lie in a finite connected union of images of $\Omega$.
A finite connected union is always convex.
We can order the images and keep adding domains one by one.
Let $\Omega_n$ be the $n$-th union. 
At each step, the union of $\Omega_n$ and a new domain 
$\Omega'$ to be added meet at a line acted upon by 
a hyperbolic transformation, which is a conjugate of 
$h_{\Sigma'}(\vth')$. Moreover, $\Omega_n$ and 
$\Omega'$ are subsets of adjacent triangles acted upon by 
the same hyperbolic transformation. 
The above supporting line argument applies. 

Since $\Gamma''$ acts properly on 
$\Omega''$, we see that $\Omega''/\Gamma''$ is 
a compact convex $2$-orbifold $\Sigma''$ obtained from $S$ 
and $\Sigma'$ by pasting along $b$ and $b'$ and leaving other 
components untouched.

%%%
%%% this last paragraph should be rewritten
%%%

\begin{figure}[ht]
\centerline{\epsfxsize=2.5in \epsfbox{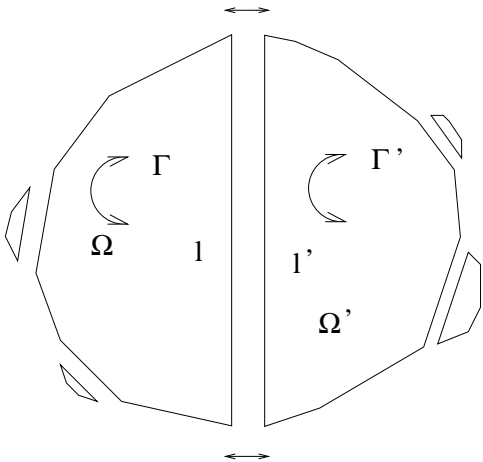}}
\caption{\label{fig:pasting} 
A convex orbifold obtained by pasting two smaller convex orbifolds.}
\end{figure}
\typeout{<<pasting.eps>>}

\begin{lem}\label{lem:subconv}
A subsurface $S'$ of a convex surface $S$ bounded by 
closed geodesics is convex.
\end{lem}
\begin{proof} 
A path in the subsurface $S'$ is homotopic to a geodesic 
in $S$. Since there are no bigons, the geodesic itself is 
in $S'$. A closed geodesic in $S$ is always principal.
Thus, $S'$ is convex (see \cite{cdcr2}).
\end{proof}

Conversely, given a convex $\rpt$-structure on 
$\Sigma''$ with a principal geodesic $b''$ 
as above, we see that the completions of 
$\Sigma'' - b''$ is a convex 
$2$-orbifold by taking a finite cover and Lemma 
\ref{lem:subconv}.
Therefore, such a structure can always be constructed 
in this manner. 

In fact, the choice of $f'$ is not unique. $f'$ can be replaced by
\begin{equation}\label{eqn:conjugating}
\hat f = g \circ f' 
\end{equation}
where $g$ is in the identity component of 
$\PGL(3, \bR)$ commuting with $h_{\Sigma'}(\vth')$.
Thus, $g$ is diagonalizable with positive eigenvalues,
and hence the space of such $g$ is $\bR^2$.

Let $h''_{\hat f}: \pi_1(\Sigma'') \ra \PGL(3, \bR)$ be given 
by amalgamating the fundamental groups and extending
homomorphisms $h_{\Sigma}$ and $h_{\Sigma'}$ in
the obvious manner with image group $\Gamma''$. 
There is an $\bR^2$-parameter space of possible $\hat f$. 
For different choices of $\hat f$, $h''_{\hat f}$
yields non-conjugate actions.
Hence the resulting $\rpt$-structures are non-isotopic
(Lemma \ref{lem:fourpts}). For 
given fixed invariants for $b$ and $b'$ and fixed 
$\Sigma'$, this construction gives
an $\bR^2$-parametrized family of non-isomorphic $\Sigma''$.

\begin{lem}\label{lem:fourpts} 
Let $h$ be a holonomy homomorphism of a convex $2$-orbifold of negative 
Euler characteristic with possibly nonempty but 
principal geodesic boundary. Let $f_t$ be a one-parameter 
family of collineations.
Then $f_t h(\cdot) f^{-1}_t$ equals
$f_0 h(\cdot) f_0^{-1}$ if and only if $f_t$ is constant. 
\end{lem}
\begin{proof}
Pass to a finite covering surface $S$ with $\chi(S)<0$.
Let $(\dev, h)$ be a development pair and $\tilde S$
the universal cover of $S$.  Since $S$ is convex, its holonomy representation
is faithful.  Let $g_1$ and $g_2$ be two
deck-transformations corresponding to simple closed curves not
homotopic to each other.  Then these curves are homotopic to principal
closed geodesics and $g_1$ and $g_2$ act on geodesic lines $l_1$ and
$l_2$ in the universal cover $\tilde S$ corresponding to the closed
geodesics respectively.  The endpoints of $\dev(l_1)$ are the
attracting and repelling fixed points of $h(g_1)$ by the principal
conditions; and the endpoints of $\dev(l_2)$ those of $h(g_2)$. No two
of these points coincide. Otherwise, $\dev(l_1)$ can be sent
arbitrarily close to $\dev(l_2)$ by $h(g_2^n)$. Since $l_1$ map to a
simple closed curve and $\dev$ is an embedding, this cannot happen.
No three of these points are collinear since $h(g_i)$ acts on
$\dev(l_i)$ freely, and by convexity.

Suppose that $f' h(\cdot) f^{\prime -1} = h(\cdot)$ for a collineation
$f'$. Then $f' h(g_i) f^{\prime -1} = h(g_i)$ for
$i=1,2$, and $f'$ acts on each of the two pairs of four noncollinear
points.  Since the images of four points, no three of which are
collinear, determine the collineations, $f'$ is the
identity map or a unique reflection determined by the four
points. Therefore $f_t = f'$ is constant.
\end{proof}

We describe this construction in a different language following
\cite{Gconv:90}: First find a slightly bigger $\rpt$-orbifold $T$
containing $\Sigma$ so that $T -\Sigma$ is a union of two annuli
parallel to $b$ and $b'$ respectively.  Then there are open tubular
neighborhoods of $b$ and $b'$ which are isomorphic as $b$ and $b'$ are
conjugate elements of $\pi_1(\Sigma)$ .  Remove from $T$ what are
outside these annuli to obtain $T'$. Now identify these two annuli by
a projective diffeomorphism $f'$. As $b$ and $b'$ are principal
geodesics, $f'$ sends $b$ to $b'$.  Then $\Sigma''$ is
independent of the choice of annuli but depends on the germ of
$f'$ near $b$.

There is a two-parameter family of projective diffeomorphisms
with corresponding annular neighborhoods of $b$ and $b'$. 
First consider the annuli as quotients of domains $D_1$
and $D_2$ in $\rpt$ and a collineation $\tilde f$
sending $D_1$ to $D_2$ lifting $f'$. Such $\tilde f$ satisfies
\[\tilde f \circ \vth_1 \circ \tilde f^{-1} 
= \vth_2\] where $\vth_1$ and $\vth_2$ are generators
of the infinite cyclic groups acting on 
$D_1$ and $D_2$ respectively. Therefore, 
there are choices of maps $f$ parameterized by $\bR^2$. 
Essentially by Lemma \ref{lem:fourpts}, different choices
of $\tilde f$ yield non-conjugate holonomy groups, 
and hence non-isotopic $\rpt$-structures. 
(This is the projective version of Fenchel-Nielsen twists
for hyperbolic surfaces.)

To summarize: a family of distinct $\rpt$-structures on
$\Sigma''$ is parametrized by $\bR^2$ when the common conjugacy class of
the holonomy of $b$ and $b'$ is fixed and $\Sigma'$ is fixed.
The group $\bR^2$ acts on the deformation space 
of $\Sigma''$ by changing the gluing map as above. 
Thus, we obtain a principal $\bR^2$-fibration 
description as in Proposition \ref{prop:principal1}.
See \S5 of \cite{Gconv:90} for details, 

We may also assume that $S = \Sigma'$ but $b$ not equal to $b'$. 
In this case, the discussions are similar with $\Omega$ and 
$\Omega'$ being equal and $\Gamma$ becoming an HNN-extension.

% (2) 
% need to organize this better
% C. rewrite
%%rewrite on set up and many bad disorganizations
%4/29 5:10

We now go over the geometric crosscapping construction.
Suppose now that $b = b'$. In this case, $S = \Sigma'$,
and $b$ corresponds to a simple closed curve $b''$
in $\Sigma''$ with a M\"obius band neighborhood. 

Let $(\dev, h)$ denote the development pair of $S$, and $\tilde S$ the
universal cover of $S$.  Then let $\Gamma = h(\pi_1(S))$.  
Let $\tilde b$ denote a component of the inverse image of $b$ and $\vth$ be the
corresponding deck transformation acting on $\tilde b$.  Let $\vth'$
be the unique projective automorphism acting on $\tilde b$ preserving
an orientation of $\tilde b$ but reversing the orientation of 
$\rpt$ so that $\vth^{\prime 2} = \vth$; that is, we want $\vth'$ to be a
hyperbolic slide reflection.  (This can obviously solved by
conjugating the hyperbolic $\vth$ to a diagonal form.)  Let $\Omega =
\dev(\tilde S)$. $\vth'(\Omega)$ and $\Omega$ meet exactly at
$\dev(\tilde b)$. Since $\Omega$ and $\vth'(\Omega)$ are
$\vth$-invariant, $\Omega \cup \vth'(\Omega)$ is a convex domain
similarly to the pasting case. Let
$\Gamma''$ denote the group generated by $\Gamma = h(\pi_1(\Sigma))$
and $\vth'$. Let $\Omega''$ be the union of images of $\Omega$ under
$\Gamma''$.  Then $\Omega''$ is a convex domain as in the pasting case, 
and $\Omega''/\Gamma''$ is an orbifold diffeomorphic to a component
orbifold of $\Sigma''$.

Any convex $\rpt$-orbifold diffeomorphic 
$\Sigma''$ can be constructed in this manner. 
Finally, we remark that given a fixed invariant on $b$, there is 
a unique $\Sigma''$ that can be constructed. 

\subsubsection{Silvering and folding {\rm (A)(II)}}\label{subsec:AII}
In this case, $f$ either 
(1) is the identity map or (2) has exactly two fixed points 
reversing the orientation of $b$ and of order two. 
%(Otherwise, 
%we cannot get manifold underlying space for $\Sigma''$.) 
% delete
Let $\tilde S$ be the universal cover of a component $S$
containing $b$ and $(\dev, h)$ the development pair of $S$. Let
$\tilde b$ be a component of the inverse image of $b$ in $\tilde S$
and $\vth$ the corresponding deck transformation acting on $\tilde b$.

(1) When $f$ is the identity, there is a unique reflection 
$F:\rpt \ra \rpt$ so that 
the line of fixed points contain $\dev(\tilde b)$ and
$F\circ h(\vth) \circ F^{-1} = h(\vth)$.  Thus the isolated fixed
point of $F$ coincides with the fixed point of hyperbolic automorphism
$h(\vth)$ not on the closure of $\dev(\tilde b)$.  As above, 
consider the group $\Gamma''$ generated by $\Gamma$ and $F$, and the
union $\Omega''$ of images of $\dev(\tilde S)$ under the action of
this group. $\Omega''$ is a convex domain.  Then 
$\Omega''/\Gamma''$ as a component of $\Sigma''$.  As $F$ is unique,
$\Sigma$ determines $\Sigma''$.
%(We use the terminology that $b$ has been transformed to a mirror, or 
%$b$ has been {\em silvered}.)
%C. don't need this also.
\label{page:silvered}

When the holonomy of $b$ is hyperbolic and $b$ is
geodesic, then $b$ can be always silvered since such an
element $F$ exists.  A boundary component with quasi-hyperbolic holonomy
cannot be silvered.

(2) When $f$ has exactly two fixed points,
there is a reflection $F$ so that the isolated
fixed point of $F$ lies on $\dev(\tilde b)$ and
\begin{equation}\label{eqn:rF}
F\circ h(\vth) \circ F^{-1} = h(\vth)^{-1}.
\end{equation}
This forces the hyperbolic $h(\vth)$ to have eigenvalues $\lambda, 1,
\lambda^{-1},$ for $\lambda > 1$, and $F$ exchanges the endpoints of
$\dev(\tilde b)$ and fixes the third fixed point of $h(\vth)$, which
is not an endpoint of $\dev(\tilde b)$.  The choice of isolated
fixed point of $F$ on $\dev(\tilde b)$ determines the intersection
of the fixed line of $F$ with  the line 
containing  $\dev(\tilde b)$.  
Let $\Gamma''$ be generated by $\Gamma$ and $F$, 
and $\Omega''$ be the union of images of $\dev(\tilde \Sigma)$ under
$\Gamma''$, which is again convex.  (This construction produces two
cone-points of order two, one corresponding to $F$ and the
other to $h(\vth)\circ F$, which by \eqref{eqn:rF} is a reflection.)

\begin{figure}[ht]
\centerline{\epsfxsize=2.5in \epsfbox{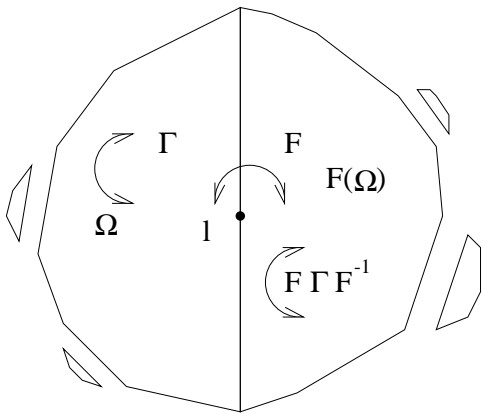}}
\caption{\label{fig:pastecone} 
How to fold an orbifold.}
\end{figure}
\typeout{<>}

The choice of the isolated fixed points in $\dev(\tilde b)$ itself
produces non-isotopic $\rpt$-structures on $\Sigma''$.
% essentially by Lemma \ref{lem:fourpts}.  
% delete
Thus, there is an $\bR$-family of non-isotopic $\rpt$-structures on
$\Sigma''$ for $h(\vth)$ in a fixed conjugacy
class of a  purely hyperbolic transformation.  
Any other choice of $F$ equals $g\circ F \circ g^{-1}$ for
unique $g$ commuting with $h(\vth)$ and with positive eigenvalues
including $1$. Such $g$ is diagonalizable with positive
eigenvalues. Thus $\bR$ acts on the deformation space of convex
$\rpt$-structures on $\Sigma''$.

\subsubsection{Pasting {\rm (B)(I)}}
Now suppose that $b$ is a full $1$-orbifold.  
Consider a diffeomorphism $f:b \ra b'$ with another
$1$-orbifold $b'$.

%Let $b'$ be another one, and a diffeomorphism 
%constructed topologically. 

The holonomy group of $b$ is generated by 
two reflections $r_1$ and $r_2$ acting on the line $l$ to which
a lift of $b$ to $\tilde \Sigma$ develops. The
fixed lines of $r_1$ and $r_2$ are transversal to $l$.
The respective fixed points $p_1$ and $p_2$ of $r_1$ and
$r_2$ lie on $l$. Letting $q_1$ and $q_2$ denote the respective
intersection points of fixed lines of $r_1$ and $r_2$ with $l$, 
the cross ratio $[p_1, p_2; q_1, q_2]$ determined the respective
conjugacy classes of $r_1$ and $r_2$. 
If $b$ and $b'$ have equal invariants, we can find 
a projective automorphism conjugating the reflections 
$r_1$ and $r_2$ to that corresponding to $b'$. 
The pasted orbifold
$\Sigma''$ carries a convex $\rpt$-structure  by a similar construction.
%The resulting  $\rpt$-manifold  $\Sigma''$ is convex. 

%The appropriate double cover
%of $\Sigma$ is pasted along a principal closed geodesic
%obtaining a double cover of $\Sigma''$. 
The $\bR$-family of
conjugating elements  determines an $\bR$-family of
non-isotopic $\rpt$-structures on $\Sigma''$
for fixed invariants for $b$ and $b'$.  

\begin{figure}[ht]
\centerline{\epsfxsize=3.2in \epsfbox{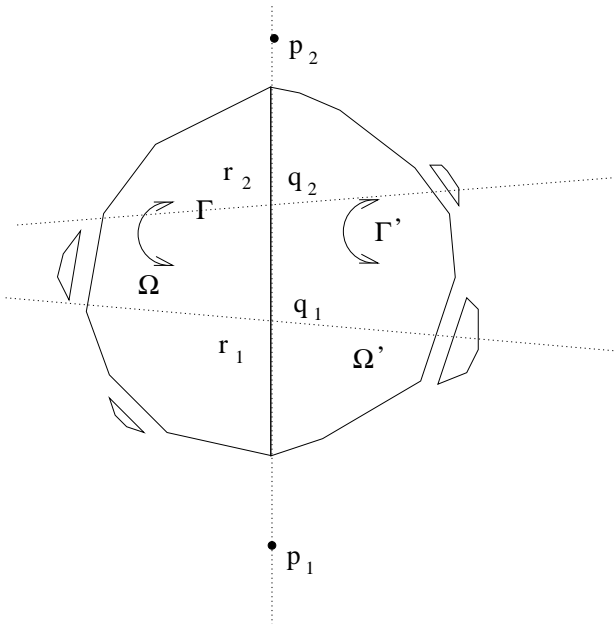}}
\caption{\label{fig:paste1orb} 
Pasting along a $1$-orbifold.}
\end{figure}
\typeout{<>}

% 4/29 1:35
\subsubsection{Silvering and folding {\rm (B)(II)}}
Now suppose $b' =b$. Then $f$ is the identity map or 
fixes a unique point. If $f$ is the identity, $b$ contains
mirror points.  Then $\Sigma$ determines  $\Sigma''$ uniquely.
% We thus {\em silvered} $b$.  (This is similar to (A)(II)(1).)
%\label{page:silvered2}
%sc. I think we need this actually.

% (2) 
If $f$ fixes a unique point, $\Sigma'' - (\Sigma - b)$ is a
$1$-orbifold with one endpoint a cone-point of order 2 and the other
endpoint a mirror point. Then a unique convex $\rpt$-structure on $\Sigma''$
exists,  given one in $\Sigma$ as above in 
(A)(II)(2). This follows since the gluing automorphism has to switch two
reflections by conjugations.  
%(We are said to have {\em folded} $b$.)
%sc. I think we need this.
A fixed invariant for $b$ determines a unique convex structure on $\Sigma''$.

These operations (A)(I), (A)(II), (B)(I), and (B)(II) can be performed for
$\rpt$-structures on possibly nonconvex orbifolds since
these operations are supported on components of the boundary. 

%The details are described in the next section.
% (We give the details next.)  
% (For convex ones with principal boundary, the good
% points are that the operations preserve convexity with principal
% boundary condition.)
(However, all these operations preserve convexity.)

\subsection{Deformation spaces for bounded orbifolds}
% Finally let us now discuss the deformation spaces for orbifolds with
% nonempty boundary: 

% From now on, we use the convention that an $\rpt$-
% orbifold with boundary is {\em convex} if it is convex and has
% principal boundary.  Hyperbolic $\rpt$-structures are convex in this 
% convention.

%% this is said in the following paragraph.
%C. OK. 

When  $\partial\Sigma\neq\emptyset$, we consider
$\rpt$-structures on $\Sigma$ so that its boundary components are
principal (see Remark \ref{rem:qh}).  (Clearly convex structures
belong here.)  % Thus, we define 
Denote by $\rpt(\Sigma)$ the subspace of
the deformation space of all $\rpt$-structures with principal
boundary components. 
$\rpt(\Sigma)$ is the quotient space by the action of $\PGL(3, \bR)$ as
before.  Given an element of $\rpt(\Sigma)$ for an orbifold $\Sigma$, one
can silver them producing an $\rpt$-structure on $\Sigma'$ with
silvered boundary since there exists a unique involution
centralizing the holonomy of the boundary component.
% ($\Sigma'$ is a boundaryless orbifold.)  
For any $\rpt$-orbifold $\Sigma'$, 
clarifying the boundary components produces
$\rpt$-structures on $\Sigma$.  
Thus,   
$\mathcal{S}(\Sigma)$ bijectively corresponds to $\mathcal{S}(\Sigma')$ 
using the silvering and clarifying operations.
Thus $\rpt(\Sigma)$ bijectively corresponds to $\rpt(\Sigma')$.
(A bit of subtlety lies with the fact that isotopies on $\Sigma'$ need
to have some symmetry for new mirror-point sets while isotopies on
$\Sigma$ need not. However, the charts of $\Sigma'$ need to have the
symmetries also. The net effect is the homomorphism. )  

We define $\mathcal{C}(\Sigma)$ as a subspace of $\rpt(\Sigma)$
of elements represented by convex structures on 
$\Sigma$ with principal boundary.

A similar bijection exists 
between a subset of $\mathcal{C}(\Sigma)$ and
$\mathcal{C}(\Sigma')$. 

%which correspond to 
%adding and removing a reflection element 
%to the holonomy group since $\mathcal{C}(\Sigma')$ consists 
%of conjugacy classes of discrete representations.

\begin{thm}\label{thm:silvering}
The process of silvering and clarifying 
induces a one-to-one correspondence between 
$\rpt(\Sigma)$ and $\rpt(\Sigma')$ 
where $\Sigma'$ is obtained from silvering 
$\Sigma$ in the topological sense. 
The same can be said for ${\mathcal{C}}(\Sigma)$ 
and ${\mathcal{C}}(\Sigma')$.
\end{thm}

When $\partial\Sigma\neq\emptyset$, 
denote by $H(\Sigma)_p$ the open subset of 
$H(\Sigma)$ consisting of structures for which 
the holonomy of each component of $\partial\Sigma$ 
is hyperbolic.

Define as above
\[U(\Sigma)_p = H(\Sigma)_p \cap U^g\]
for appropriate $U^g$. We can form a one-to-one 
correspondence between $U(\Sigma)_p$ and 
$U(\Sigma')$ by adding or removing reflections
corresponding to the boundary components of $\Sigma$.
The correspondence is obviously a homeomorphism,
which is not proper.

\begin{thm}\label{thm:Haus2} 
Let $\Sigma$ be a connected compact $2$-orbifold with 
$\chi(\Sigma)<0$ and $\partial\Sigma\neq\emptyset$.
Then $\rpt(\Sigma)$ has the structure of 
a Hausdorff real analytic variety modeled on $U(\Sigma)_p/\PGL(3, \bR)$,
for which the induced map 
\[\mathcal{H}: \rpt(\Sigma) \ra U(\Sigma)_p/\PGL(3, \bR)\] 
is a local homeomorphism onto its image.
\end{thm}
\begin{proof}
Since the silvered $\Sigma'$ has empty boundary, 
\[\mathcal{H}': \rpt(\Sigma') \ra U(\Sigma')/\PGL(3, \bR)\]
is a local homeomorphism. The one-to-one correspondences 
discussed above obviously makes a commutative diagram: 
\begin{eqnarray*}
\rpt(\Sigma)&  \stackrel{\mathcal{H}}{\longrightarrow} 
& U(\Sigma)_p/\PGL(3, \bR) \\
\updownarrow &  & \updownarrow \\
\rpt(\Sigma') & \stackrel{\mathcal{H}'}{\longrightarrow} 
& U(\Sigma')/\PGL(3, \bR).
\end{eqnarray*}
The desired conclusion follows.
\end{proof}

% 
% this is not an acceptable proof
%
%C. Write details!

\begin{prop}\label{prop:openness2} 
Let $\Sigma$ be a compact connected orbifold 
with $\chi(\Sigma) < 0$ and $\partial \Sigma \ne \emp$. 
Then $\mathcal{C}(\Sigma)$ is 
an open subset of $\rpt(\Sigma)$. 
\end{prop}
\begin{proof} 
By Proposition \ref{prop:openness}, $\mathcal{C}(\Sigma')$ is 
an open subset of $\rpt(\Sigma')$ where $\Sigma'$ is obtained 
from $\Sigma$ by silvering boundary components. 
By the above one-to-one correspondence between 
$\mathcal{C}(\Sigma)$ and $\mathcal{C}(\Sigma')$ 
and that $\rpt(\Sigma)$ between $\rpt(\Sigma')$, 
we obtain the result.
\end{proof}

\begin{prop}\label{prop:closedness2}
The image of ${\mathcal{C}}'(\Sigma)$ 
under $\mathcal{H}'$ to $U(\Sigma)_p$
is a {\rm (}relatively\/{\rm )} closed subset.
\end{prop}
\begin{proof}
We have the following commutative diagram:
\begin{eqnarray*}
{\mathcal{C}}(\Sigma)&  \stackrel{\mathcal{H}}{\longrightarrow} 
& U(\Sigma)_p/\PGL(3, \bR) \nonumber \\
\updownarrow &  & \updownarrow \nonumber \\
{\mathcal{C}}(\Sigma') & \stackrel{\mathcal{H}'}{\longrightarrow} 
& U(\Sigma')/\PGL(3, \bR).
\end{eqnarray*}
Since the image of $\mathcal{H'}$ is closed, our result follows.
\end{proof}

%4/29 8:55 
%% I need to look again to streamline here.

\subsubsection{A question of Fock}

We will not need the following theorem in this paper, 
but to answer a question (for bounded surfaces) of Vladimir Fock, 
we state:
\begin{thm}\label{thm:bdHitchin} 
Let $\Sigma$ be a compact connected orbifold of
negative Euler characteristic and with nonempty boundary. 
Then 
\[\mathcal{H}: \mathcal{C}(\Sigma) \ra U(\Sigma)_p\]
is a homeomorphism onto a component of $U(\Sigma)_p$ 
containing the image of $\mathcal{T}(\Sigma)$
said to be the {\em Hitchin-Teichm\"uller component}. 
\end{thm}
\begin{proof}
Theorem \ref{thm:Haus} and Proposition \ref{prop:openness2} 
show that the image under $\mathcal{H}$ is an open set. 
The proof follows from Proposition \ref{prop:closedness2}
and the fact that $\mathcal{C}(\Sigma)$ is connected
by Theorem \ref{thm:finaldim}.  
\end{proof}

\subsection{Fibrations of 
deformation spaces %corresponding to the operations
}

\begin{prop}\label{prop:principal1} 
Consider only compact $2$-orbifolds of negative Euler characteristic.
% with or without boundary

\begin{description} 
\item[(A)(I)(1)] 
Let the $2$-orbifold $\Sigma''$ be obtained 
from pasting along two closed curves $b, b'$ 
in a $2$-orbifold $\Sigma'$.
The map resulting from splitting
\[\mathcal{SP}: \mathcal{C}(\Sigma'') \ra 
\Delta \subset \mathcal{C}(\Sigma')\]
is a principal $\bR^2$-fibration, 
where $\Delta$ is the subset of $\mathcal{C}(\Sigma')$ 
where $b$ and $b'$ have equal invariants.
% determined by letting the invariants of $b$ and $b'$ be equal.
% Moreover, $\bR^2$ acts transitively on each fiber.
\item[(A)(I)(2)] 
Let $\Sigma''$ be obtained from $\Sigma'$ by cross-capping.
The resulting map
\[\mathcal{SP}: \mathcal{C}(\Sigma'') \ra \mathcal{C}(\Sigma')\]
is a diffeomorphism.
\item[(A)(II)(1)] Let $\Sigma''$ be obtained from $\Sigma'$ by 
silvering. The clarifying map
\[\mathcal{SP}: \mathcal{C}(\Sigma'') \ra \mathcal{C}(\Sigma')\]
is a diffeomorphism.
\item[(A)(II)(2)] Let $\Sigma''$ be obtained from $\Sigma'$ 
by folding a boundary closed curve $l'$. 
The unfolding map 
\[\mathcal{SP}: \mathcal{C}(\Sigma'') \ra \Delta \subset
\mathcal{C}(\Sigma')\]
is a principal $\bR$-fibration, 
where $\Delta$ is a subspace of $\mathcal{C}(\Sigma')$ consisting 
of $\rpt$-structures with hyperbolic holonomy for $l'$.
\item[(B)(I)] Let $\Sigma''$ be obtained by pasting along 
two full $1$-orbifolds $b$ and $b'$ in $\Sigma'$.
The splitting map
\[\mathcal{SP}: \mathcal{C}(\Sigma'') \ra 
\Delta \subset \mathcal{C}(\Sigma')\] 
is a principal $\bR$-fibration where
$\Delta$ is a subset of $\mathcal{C}(\Sigma')$ where 
the invariants of $b$ and $b'$ are equal.
\item[(B)(II)] Let $\Sigma''$ be obtained by silvering 
or folding a full $1$-orbifold.
The clarifying or unfolding map
\[\mathcal{SP}: \mathcal{C}(\Sigma'') \ra \mathcal{C}(\Sigma')\]
is a diffeomorphism.
\end{description} 
\end{prop}

%\marginpar{Should we explain the $\bR^2$-action? and 
%naturality ? }

\section{Decomposition of convex orbifolds into  elementary orbifolds}
\label{sec:decomp}

This section describes how to decompose a compact convex $\rpt$-orbifold
with negative Euler characteristic and principal boundary
into elementary orbifolds along 
disjoint simple closed geodesics and geodesic full $1$-orbifolds.  
Crucial is Lemma \ref{lem:disjointcut},
whereby simple closed curves and $1$-orbifolds are realized
by simple closed geodesics and geodesic $1$-orbifolds.  
%Elementary orbifolds are defined in \S???.
Next, we define elementary orbifolds and prove that 
we can decompose convex orbifolds into elementary ones. 
%This is correct! SC.

% The proof
% essentially follows from choosing the simple closed curves and
% $1$-orbifolds well.

\subsection{The existence of $1$-suborbifolds.}

%An {\em isotopy} is a diffeomorphism 
%$g:\Sigma\longrightarrow\Sigma$ of an orbifold $\Sigma$ 
%is an orbifold map
%$F:\Sigma \times [0,1] \ra \Sigma$ such that 
%the maps $f_t:\Sigma \ra \Sigma$ defined by $f_t(x) = F(x, t)$
%are orbifold diffeomorphisms such that $f_0$ the identity,
%and $f_1=g$.

The $\rpt$-orbifolds are decomposed using the following lemma:
\begin{lem}\label{lem:disjointcut} 
Suppose that $\Sigma$ is a compact convex $\rpt$-orbifold with
$\chi(\Sigma) < 0$ and principal geodesic boundary.  Let $c_1, \dots,
c_n$ be a mutually disjoint collection of simple closed curves or
$1$-orbifolds so that the orbifold Euler characteristic of the
completion of each component of $\Sigma - c_1 - \cdots - c_n$
is negative.  Then $c_1, \dots, c_n$ are isotopic to principal simple
closed geodesics or principal geodesic $1$-orbifolds $d_1, \dots, d_n$
respectively. %Here $c_i$ is isotopic to $d_i$ for each $i$, and hence
% $c_i$ is a $1$-orbifold if and only if $d_i$ is.
%  Also, the isotopy could be chosen simultaneously.
%%C. delete
\end{lem}
\begin{proof}

Consider first the case $n=1$.  The completions of components of
$\Sigma - c_1$ have negative Euler characteristic by Proposition
\ref{prop:eulerspli} since they are obtained by sewing the components
of $\Sigma - c_1 - c_2 -\dots - c_n$.  Let $\Sigma'$ be a finite cover
of $\Sigma$ which is an orientable surface with principal closed
geodesic boundary.  Let $G$ be the group of automorphisms of $\Sigma'$
so that $\Sigma'/G$ is projectively diffeomorphic to $\Sigma$ in the
orbifold sense. (We showed the existence of such a cover in 
\S\ref{sec:orbifolds}.)

Let $c_{1i}$, $i=1, \dots, k$, be the components of the inverse image
of $c_1$.  Since $\Sigma'$ is convex, orientable, and
$\chi(\Sigma')<0$, each $c_{1i}$ is homotopic to a simple closed
geodesic (see \cite{Gconv:90}).  Let $d_{1j}$, $j=1, \dots, l$, be the
union of all simple closed geodesics in $\Sigma'$ homotopic to some 
$c_{1i}$. Then clearly $G$ acts on
\[\bigcup_{j=1, \dots, l} d_{1j}.\]
Let us take a component $d_{1j}$, and let $G_j$ be the subgroup
stabilizing $d_{1j}$. Let $c_{1, j_1}, \dots, c_{1, j_m}$ be
the curves homotopic to $d_{1j}$ in $\Sigma'$,
on the union of which $G_j$ acts on. $G_j$ must act transitively
on the set \[C_j = \{c_{1, j_1}, \dots, c_{1, j_m}\}\]
as $c_1$ is the unique image in $\Sigma$.
If $C_j$ has a unique element, then $G_j$ also act on
the unique element, and the image of $d_{1j}$ is
the unique principal closed geodesic or principal $1$-orbifold
$d_1$ isotopic to $c_1$. % since we can build a $G$-equivariant 
% homotopy using and we are done. 
% C. Check it again.
(A $G$-equivariant isotopy 
in $\Sigma'$ can be constructed by perturbations of $c_1$ and $d_{1j}$s
and using innermost-bigons which occur $G$-equivariantly.)

If $C_j$ has more than one element, then let $A_j$ be the unique
maximal annulus bounded by elements in $C_j$.  Then it is easy to see
that $G_j$ acts on $A_j$ and $G$ permutes $A_j$s by projective
diffeomorphisms.  $A_j$ covers $A_j/G_j$ embedded in $\Sigma$, and
$A_j/G_j$ is a suborbifold of $\Sigma$ bounded by $c_1$.  The Euler
characteristic $\chi(A_j/G_j) = 0$ since $\chi(A_j) = 0$.  Since the
completions of components of $\Sigma - c_1$ have negative Euler
characteristic and the embedded suborbifolds $A_j/G_j$ are such
completions, we obtain a contradiction.  
Hence, the element of $C_j$ is unique.

We obtain $d_2$ in a similar manner. If $d_1$ and $d_2$ are
principal closed geodesics or geodesic full $1$-orbifolds, subsegments
of $d_1$ and $d_2$ do not bound a bigon. (To prove this, simply lift
to the universal cover which may be considered a convex domain in an
affine patch.) Thus, $d_1$ and $d_2$ meet in the least number of
points or $d_1 = d_2$.  If $d_1 = d_2$, then a power of $c_{1,i}$ is
homotopic to that of $c_{2,j}$ for some $i, j$. Since they are
disjoint simple closed curves, they bound an annulus $A$ in the closed
surface $\Sigma'$. If there are $c_{1, k}$ or $c_{2, l}$ for some $k$
and $l$ in the interior of $A$, they are essential in $A$. We find an
annulus $A'\subset A$ bounded by them and with interior disjoint from
such curves. (Possibly $A' = A$.) Now, $A'$ covers a suborbifold in
$\Sigma$ with zero Euler characteristic, which is not possible.  We
conclude that $d_1$ and $d_2$ are disjoint.
By induction, % we can show that 
there exist 
principal closed geodesics or principal 
$1$-orbifolds $d_1, \dots, d_k$ disjoint from each other
so that $c_i$ is isotopic to $d_i$ for each $i$.
\end{proof}

\subsection{Elementary Orbifolds.}\label{sec:elementary}
The following orbifolds are said to be {\em elementary}.  They are
required to be convex and have principal geodesic boundary components
and have negative Euler characteristic.  We require that no closed
geodesic is in their singular locus.  We give nicknames 
and Thurston's notations in parentheses:
\begin{enumerate} %Elementary orbifolds

\item[(P1)] A pair-of-pants.

\item[(P2)] An annulus with one cone-point of order $n$. ($A(\smc n)$)

\item[(P3)] A disk with two cone-points of order $p, q$, one of 
which is greater than $2$. ($D(\smc p, q)$) 
\item[(P4)] A sphere with three cone-points of order 
$p, q, r$ where $1/p + 1/q + 1/r < 1$. ($\SI^2(\smc p, q, r)$)

\item[(A1)] An annulus with one boundary component 
% of the underlying space 
a union of a singular segment 
and one boundary-orbifold. The other boundary component 
is a principal closed geodesic. (We call it a 
{\em $2$-pronged crown  \/} and denote it
$A(2, 2 \smc)$.)
It has two corner-reflectors of order $2$
if the boundary components are silvered.

\item[(A2)] An annulus with one boundary component of
the underlying space in a singular locus with 
one corner-reflector of order $n$, $n \geq 2$. 
(The other boundary component is a principal closed geodesic
which is the boundary of the orbifold.) (We call it 
a {\em one-pronged crown} and denote it $A(n \smc)$.) 

\item[(A3)] A disk with one singular segment and one boundary $1$-orbifold 
and a cone-point of order greater than or equal to three
($D^2(2,2 \smc n)$). %xi
\item[(A4)] A disk with one corner-reflector of order $m$ 
and one cone-point of order $n$ so that $1/2m + 1/n < 1/2$
(with no boundary orbifold). ($n \geq 3$ necessarily.)
($D^2(m \smc n)$.) %xii

\item[(D1)] A disk with three edges and three boundary
$1$-orbifolds. No two boundary $1$-orbifolds are adjacent.
(We call it a {\em hexagon} or $D^2(2,2,2,2,2,2 \smc)$.) %viii

\item[(D2)] A disk with three edges and two boundary $1$-orbifolds 
on the boundary of the underlying space. 
Two boundary $1$-orbifolds are not adjacent, and 
two edges meet in a corner-reflector of order $n$, and 
the remaining one a segment. (We called it a {\em pentagon} and denote it
by $D^2(2,2,2,2,n \smc)$.) 

\item[(D3)] A disk with two corner-reflectors of order $p$, $q$, 
one of which is greater than or equal to $3$, and one boundary $1$-orbifold.
The singular locus of the disk is a union of three edges and two 
corner-reflectors. (We call it a {\em quadrilateral} or
$D^2(2,2,p,q \smc)$.) %ix
\item[(D4)]  A disk with three corner-reflectors of order $p, q, r$ where
$1/p + 1/q + 1/r < 1$ and three edges (with no boundary orbifold).
(We call it a {\em triangle} or $D^2(p, q, r \smc)$.) 
\end{enumerate}

\begin{figure}[ht] %%
\centerline{\epsfxsize=4.5in \epsfbox{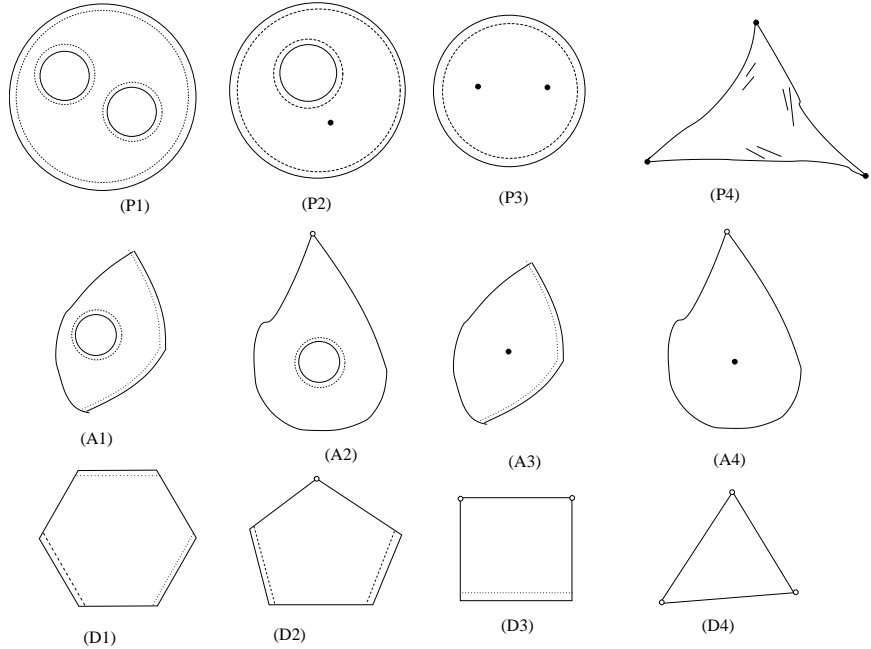}}
\caption{\label{fig:eleorb} 
The elementary orbifolds. % which we reproduced 
% from Thurston's notes. 
Arcs with dotted arcs next 
to them indicate boundary components. 
Black points indicate cone-points
and white points the corner-reflectors.}
\end{figure}
\typeout{<>}

We justify our notations. 
P1, P2, P3, and P4 type $2$-orbifolds are all obtained from a 
pair-of-pants by ``degenerating'' one, two, and three boundary 
components to cone-points respectively. (Such processes 
are realizable as deformations of hyperbolic 
structures first making the boundary components
into cusps and from cusps to cone-points
as can be accomplished by Kerckhoff's paper 
\cite{Kerck}. See also Cooper-Hodgson-Kerchkoff 
\cite{CHK:2000} in three-dimensional cases.)
A1, A2, A3, and A4 type $2$-orbifolds 
are doubly-covered by P1, P2, P3, and 
P4 type $2$-orbifolds respectively
where for each of them the order-two deck-transformation acts 
with the set of fixed points equal to an arc 
from a single boundary component
or a single cone-points to itself.
%For A3 and A4 orbifolds, the orders of 
%the cone-points of covering orbifolds have to match. 
%Here, one boundary component or one cone-point
%lifts upstairs. 
%%%%! clsrify  C. OK.
D1, D2, D3, and D4 are also doubly-covered by 
P1, P2, P3, and P4 respectively where 
for each of them the order-two deck-transformation 
acts with the set of fixed points the union of
mutually disjoint three arcs obtained by 
connected any pair of cone-points or boundary components.
%Here, each boundary component or cone-point is doubly covered
%by a boundary component or a cone-point in the covering 
%orbifold. 
(Again obtaining an A2 orbifold from an A1 orbifold 
can be realized as a deformation of hyperbolic structures.
The same can be said for other types here.)

Since a unique complete hyperbolic structure exist on 
a $3$-punctured spheres with three punctures, 
an orientation-preserving self-homeomorphism of
a pair-of-pants fixing each of the ends is isotopic to the identity.
% (from a simple argument using properly embedded arcs).  
% this vague statement should either be explained or deleted
Automorphism groups of P1, P2, P3, or P4 type orbifolds are thus
determined by the action on boundary components and cone-points.
They are realized as permutation groups of 
oriented boundary components and cone-points.  It follows
that the above elementary orbifolds are all of the quotient orbifolds
orbifolds of type P1, P2, P3, and P4.

\begin{rem}
None of these have  $1$-orbifolds decomposing them 
into unions of negative Euler characteristic 
orbifolds. % This is straightforward 
\end{rem}

%%%% let's take the next stuff out
% to see using case-by-case approach.
% But we do not need it in this paper. 
% if we don't need it, why distract the reader?
%C. I think we use it for hyperbolic stuffs. 
%C. However, we may move it?

%If we just consider them topologically, then they should be called
% {\em topological elementary orbifolds}.  
% They play a role of
% pair-of-pants decomposition of surfaces for $2$-orbifolds.
%%%% do we ever use this?
%%% it seems pretty obvious and I don't think saying this
%%% adds anything to the exposition

% Thus, a topological
% proof of showing that a $2$-orbifold of negative Euler 
% characteristic decomposes into these is easier than
% our geometric version, Theorem \ref{thm:decomp}
% since we do not need to obtain geometric realizations of 
% $1$-orbifolds used in the geometric operations.
% \end{rem}
%%%%! strengthen 10/16

% \begin{rem}\label{rem:teichmuller}
% A compact $2$-orbifold of negative Euler
% characteristic always admits a hyperbolic structure so that 
% the boundary orbifolds are geodesics. 
% This is proved in Chapter 5 of Thurston's note
% \cite{Thnote}. Hence, above elementary orbifolds (in topological
% sense) admit hyperbolic projective structures; that is, 
% there are some convex projective structures on each type of 
% elementary orbifolds. Actually, from hyperbolic
% geometry, one can easily construct them all, without
%any trouble, as we do in \S \ref{sec:teichmuller}.

% I think we've already said this several times.
% Maybe delete unless needed.

If $\Sigma$ has a convex $\rpt$-structure, and each boundary 
component of $\Sigma$ is geodesic 
and principal, 
then the decomposed orbifolds have $\rpt$-structures 
with principal geodesic boundary.

% this on the other hand is important
%

\subsection{Decomposition into elementary orbifolds} 
% Now we prove the decomposition 
% into elementary orbifolds.

\begin{thm}\label{thm:decomp}
Let $\Sigma$ be a compact convex $\rpt$-orbifold with 
$\chi(\Sigma) < 0$ and principal geodesic
boundary. Then there exists a mutually disjoint collection of 
simple closed geodesics and mirror- or cone- or mixed-type
geodesic $1$-orbifolds so that $\Sigma$ 
decomposes along their union to a union of elementary $2$-orbifolds
or such elementary $2$-orbifolds with some boundary $1$-orbifolds silvered
additionally.
\end{thm}
\begin{proof} 
We essentially follow the sketch of 
the proof of Chapter 5 of Thurston \cite{Thnote}.
We remark that
we don't really need to have silvered elementary $2$-orbifolds 
in the conclusion since we can define 
clarifying as ``decomposing" also; however, we don't 
need this fact in this paper.
In this proof, we won't distinguish between
silvered elementary annuli and elementary annuli
mainly in indicating their types.

Suppose first that $\Sigma$ has no corner-reflectors
or boundary full $1$-orbifolds; that is, singular points of $\Sigma$ 
are cone-points or in a closed geodesic of mirror points in 
the boundary of $X_\Sigma$, and the boundary components of 
$\Sigma$ are principal closed geodesics. 

Let $p_1, \dots, p_m$ be the cone-points of $\Sigma$.  Assume for the
moment that the number of cone-points $m$ is greater than or equal to
$3$, and $\chi(X_\Sigma) \leq 0$.  Let $D$ be a disk in $\Sigma$
containing all $p_1, \dots, p_m$. Then $D$ has a structure of a
$2$-orbifold of negative Euler characteristic (by equation
\eqref{eqn:riehur}).  If at least one $p_i$ has an order greater than
two, then we find a simple closed curve bounding a disk $D_1$
containing $p_i$ and another cone-point, say $p_j$. We find simple
closed curves that bounds a disk $D_2$ including $D_1$ and another
cone-point.  We find a disk $D_3$ including $D_2$ and another
cone-point and so on. Therefore, $D$ is divided into a disk $D_1$ and
annuli $A_i$ containing unique cone-points. Each of $D_1$ and $A_i$s
has negative Euler characteristic. $X_\Sigma - D$ contains no singular
points and has negative Euler characteristic.  
%
% Apply the pair-of-pants
% decomposition to $X_\Sigma - D^o$, and obtain pairs of pants.  
%
Decompose $X_\Sigma - D^o$ into pairs-of-pants.
By Lemma \ref{lem:disjointcut}, we can find a collection of principal
closed geodesics in $\Sigma$ decomposing $\Sigma$ into convex 
$\rpt$-suborbifolds of negative Euler characteristic 
of type (P1), (P2), and (P3).

%\marginpar{A substantial change here. Should shorten the first part as well?}

If all cone-points have order two, then we choose a disk $D_1$ in $D$
containing three cone-points.  Split along a cone-type $1$-orbifold
connecting the two cone-points in $D_1$.  Now find $D_2, D_3, \dots,$
as above to obtain the decomposition into elementary $2$-orbifolds of
type (P1) or (P2).

If $\chi(X_\Sigma) =1$, then $X_\Sigma$ is homeomorphic to a disk or
$\rpt$.  In the former case, $X_\Sigma$ decomposes as above.  In the latter
case, the decomposition along a one-sided simple closed curve
produces a disk with cone-points 
In particular the resulting orbifold has negative Euler characteristic.
Now apply Lemma~\ref{lem:disjointcut}. 

Suppose that $\chi(X_\Sigma) = 2$.  When $m=3$, $\Sigma$ is an
elementary $2$-orbifold of type (P4).  

Next suppose $m\geq 4$ and at least two cone-points of order greater than $2$. 
Then there exists a simple closed curve in $X_\Sigma$ 
bounding two disks containing at
least two cone-points, each containing a cone-point of order greater
than $2$, and with negative Euler characteristic. Divide the
complement of the union of these disks to annuli as above.

Next suppose $m\geq 4$ and a single cone-point of order greater
than two.  Then there exists a cone-type $1$-orbifold connecting two
cone-points of order two. Cutting along the $1$-orbifold produces
a disk with principal geodesic boundary.  Decompose the
disk into union of a disk containing the single cone-point of order
greater than two, and annuli with single cone-points.  If all
cone-points are of order two, then $m\geq 5$.  Then decompose
along two disjoint cone-type $1$-orbifolds, obtaining  an
annulus with cone-points of order two. Now proceed as above.

Suppose that $m=1$ or $2$. If $\chi(X_\Sigma) < 0$, then we can
introduce one or two annuli or M\"obius bands containing the one
singular point each and we obtain a decomposition into (P1), (P2),
(P3) type orbifolds.  If $\chi(X_\Sigma) = 0$, then an
essential simple closed curve cuts $X_\Sigma$ into an annulus.
The decomposed orbifold is connected 
and of negative Euler characteristic.  
If $m=1$ and $X_\Sigma$ is homeomorphic to $\rpt$ or $\SI^2$,
then $\chi(\Sigma)\ge 0$.
If $m=2$ and $\chi(X_\Sigma) \geq 0$, then $X_\Sigma$
is homeomorphic to $\rpt$. Split $\Sigma$ along a
simple closed curve into a disk with two cone-points, leading to 
the next case.  Finally if $X_\Sigma$ is homeomorphic to a disk, 
then $\Sigma$ is an elementary $2$-orbifold of type (P3) and $m=2$.  
(This settles case (d) below.)

We now suppose that $\Sigma$ has corner-reflectors
and/or boundary full $1$-orbifolds.
% and  some or no
% boundary full $1$-orbifolds
%% does that last phrase have any content??
Let $b_1, b_2, \dots, b_l$ be the boundary components of 
$X_\Sigma$ containing corner-reflectors or boundary 
$1$-orbifolds, and $b_{l+1}, \dots, b_k$ denote the remaining
components. 

Let $c_i$ be a simple closed curve in $X_\Sigma^o$ homotopic to $b_i$
for $i=1,\dots, l$.  The component of $\Sigma - c_i$ containing $b_i$
has negative Euler characteristic by equation \eqref{eqn:riehur2}.

Suppose that $\Sigma$ split along the union of $c_i$s has a component $C$
with nonnegative Euler characteristic. Then the underlying space of
this component is homeomorphic to either a disk, an annulus, or a
M\"obius band. If $C$ is homeomorphic to an annulus or M\"obius band, 
no singular point of $\Sigma$ lies on the component.
When $C$ is an annulus, cut along the center closed curve obtaining 
$C$ decomposed into two annuli with singular points.
They can be pasted with other components of $\Sigma$ without 
changing their Euler characteristics, obtaining a different 
decomposition of $\Sigma$.
When $C$ is a M\"obius band, cut along a one-sided closed curve obtaining 
an annulus with singular points. The annulus is then pasted
with the adjacent component of $\Sigma$ 
without changing the Euler characteristic, obtaining 
a new decomposition. If $C$ is homeomorphic to a disk, then there are 
either two cone-points of order two or a single cone-point
of arbitrary order by the Euler characteristic condition 
or no cone-points, and $X_\Sigma$ is homeomorphic to a disk. 
If there are two cone-points of order two, 
we can split $\Sigma$ along a cone-type $1$-orbifold in $C$
to obtain an annulus with corner-reflectors.
If there is a single cone-point, then $\Sigma$ 
must have been a disk with one cone-point.
If there is no cone-point in $C$, then $\Sigma$ is 
a disk. 

These facts and Lemma \ref{lem:disjointcut} imply that 
$\Sigma$ decomposes into
convex $2$-suborbifolds listed below:
%%%%! continue cleaning this up
%%%% C. rewrite????
\begin{itemize}
\item[(a)] A disk with corner-reflectors and boundary full $1$-orbifolds
and no cone-points.
\item[(b)] An annulus with corner-reflectors and boundary full
$1$-orbifolds in one boundary component of
its underlying space and no cone-points.
\item[(c)] A disk with corner-reflectors, boundary
$1$-orbifolds, and one cone-point of arbitrary order.
\item[(d)] A $2$-orbifold without corner-reflectors or boundary 
full $1$-orbifolds.
\end{itemize}

When there are boundary full $1$-orbifolds, we silver them.  This does
not change the Euler characteristic of any $2$-suborbifold containing
it by equation \eqref{eqn:riehur2}.  We mark the full $1$-orbifold to be
of {\em boundary type}, and {\em we will never cut them apart}.  After
our decomposition, we will clarify them.

Let $\Sigma'$ be one of the above (a)-(c).  
Let $m$ denote the number of corner-reflectors.
In case (a), assume $m \ge 3$.  If $m=3$, 
the orbifold $\Sigma'$ is an elementary one of type
(D4). Suppose that $m \geq 4$. From a side $A$, 
we choose one $1$-orbifold, say $a_1, \dots, a_{m-2}$, 
ending at $A$ and each of the other sides except the 
sides adjacent to $A$. We choose these to be disjoint.
If $a_i$ ends at a side was a boundary full $1$-orbifold,
we delete $a_i$. 
If any component of $\Sigma'$ split along the remaining 
$1$-orbifolds $a_i$ has Euler characteristic greater than or equal to zero, 
then we delete one of the $1$-orbifold adjacent to the component. 
Continuing in this manner, we obtain a collection 
$a_i$ so that if $\Sigma'$ is split along them, each
of the component is an orbifold of type (D1), (D2), or (D3).
Finally, clarifying the boundary types yields the desired decomposition.

In case (b), the underlying space of $\Sigma'$ is a disk
minus an open disk in its interior.  As before let $p_1, p_2,
\dots, p_m$ denote the corner-reflectors and $e_1, e_2, \dots, e_m$
edges.  If $m=1$, then $\Sigma'$ is an elementary $2$-orbifold of type
(A2). 
%If $m=2$, and the corner-reflectors are of order two, then
%$\Sigma'$ is an elementary $2$-orbifold of type (A1) with boundary
%silvered.  

Suppose $m \geq 2$. From an edge $\kappa$ which is not boundary type, 
we choose $1$-orbifolds $a_1, \dots, a_m$ ending at the sides
$e_1, \dots, e_m$ so that the two endpoints of $a_m$ are in $\kappa$
and other $a_i$ has endpoints in distinct edges. 
We eliminate $a_i$s that ends in the edges adjacent to $\kappa$
and ones of boundary type. If we split along $a_i$s
and obtain a component of Euler characteristic greater than 
zero, then we eliminate one of the adjacent $a_i$s. 
As in case (a), by continuing in this manner, 
we decompose $\Sigma'$ into an elementary 
orbifold of type (A2) and ones of type (D1), (D2), or (D3). 

%If $m=2$, and at least one corner-reflector has order $\geq
%3$, then there exists a $1$-orbifold with endpoints in an edge, along
%which $\Sigma'$ decomposes into a quadrilateral and an elementary
%$2$-orbifold of type (A1).

%Suppose $m \geq 2$.  Choose an edge $k$ which is not a boundary type
%full $1$-orbifold.  Choose a $1$-orbifold $l$ with endpoints in $k$ so
%that a segment in $k$ and itself bound an annulus including the
%boundary component without corner-reflectors. (We require the
%topological interior of $l$ to lie in the topological interior of
%$\Sigma'$.)  Then $l$ decomposes $\Sigma'$ to an elementary
%$2$-orbifold of type (A1) and one whose underlying space is
%homeomorphic to a disk.  For a disk-type one, we silver $l$ to become
%a singular-mirror-type-$1$-orbifold.  Applying above arguments again,
%we obtain a desired decomposition of $\Sigma'$. (We clarify copies of
%$l$ back at the end.)

In case (c), suppose that the cone-point has order greater than
two. If the number of corner-reflectors is one, then $\Sigma'$ is an
elementary $2$-orbifold of type (A4). As above, we find 
$1$-orbifolds $a_1, \dots, a_m$ and do similar constructions. 

%If it is two, and one of the
%corner-reflectors has an order higher than two, then we can again
%decompose $\Sigma'$ to a quadrilateral and an elementary $2$-orbifold
%of type (A3).  If it is two, and so are the orders of the
%corner-reflectors, then $\Sigma'$ is elementary one of type (A3).
%Assume that the number is greater than or equal to three.  We find a
%$1$-orbifold $l$ with endpoints in an edge along which $\Sigma'$
%decomposes to a disk-type orbifold and an elementary $2$-orbifold of
%type (A3).  Again, we silver $l$ and decompose the disk-type orbifold
%as above, and finally clarify copies of $l$ to get the decomposition
%into quadrilaterals, pentagons, and type-(A3)-elementary-orbifolds.

Suppose that the cone-point has order two. Draw
a mixed-type $1$-orbifold $\tau$ from the cone-point to a mirror point
in the nonboundary type full $1$-orbifold. The 
decomposition results in a disk. Silver the $1$-orbifold 
which folds to $\tau$. Decompose the disk and clarify it back.
\end{proof}

%%4/20 8:45
\section{The Teichm\"uller spaces of elementary orbifolds}
\label{sec:teichmuller}

This section contains a proof of part of Thurston's 
Theorem \ref{thm:Thur} for elementary $2$-orbifolds . 

% (see Remark \ref{rem:teichmuller}). 
We will need in the next section that the space of
convex $\rpt$-structures with principal boundary
on each elementary $2$-orbifold is nonempty. 
Since hyperbolic
projective structures on elementary $2$-orbifolds are certainly
convex with principal boundary, 
the following proposition implies this nonemptyness. 
 Recall that the Hilbert metric of a hyperbolic $\rpt$-structure
is the hyperbolic metric.
%By ``length'' in this section, we mean the length in
% hyperbolic metric or the Hilbert metric of the conic,
% which are the same metric.

\begin{prop}\label{prop:hypele}
Let $\Sigma$ be a compact $2$-orbifold with empty boundary
and negative Euler characteristic diffeomorphic to
an elementary $2$-orbifold. 
Then the deformation 
space $\mathcal{T}(\Sigma)$ of hyperbolic $\rpt$-structures 
on $\Sigma$ is homeomorphic to 
a cell of dimension $-3\chi(X_\Sigma) + 2k + l + 2n$ 
where $X_\Sigma$ is the underlying space and 
$k$ is the number of cone-points, $l$ is 
the number of corner-reflectors, and $n$
is the number of boundary full $1$-orbifolds.
\end{prop} 

% Here, $\mathcal{T}(\Sigma)$ is a subspace of $\rpt(\Sigma)$
% representable by hyperbolic projective structures. 
% It is easy to see that $\mathcal{T}(\Sigma)$ is identifiable
% with the Teichm\"uller space in the ordinary sense. 
%%%% i think we already said this
%C. delete

The strategy of the proof is to show that
for each elementary $2$-orbifold $S$, 
$\mathcal{T}(S)$ is homeomorphic to $\mathcal{T}(\partial S)$, 
where $\mathcal{T}(\partial S)$ is the product of $\bR^+$ for
each component of $\partial S$ corresponding to the 
hyperbolic-metric lengths
of components of $\partial S$. Then for hyperbolic
structures, to obtain a bigger orbifold in the 
sense of this paper, we simply need to match lengths
of boundary components. 

A generalized triangle in the hyperbolic plane is one of following:
\begin{itemize}
\item[(a)] A hexagon: a disk bounded by six geodesic sides meeting in
right angles labeled $A, \beta, C, \alpha, B, \gamma$.
\item[(b)] A pentagon: a disk bounded by five geodesic sides 
labeled $A, B, \alpha, C, \beta$ where $A$ and $B$ meet in
an angle $\gamma$, and the rest of the angles are right angles. 
\item[(c)] A quadrilateral: a disk bounded by four geodesic sides 
labeled $A, C, B, \gamma$ where $A$ and $C$ meet in an angle $\beta$,
$C$ and $B$ meet in an angle $\alpha$ and the two remaining
angles are right angles. 
\item[(d)] A triangle: a disk bounded by three geodesic sides labeled
$A, B, C$ where  $A$ and $B$ meet in an angle $\gamma$ and 
$B$ and $C$ meet in an angle $\alpha$ and $C$ and $A$ meet in
angle $\beta$. 
\end{itemize}

\begin{figure}[ht] %%
\centerline{\epsfxsize=3.5in \epsfbox{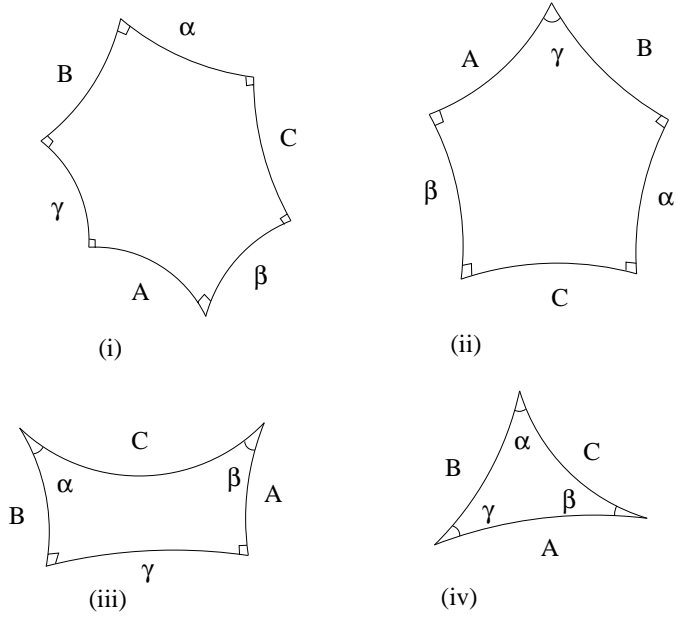}}
\caption{\label{fig:gentri} Generalized triangles 
and their labels.}
\end{figure}
\typeout{<>}

%%%%%%%%%%%%%%%%%%%%%%%%%%%%%%%%%%%%%%%%%%%%%%%%%%%%%%%%%%%%%%%%%%%%%%%%%%%%%%
% the indices in the figures don't match the notation in Lemma 5.2!!!
%%%%%%%%%%%%%%%%%%%%%%%%%%%%%%%%%%%%%%%%%%%%%%%%%%%%%%%%%%%%%%%%%%%%%%%%%%%%%%

\begin{lem}\label{lem:polygons} 
For generalized triangles in the hyperbolic plane, the following
equations hold. 
\begin{eqnarray}
{\mathrm{(a)}} & \cosh C 
&= \frac{\cosh \alpha \cosh \beta + \cosh \gamma}{\sinh \alpha \sinh \beta}
\nonumber \\
{\mathrm{(b)}} & \cosh C 
&= \frac{\cosh \alpha \cosh \beta + \cos \gamma}{\sinh \alpha \sinh \beta}
\nonumber \\
{\mathrm{(c)}} & \sinh A 
&= \frac{\cosh \gamma \cos \beta + \cos \alpha}{\sinh \beta \sin \gamma}
\nonumber \\
{\mathrm{(d)}} & \cosh C 
&= \frac{\cos \alpha \cos \beta + \cos \gamma}{\sin \alpha \sin \beta}
\end{eqnarray}
In {\rm (a),} $(\alpha, \beta, \gamma)$ can be any triple of positive 
real numbers.
In {\rm (b),} $(\alpha, \beta)$ can be any pair of
positive numbers and $\gamma$ a real number in $(0, \pi)$
In {\rm (c),} $(\alpha, \beta)$ can be any pair of
positive real numbers in $(0, \pi)$ satisfying $\alpha + \beta < \pi$.
In {\rm (d),} $(\alpha, \beta, \gamma)$ 
can be any triple of real numbers in $(0, \pi)$
satisfying $\alpha + \beta + \gamma < \pi$.
\end{lem} 
\begin{proof}
For a proof, see Coxeter \cite{Cox:99} or Ratcliffe \cite{RT:94}.
\end{proof}

The following lemma implies Proposition \ref{prop:hypele}
for elementary $2$-orbifolds of type (D1), (D2), (D3), and (D4). 
\begin{lem}\label{lem:hypele2} 
Let $P$ be an orbifold whose underlying space is a disk constructed
as follows.
Silver edges labeled by the capital letters $A, B, C$.
Assign to each vertex an angle of the form $\pi/n$ (where
{\rm (}$n>1$ is an integer{\rm ),} 
for which it is a corner-reflector of that angle.
Each edge labeled by Greek letters $\alpha,\beta,\gamma$ is 
a boundary full $1$-orbifold.
Then in cases {\rm (a), (b), (c), (d)}
$\mathcal{F}: \mathcal{T}(P) \ra \mathcal{T}(\partial P)$ for
each of the above orbifolds $P$ is a homeomorphism{\rm ;}
that is, $\mathcal{T}(P)$ is homeomorphic to a cell of dimension $3$,
$2$, $1$, or $0$ respectively.
\end{lem}
\begin{proof} 
We see that the edge lengths or angle measures of
$\alpha, \beta, \gamma$ completely
determine the unique disks. 
The edges labeled by Greek letters can be made 
arbitrarily large or small after the angles labeled by
Greek letters are assigned 
as the reader can easily verify. 
By the above formulas, $\mathcal{F}$ is a homeomorphisms. 
\end{proof}

\begin{lem}\label{lem:hypele3} 
Let $S$ be an elementary $2$-orbifold of 
type {\rm (A1), (A2), (A3),} or {\rm (A4).} 
Then $\mathcal{F}:\mathcal{T}(S) \ra \mathcal{T}(\partial S)$ is
a homeomorphism. Thus, $\mathcal{T}(S)$ is a cell of 
dimension $2, 1, 1,$ or $0$ when $S$ is of type
{\rm (A1), (A2), (A3)} or {\rm (A4)} respectively.
%In case {\rm (A4),} $\mathcal{T}(S)$ is a single point.
\end{lem}
\begin{proof}
In case (A1), find the shortest segment $s$ from 
one boundary component to the other. 
Cutting along it, we obtain a hexagon, where
the boundary is cut into three alternating sides of
the hexagon $\alpha, \beta, \gamma$. 
Let $\alpha$ and $\beta$ be from the boundary 
full $1$-orbifold. Let the length of $\alpha$ and 
$\beta$ be equal. By symmetry, the lengths of 
two sides corresponding to $s$ become equal.
Hence, we can glue back such a hexagon to
obtain an elementary $2$-orbifold of type (A1) always.
Thus, we see that given $2\alpha$ and $\gamma$,
we can obtain an elementary $2$-orbifold.
Since $2\alpha$ and $\beta$ are lengths of
the boundary components, we have shown that
$\mathcal{F}$ is a homeomorphism.

In case (A2), we silver the boundary component
temporarily, find a shortest segment $s$ from
the mirror edge to the boundary component. 
Then cutting along $s$, we obtain a pentagon. 
Let $\alpha$ and $\beta$ be the boundary
full $1$-orbifolds of the pentagon. 
Letting the lengths of $\alpha$ and $\beta$
equal, we can always glue back to obtain
an elementary $2$-orbifold of type (A2). Since we can 
change the length of $C$ arbitrarily
by changing the common length of $\alpha$ and $\beta$,
we see that $\mathcal{F}$ is a homeomorphism.

In case (A3), we draw a shortest segment $s$
from the cone-point to the boundary $1$-orbifold.
We obtain a pentagon where $s$ corresponds 
to edges labeled $A$ and $B$ above. If the lengths
of $A$ and $B$ are equal, we can glue back to
obtain an elementary $2$-orbifold of type (A3). 
Here, $\alpha$ and $\beta$ are from the boundary
$1$-orbifold by cutting, and their lengths are the same. 
If the lengths of $\alpha$ and $\beta$ are the same,
then those of $A$ and $B$ are equal, and we can glue back.
Since $2\alpha$ is the length of the boundary
$1$-orbifold, we see that $\mathcal{F}$ is a homeomorphism. 

In (A4), we draw a shortest segment $s$ from
the cone-point to the mirror edge. Then 
we obtain a quadrilateral with angles $\pi/p$, $\pi/2$,
$2\pi/q$, and $\pi/2$ where $p$ is the order of the corner-reflector
and $q$ is that of the cone-point. Such a quadrilateral 
is unique, and we can glue back to obtain
an elementary $2$-orbifold of type (A4) always. 
Thus the Teichm\"uller space is a single point.
\end{proof}

Finally, if $S$ is an elementary $2$-orbifold of type (P1),(P2),(P3), or
(P4), then there exists an order-two self-isometry so 
that the quotient $2$-orbifolds are of type (D1), (D2), (D3), or (D4)
respectively. 
This can be seen by choosing shortest 
segments connecting any pair of 
the boundary components or cone-points. 
Proposition \ref{prop:hypele} in these cases follows from 
Lemma \ref{lem:hypele2}. To conclude, 
Lemmas \ref{lem:hypele2} and \ref{lem:hypele3}
imply Proposition \ref{prop:hypele}.
(In the final part of \S\ref{sec:deformelement}, we will prove 
Theorem \ref{thm:Thur}.) 

%%5/1 10:30

\section{The deformation spaces of elementary $2$-orbifolds \\
and the proof of Theorem A.}
\label{sec:deformelement}

This section concludes the proof of Theorem~A,
that the deformation spaces of 
convex $\rpt$-structures on elementary $2$-orbifolds
are homeomorphic to cells.
The inductive proof begins from elementary 
orbifolds: We show that the map from the deformation 
space to the deformation space of the boundary is 
a fibration with base and fiber homeomorphic to cells. 
After this, when we build bigger orbifolds, this facts 
automatically hold. 
%Thus, the main result follows. 
%C. delete

The method we use for elementary orbifolds of type (P1)-(P4)
is from \cite{Gconv:90} for a pair-of-pants. 
In fact, we need very small changes for our purposes.
That is, we first show that the $2$-orbifold of 
type (P1)-(P4) is built up from 
two triangles. We study the quadruples of triangles in order 
to understand how they are glued. By converting the geometric
conditions that they assemble to convex elementary orbifolds
into algebraic relations and solving the relations, 
we show that the deformation
spaces are described by cells. 
For some annular elementary $2$-orbifolds (A1)-(A4), 
we need Steiner's theorem defining conics
as the set of intersection points of pencils of lines 
through two given points related by a projectivity
(see Chapters 6 and 7 of Coxeter \cite{Cox:60}).
%%%! what is Steiner's theorem? Reference?
%%%Reference given below.
To study elementary $2$-orbifolds with corner-reflectors
(D1)-(D4), we generalize the methods in \cite{G:77}. 
%However, some parts of this work are new. 
Finally, Theorem A is proved by noticing that the above fibration 
property can be inductively proved as we build up 
$2$-orbifolds from elementary $2$-orbifolds.
Finally, we prove Theorem \ref{thm:Thur}
(originally due to Thurston)
for the sake of completeness of the paper.

%Using the results on elementary $2$-orbifolds, 
%we prove the final result we seek, 
The result proving Theorem A is as follows:
\begin{thm}\label{thm:finaldim}
Let $\Sigma$ be a compact $2$-orbifold with negative Euler 
characteristic. Then the deformation space of convex 
$\rpt$-structures $\mathcal{C}(\Sigma)$ is homeomorphic to
a cell of dimension 
\begin{equation}\label{eqn:finaldim}
d(\Sigma) := -8\chi(X_\Sigma) + (6k_c -2 b_c) + (3 k_r - b_r) + 4n_0 
\end{equation}
where $k_c$ is the number of cone-points, 
$k_r$ the number of corner-reflectors, 
$b_c$ the number of cone-points of order two, 
$b_r$ the number of corner-reflectors of order two,
and $n_o$ is the number of boundary full $1$-orbifolds.
\end{thm}

% We determine the topology
% and the dimensions of deformation spaces of convex
% $\rpt$-structures on elementary $2$-orbifolds with
% boundary orbifold invariants fixed. 
% Then we obtain
% the total dimensions by varying the invariants and the gluing
% freedom based on conjugating maps. The dimension
% of the space of invariants of the boundary is obviously
% the sum of dimensions of those of the boundary components.
%C. outline needed somewhere?

\subsection{The Fibration Property}

\begin{defn}\label{defn:dimension}
We say that the deformation space of
$\Sigma$ satisfies the {\em fibration property} 
if:
\begin{itemize}
\item[(i)] the deformation space of convex 
$\rpt$-structures $\mathcal{C}(\Sigma)$ is homeomorphic to
a cell of dimension $d(\Sigma)$ as in
equation \eqref{eqn:finaldim}, and
\item[(ii)] there exists a principal fibration
${\mathcal{F}}: \mathcal{C}(\Sigma) \ra \mathcal{C}(\partial \Sigma)$ 
with the action of a cell of dimension 
$\dim \mathcal{C}(\Sigma) - \dim \mathcal{C}(\partial \Sigma)$,
equal to
\[-8\chi(X_\Sigma) + (6k_c - 2b_c) +  (3 k_r - b_r) + 3n_0
- 2b_\Sigma = d(\Sigma) - n_0 - 2 b_\Sigma\] 
where $b_\Sigma$ is the number of boundary components of
$\Sigma$ homeomorphic to $S^1$.
\end{itemize}
\end{defn}

We will prove first the following:
\begin{prop}\label{prop:elemdim}
Let $\Sigma$ be an elementary $2$-orbifold of negative Euler 
characteristic.
Then the deformation space of $\Sigma$ has the fibration property.
\end{prop}

%\begin{cor}\label{cor:elemdim}
%Suppose that $\Sigma$ is obtained from an elementary $2$-orbifold 
%$\Sigma'$ by silvering some boundary
%component $1$-orbifolds. Then $\Sigma$ satisfies the fibration property.
%\end{cor}
%\begin{proof}
%Theorem \ref{thm:silvering} implies that
%$\mathcal{C}(\Sigma)$ is homeomorphic to $\mathcal{C}(\Sigma')$.
%We check the dimension formula:
%Silvering a principal boundary curve does not change 
%any of the numbers $k_c, b_c, k_r, b_r, n_o$.
%Silvering a principal full $1$-orbifold 
%$b_r$ increases by two and $k_r$ by two and $n_o$ decreases 
%by $1$. Thus, $d(\Sigma) = d(\Sigma')$. 
%The fibration property (i) follows.
%
%The map
%$\mathcal{F}:\mathcal{C}(\Sigma) \ra \mathcal{C}(\partial \Sigma)$ 
%is obtained from the fibration
%${\mathcal{F}}^{\prime}:
%\mathcal{C}(\Sigma') \ra \mathcal{C}(\partial \Sigma')$ 
%post-composed with the projection 
%$\mathcal{C}(\partial \Sigma') \ra \mathcal{C}(\partial \Sigma)$ 
%by forgetting silvered boundary components.
%Thus $\mathcal{F}$ is a principal fibration. 
%(ii) follows by checking the dimensions. 
%\end{proof}

\subsection{Simple geodesics on $2$-orbifolds}

\begin{prop}\label{prop:imbarc} 
Let $l$ be an geodesic arc in an orientable 
convex $2$-orbifold of negative Euler characteristic 
with boundary points in the boundary of the good orbifold
or an interior point. Suppose that $l$ is homotopic
to a simple arc by a homotopy fixing the boundary
points in the boundary of the orbifold or fixing the point
if it is an interior point. 
Then $l$ is simple.  
\end{prop}
\begin{proof} 
The universal cover of the orbifold is a convex domain and
contains no bigons or monogons with geodesic boundary. 
Any nonsimple arc which is homotopic to a simple 
arc must have some bigon or monogon in the universal
cover. Thus $l$ is simple.
\end{proof}

We say that a sequence $l_i$ of segments or line in a convex domain
{\em converges to a line $l_\infty$\/} in 
$\Omega$, $l_i \cap K$ converges to $l_\infty \cap K$ in the Hausdorff sense
for every compact subset $K$. of

%%% we should have a notation for the regular set of an orbifold,
%%% maybe S_reg ?
% Do you need? 
% 03/03 8:55

A {\em lamination} in a convex $2$-orbifold $S$ is a subset of the
regular set of $S$ such that for each
coordinate neighborhood $U$ not meeting the singular points meets the
set in a closed subset that is a union of disjoint lines that pass
through the open set completely.  (We allow half-open lines 
as long as they are complete in the neighborhood.)  A {\em leaf} of a
lamination is an arcwise connected subset.

A lamination is {\em finite} if each open set meets finitely many
leaves. We say that an end of an infinite lamination 
{\em winds around} a simple closed curve, 
if a half-infinite arc corresponding to the end 
is so that all of its accumulation points lie in the curve.

\begin{lem}\label{lem:twotriangle}
Suppose that $S$ is a convex $2$-orbifold which either 
is diffeomorphic to a pair-of-pants, an annulus
with a cone-point, a disk with two cone-points,
or a sphere with three cone-points.
We choose orientations on the boundary components. 
Then there exists a geodesic lamination $l$ with three leaves 
so that $S$ with the boundary components 
and the cone-points and the lamination $l$
removed is a union of two convex 
triangles. Each leaf of $l$ has two ends, 
each of which either ends in a cone-point or 
winds around a boundary component of $S$ 
infinitely often. We can choose $l$ so
that the direction of the winding follows 
the arbitrarily chosen orientation.
\end{lem}
\begin{proof}
The universal covering of $S$ is a convex domain with a projectively
invariant Hilbert metric. This metric induces a Finsler metric
on $S$ with a corresponding notion of arclength.

% Let us define a Hilbert metric in the interior of $S$ by first
% defining the Hilbert metric on the convex domain which is the
% universal cover of $S$ and inducing one on the interior of $S$. Since
% this is a Finsler metric, one can define a notion of arclength.
% C. delete

Define $S_\eps$ for each $\eps >0$ to be the compact subset of $S$
obtained by removing from $S$ the union of convex open neighborhoods
of the boundary components and the cone-points so that $S_\eps$ is
homeomorphic to a pair-of-pants.  We assume that as $\eps \ra 0$,
$S_\eps$ is strictly increasing and contains any interior nonsingular
point of $S$ eventually.

First, there exists a topological lamination $\lambda$ with this
property. We can approximate it by a finite-length three-leaf
topological lamination $\lambda_t$ with six (distinct) endpoints
either in the boundary components of $S$ or the cone-points.  We
assume that for each $\eps > 0$, as $t \ra \infty$, $S_\eps \cap
\lambda_t$ approximates $S_\eps \cap \lambda$ very closely in
$C^1$-sense. Moreover, we assume that the endpoints of $\lambda_t$ 
winds around $\partial S$ along the chosen orientation infinitely as
$t\ra \infty$.

Since $S$ is convex, 
there exists a geodesic lamination $\hat \lambda_t$ homotopic to
$\lambda_t$ by a homotopy fixing the endpoints.  We may choose
$S_\eps$ for each $\eps > 0$ so that $\hat \lambda_t$ meets it always
in a union of six connected arcs near the six end points of $\hat
\lambda_t$.  Thus $\hat \lambda_t \cap S_\eps$ is a finite-length
finite lamination with endpoints in $\partial S_\eps$ by Lemma
\ref{lem:ess}.  Since $\hat \lambda_t \cap S_\eps$ has six endpoints,
it has three leaves.

The length of $\hat \lambda_t \cap S_\eps$ is bounded for fixed $\eps$
as $t\ra \infty$. Otherwise the Hausdorff limit set is a lamination in $S_\eps$
containing an infinite leaf. This contradicts Lemma \ref{lem:ess}.

Thus for each $\eps$, we can choose a subsequence so that $\hat
\lambda_{t_i} \cap S_\eps$ converges to a three-leaf finite-length
finite lamination $l_\eps$. By using a diagonal subsequence argument,
we obtain a sequence $\hat \lambda_{t_j}$ so that $\hat \lambda_{t_j}
\cap S_\eps$ converges to such a lamination for each $\eps >0$.

The direction of the winding is the same as the topological lamination
by Lemma \ref{lem:repelfix}.

The limit is now a geodesic lamination with three leaves. 
The intersection of the geodesic lamination with the compact set $S_\eps$ has
fixed topological type. Therefore the complement in $S$ 
of the union of the boundary, the lamination,
and the cone-points is precisely the union of two disjoint convex open
disks, each of whose boundary are leaves of the lamination.

Lifting to the universal covering convex domain, 
such disks develop to convex open triangles.
\end{proof}

%\begin{lem}\label{lem:geomreal}
%Let $l$ be a topological lamination in $S$ consisting of three finite
%length arcs.  Then there exists a self-homeomorphism $S \longrightarrow S$
%sending $l$ to a lamination with geodesic leaves.
%\end{lem}
% \begin{proof}
% 
% self-homeomorphism exists follows from the uniqueness of laminations
% with three finite length leaves up to isotopy.
% 
% 
% Let $\tilde l$ be the inverse images in the universal covering $\tilde S$ of
% $S$.  Each leaf of $\tilde l$ is
% a finite length arc with well-defined endpoints.  Replace the
% leaves by geodesics with the same endpoints.  Since two geodesics in a
% convex domain $\tilde S$, either agree, meet exactly at one point, or
% are disjoint, the geodesic leaves are mutually disjoint. The union of
% geodesic leaves descends to a lamination. 
% \end{proof}
%  I don't think is really a proof.
%C. leave out proof? 

\begin{lem}\label{lem:ess}
Suppose that $S$ is as above.  Then every leaf $L$ of a geodesic
lamination in $S$ ends in a cone of point of $S$ or winds around a
boundary component, and meets $S_\eps$ in a finite-length finite
lamination.
\end{lem}
\begin{proof}
% this next phrase is vague: which ``topological considerations''
% From topological considerations, we know that an 
An end of any embedded
infinite arc in $S$ which lifts to a properly embedded arc in a universal
cover winds around a simple closed curve parallel to a boundary
component or a simple closed curve around a cone-point or end at a
cone-point.  (That is, its limit points comprise one of these.)
%(See Thurston \cite{Thnote} for the case 
%of hyperbolic metric spaces.) % but this is a purely topological result.)  
%C. OK. leave out
Thus, $L$ winds around a closed geodesic parallel
to a boundary component. Since $S$ is convex and with principal
boundary and $S$ can be covered by a convex surface with principal
boundary, the closed curve is a boundary component (see \cite{cdcr2}).
\end{proof}

\begin{figure}[ht] %%
\centerline{\epsfxsize=3.5in \epsfbox{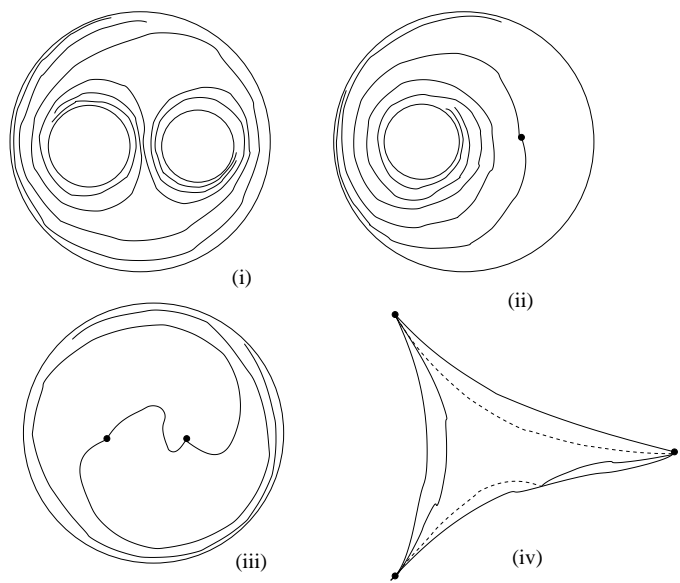}}
\caption{\label{fig:laminations} The finite laminations
in elementary $2$-orbifolds of type (P1)-(P4).}
\end{figure}
\typeout{<>}

\begin{lem}\label{lem:repelfix}
A leaf $l$ of a geodesic lamination in $S$ winds around a boundary
component $\gamma$ in a given direction if and only if $l$ lifts to a curve in 
$\tilde S = \Omega$ ending in an attracting fixed
point of the holonomy of $\gamma$.
\end{lem}
% \begin{proof}
% Let $\tilde l$ be a lift of $l$, and 
% $\alpha$ the oriented boundary component of $S$.
% There exists an element $g$ of the holonomy group of 
% $S$ so that $g$ acts on the lift $\tilde \alpha$ of
% $\alpha$ in the oriented direction. 
% Since $S$ has principal boundary,
% a triangle with vertices the attracting, repelling, 
% and saddle-type fixed points of $g$ contains 
% $\Omega$ so that the attracting fixed point and 
% the repelling fixed point are the endpoints of 
% $\tilde \alpha$. The considerations of the images
% $g^n(l)$ for $n \geq 0$ proves this lemma. 
% \end{proof}
% C. no need for a proof?

\subsection{Proof of Proposition \ref{prop:elemdim}}
We begin the case by case proof of 
Proposition \ref{prop:elemdim}. 
Case (P1) was already treated in \cite{Gconv:90}. 

\subsection{Annuli with one cone-point (P2)}
We will decompose the elementary orbifold of type (P2) 
into two triangles following \cite{Gconv:90}.
We make a one-to-one correspondence of 
$\mathcal{C}'(P)$ with configuration spaces of 
four adjacent triangles in the projective plane. 
Since $\mathcal{C}(P)$ is the quotient of 
$\mathcal{C}'(P)$ by $\PGL(3, \bR)$, 
we show that the configuration space quotient by
$\PGL(3, \bR)$ is the configuration space with 
a triangle with standard vertices quotient by 
a group of diagonal matrices. Finally, we show that 
the latter configuration space fibers over 
the boundary invariants by solving algebraic equations. 

First, we will reduce an elementary orbifold of type (P2) to 
a configuration:
Let $P$ be an annulus with a cone-point $c$ of order $n$, $n \geq 2$ 
and two boundary components $a, b$. The fundamental group of $P$ 
has presentation 
\[ \pi = \langle A, B, C | C^n = 1, ABC = 1 \rangle \]
where $A$ and $B$ are loops around two boundary components 
with boundary orientation and $C$ loops around the cone-point.
We find simple arcs $e_b$ from $c$ to $a$ and $e_b$ from $c$ to $a$.
We also find a disjoint simple arc $e_c$ from $a$ to $b$ avoiding $c$.
We will spiral $e_a, e_b,$ and $e_c$ positively in the boundary 
components $a$ and $b$, with respect to the boundary orientation, to create
geodesic laminations with leaves $l_a, l_b, l_c$
by using above results on realizations.

\begin{figure}[ht] %%
\centerline{\epsfxsize=3.8in \epsfbox{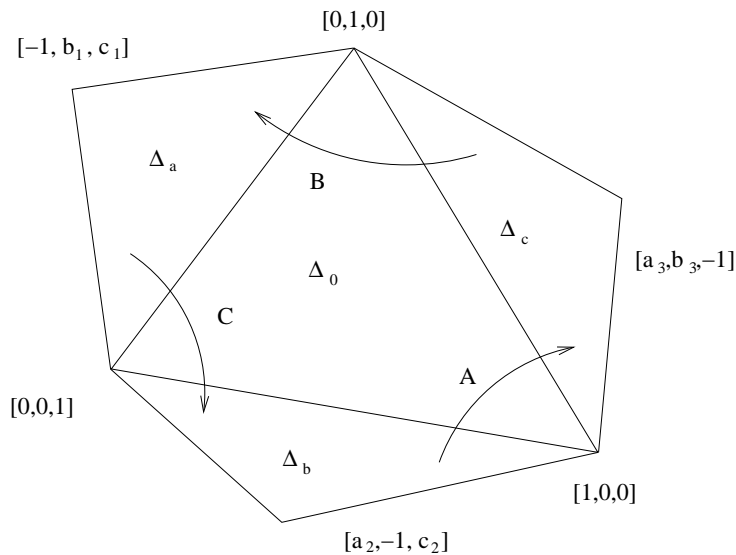}}
\caption{\label{fig:triangles} 
The four triangles needed to understand the deformation
space of an annulus with one cone-point.}
\end{figure}

We denote the geodesics by same notations.  
The complement of the lamination in $P$ with boundary
and the cone-point removed is 
a union of two open triangles. Let us denote them by $T_0$ and $T_1$.
Choose a base point $p \in T_0$, and a point $\tilde p$ 
in the inverse image in the universal $\tilde P$ of $P$. Then 
the components of inverse images of $T_0$ and $T_1$ develop
in $\rpt$ to disjoint triangles. Let $\tilde T_0$ be 
the triangle containing $\tilde p$ and let $T_a$, $T_b$, and 
$T_c$ be the components of the inverse images of $T_1$ 
adjacent to $\tilde T_0$ along the lifts of $l_a, l_b,$ 
and $l_c$ respectively. 
Therefore developing these triangles 
yields four triangles $\tri_0, \tri_a, \tri_b, \tri_c$ in $\rpt$
and collineations $A, B, C$ satisfying: 
% the following three conditions:
\begin{itemize}
\item[(i)] $\clo(\tri_a)$, $\clo(\tri_b)$, $\clo(\tri_c)$ meet 
$\clo(\tri_0)$ in three edges of $\clo(\tri_0)$. 
\item[(ii)] The union 
\[\clo(\tri_a)\cup \clo(\tri_b) \cup \clo(\tri_c) \cup \clo(\tri_0)\] 
is an embedded polygon with six or five vertices. 
(It is a convex hexagon if $n \geq 4$, and a pentagon may occurs
when $n=3$ but only rarely.)
(More precisely, any two of $\clo(\tri_a), \clo(\tri_b),$ and 
$\clo(\tri_c)$ meet exactly in a singleton as shown in
Figure \ref{fig:triangles}.)
\item[(iii)] $ABC = \idnt$, and $A(\clo(\tri_b)) = \clo(\tri_c)$, 
$B(\clo(\tri_c)) = \clo(\tri_a)$, 
and $C(\clo(\tri_a)) = \clo(\tri_b)$.
\item[(iv)] $A$ and $B$ are hyperbolic and the vertex where 
$\clo(\tri_b)$ and $\clo(\tri_c)$ meet is the repelling fixed point of $A$, 
and the vertex where $\clo(\tri_c)$ and $\clo(\tri_a)$ meet is one of $B$. 
$C$ is a rotation of order $n$ with the isolated fixed point 
the vertex at $\clo(\tri_a) \cap \clo(\tri_b)$. 
\end{itemize}

% We give a proof of (ii) only since others are obvious: 
Only the proof of (ii) is not immediate.

The complement
of the union of lines containing the three sides of $\clo(\tri_0)$ is
a union of three disjoint open triangles. Denote their closures by
$S_0, S_a, S_b, S_c$ where $S_0 = \clo(\tri_0)$.  The existence of the
actions by $A, B, C$ implies the union of $\clo(\tri_0)$ with any one of
$\clo(\tri_a), \clo(\tri_b), \clo(\tri_c)$ is a convex
quadrilateral-(*).  Thus $\clo(\tri_a), \clo(\tri_b),$
and $\clo(\tri_c)$ lie in $S_a, S_b,S_c$ respectively.  

%%% Is this an important Definition? I think ``patch''
%%% is used later in this context.
%%% previous occurrences of this term refer to ``affine patches'' 
%%% projective space
%C. I think it is important

By a patch of these triangles, we mean the disk obtained
by identifying the sides as specified by the universal cover, and it
is a disk with corners with an $\rpt$-structure.  If $n \geq 4$, then
the interior angles of the patch of $\clo(\tri_0), \clo(\tri_a),
\clo(\tri_b),$ and $\clo(\tri_c)$ are always less than $\pi$.
Hence the patch is an embedded hexagon.

Let $p_1, p_2,$ and $p_3$ denote the vertices of $\tri_0$ so that the
triangles $\clo(\tri_a)$ and $\clo(\tri_0)$ meet in a segment with
vertices $p_2$ and $p_3$, the triangles $\clo(\tri_b)$ and
$\clo(\tri_0)$ in that with $p_1$ and $p_3$, and the triangles
$\clo(\tri_c)$ and $\clo(\tri_0)$ in that with $p_1$ and $p_2$. 
(See figure \ref{fig:triangles}.)

Suppose that $n = 2$ or $3$.  If $n=3$, then the vertices of
$\clo(\tri_a), \clo(\tri_b),$ and $\clo(\tri_c)$ other than the
vertices of $\clo(\tri_0)$ are in the open triangles $S_a^o, S_b^o,
S_c^o$ as $C$ has order three. Thus, the embeddedness of the patch
follows from this fact. (Notice that the angle at $p_3$ could be $\pi$
when the vertex of $\tri_a$, $p_3$, and the vertex of $\tri_b$ are
collinear, which could happen.)

If $n=2$, then $C$ is of order two and this forces that the vertex of
$\clo(\tri_a)$ lies on the line containing $p_1$ and $p_3$, and the
vertex of $\clo(\tri_b)$ lies on the line containing $p_2$ and $p_3$.
Since $A$ and $B$ are hyperbolic, we see that the angle of the patch
is less than $\pi$ at $p_1$ and $p_2$. This implies that the patch
corresponds to an embedded hexagon with a concave vertex at $p_3$.
(See figure \ref{fig:triangles2}.)

\typeout{<>}

Let ${\mathcal{C}}'$ denote the space of all $7$-tuples
$(\tri_0, \tri_a, \tri_b, \tri_c, A, B, C)$ 
satisfying conditions (i)-(iv) with 
topology induced from the product 
\begin{equation*}
\big(\rpt\big)^{12} \times  \big(\PGL(3, \bR)\big)^3.
\end{equation*}

We showed that each element of ${\mathcal{C}}'(P)$
gives rise to one of the configurations in ${\mathcal{C}}'$.
Our aim is to identify them now. The process is 
similar to \cite{Gconv:90}:
\begin{itemize}
\item[(i)] We show that a configuration 
corresponds to an $\rpt$-structure. 
\item[(ii)] We embed the configuration space 
to the deformation space as an open subset.
\item[(iii)] We show that the deformation space 
of convex $\rpt$-structures is a subset of 
the image, i.e., our construction above.
\item[(iv)] Show that the deformation space is 
an open and closed subset of the image 
using the Koszul openness and the closedness 
in their representation space. 
\end{itemize}

Let $\hat H$ be the space of triples $(A, B, C)$ 
that occur in ${\mathcal{C}}'$.  
Give $\hat H$ the subspace topology 
of the triple product $\PGL(3, \bR)^3$ of $\PGL(3, \bR)$.  
Then $\hat H$ is open: For $A'$, $B'$, $C'$ sufficiently
close to $A, B, C$ in $\hat H$, the set of attracting fixed 
points of $A'$, $B'$ and the isolated fixed point of $C'$ 
form vertices of a convex triangle. Any other collection of the triple 
points also form vertices of convex triangles. Thus, the 
seven tuple is well-defined.

We claim that an element $(A, B, C)$ of $\hat H$
determines an element of $\mathcal{C}'$, and the
projection $p_1:{\mathcal{C}'} \ra \PGL(3, \bR)^3$ defined by
\[(\tri_0, \tri_a, \tri_b, \tri_c, A, B, C) \mapsto (A, B, C)\] 
is a homeomorphism onto $\hat H$:
The eigenvectors of $A, B, C$ and their images under $A, B, C$
determine the triangles $\tri_0, \tri_a, \tri_b, \tri_c$
and hence we can define a section of this fibration from their image.
The continuity of the section follows
from that of the map from the subspace of hyperbolic matrices 
in $\PGL(3, \bR)$ to the space of fixed points of their 
largest eigenvalue is a continuous one and that of the map from 
the subspace of rotations in $\PGL(3, \bR)$ 
to their unique isolated fixed points. 

We can identify the space $\PGL^3(\bR)$ satisfying $ABC = \idnt$ 
as the space of representations $U(P)_p$.
Therefore, $p_1$  embeds ${\mathcal{C}}'$ to the subspace 
$\hat H$ of $U(P)_p$. 

The group $\PGL(3, \bR)$ acts properly on the subset of $(\rpt)^4$ consisting
of $4$-tuples of points, no three of which are collinear. Therefore 
the $\PGL(3, \bR)$-action on ${\mathcal{C}}'$, defined by
\begin{equation*}
\big(\tri_0, \tri_a, \tri_b, \tri_c, A, B, C\big)  
\stackrel{\vth}\mapsto
\big(\vth(\tri_0), \vth(\tri_a), \vth(\tri_b), \vth(\tri_c), \vth A\vth^{-1}, 
\vth B\vth^{-1}, \vth C \vth^{-1}\big),
\end{equation*}
for $\vth \in \PGL(3, \bR)$, is proper and free.
Let $\mathcal{C}$ denote the quotient space.  

Now, we demonstrate that
$\mathcal{C}(P)$ is diffeomorphic to $\mathcal{C}$ as in Proposition
4.4 of \cite{Gconv:90}: Given an element of ${\mathcal{C}}'$, consider
the $\rpt$-orbifold with the same underlying space as $P$,
but with the cone-point, boundary points, and mirror points removed. 
The open $2$-orbifold $P'$ 
extends to a homeomorph $P''$ of $P$, whose development pair is already
determined by the $7$-tuple:
Let $E$ be an end of $P'$ corresponding to a boundary component of
$P''$ and say $A$. The universal cover $\tilde P'$ of $P'$ develops 
to a convex domain tessellated by the images of the triangles 
under the holonomy group action. Then a component of 
the inverse image of an open neighborhood of $E$ in $\tilde P'$ is 
filled with triangles developing into triangles converging 
to a segment connecting the attracting and repelling fixed points of 
$A$. Adding the interior of the segment corresponds to 
the completion of $E$. If $E$ is an end of $P'$ corresponding 
to a cone-point. Then triangles meeting $E$ develops periodically
around the fixed point of a conjugate of $C$. Adding the fixed 
point corresponds to the completion of $E$.
Therefore the homeomorph $P''$ carries 
a real projective structure since the completion 
is a geometric operation of attaching a principal
geodesic boundary component and a cone-point. 
Therefore, such a $7$-tuple determines an
element of ${\rpts(P)}$, 
that is, an isotopy class of $\rpt$-structures on $P$.

%%5/1 9:00 I need to revise further below.
The embedding $\iota:{\mathcal{C}}'\longrightarrow{\rpts}(P)$ 
defined in this way is an imbedding onto 
an open subset of ${\rpts}(P)$: 
The composition $\mathcal{H}'\circ\iota$ 
equals projection $p_1$ from $\mathcal{C}'$ to 
the open subspace $\hat H$ in 
the space of representations $U(P)_p$.  
Because $\mathcal{H}'$ is a local homeomorphism 
by Theorem \ref{thm:Haus2} and $p_1$ is
continuous, $\iota$ is continuous.
Because $\mathcal{H}'$ is a local homeomorphism and $p_1$ is open,
$\iota$ is open and its image is an open subset of ${\rpts}(P)$.  
Because $p_1$ is injective, $\iota$ is injective.
Therefore, $\iota$ is an embedding onto 
a subset of ${\rpts}(P)$.

The set ${\mathcal{C}}'(P)$ corresponding to 
convex structures on $P$ is an open subset of $\rpts(P)$
by Proposition \ref{prop:openness2}.
Since a convex $\rpt$-structure on $P$ determines 
a $7$-tuple by our construction, 
${\mathcal{C}}'(P)$ is open in the image of $\iota$. 

By Proposition \ref{prop:closedness2}, the image of ${\mathcal{C}}'(P)$
under $\mathcal{H}'$ is a closed subset of $U(P)_p$.  Since
$\mathcal{H}'$ is a local homeomorphism, ${\mathcal{C}}'(P)$ is 
locally closed in ${\rpts}(P)$.  Therefore ${{\mathcal{C}}'(P)}$
is open and locally closed in the image of $\iota$.

The Teichm\"uller space $\mathcal{T}(P)$ of hyperbolic $\rpt$-structures
is a subset of $\mathcal{C}(P)$. 
$\mathcal{C}(P)$ is not empty since it contains $\mathcal{T}(P)\neq\emptyset$.
Therefore ${\mathcal{C}}'(P)\neq\emptyset$.  Since
the image of $\iota$ is connected and open and relatively closed, 
${\mathcal{C}}'(P)$ is the entire image of $\iota$.
This completes our identification of ${\mathcal{C}}(P)$
with ${\mathcal{C}}'$.

Now, we act by $\PGL(3, \bR)$ on ${\mathcal{C}}'(P)$ to obtain a
quotient space $\mathcal{C}(P)$ and act by the same group on
${\mathcal{C}}'$ to obtain $\mathcal{C}$. Since the above
correspondence is $\PGL(3, \bR)$-equivariant, we see that
$\mathcal{C}$ and $\mathcal{C}(P)$ are diffeomorphic
by the map induced by $\iota$.

The following proposition proves Proposition \ref{prop:elemdim} when
$P$ is an elementary $2$-orbifold of type (P2):
\begin{prop}\label{prop:elem(P2)} 
$\mathcal{C}$ is an open cell of dimension $6$ or $6 - 2 = 4$ 
depending on whether or not the order of the cone-point is $2$.  
The map
\begin{eqnarray}
\mathcal{C} & \ra & {\mathcal{R}}_A \times {\mathcal{R}}_B \\
(\tri_0, \tri_a, \tri_b, \tri_c, A, B, C) & \mapsto & 
((\lambda, \tau)_A, (\lambda, \tau)_B )
\end{eqnarray} 
is a principal fibration with fiber an open two-cell over the four-cell
${\mathcal{R}}^2$ for $n \geq 3$. 
If $n=2$, the map is a diffeomorphism.
\end{prop}

%11/13 Diagonal matrices, outline for P2 needed.
We begin the proof of Proposition \ref{prop:elem(P2)}:
We may put $\tri_0$ to a standard triangle with vertices $[1,0,0]$,
$[0,1,0]$, and $[0,0,1]$ in the homogeneous coordinates of $\rpt$.
The remaining vertices of $\tri_a, \tri_b,$ and $\tri_c$ are 
$[-1, b_1, c_1]$, $[a_2, -1, c_2]$, and $[a_3, b_3, -1]$ respectively.  
If $n > 2$, then $b_1, c_1, a_2, c_2, a_3, b_3$ are positive. 
If $n=2$, $C$ sends lines through its isolated fixed point to the same
lines with orientation reversed.  Thus,
\begin{equation}\label{eqn:n=2} 
a_2 = 0, b_1 = 0,
\end{equation}
and the rest are positive.

Denoting by $\mathcal{C}''$ the configurations of four triangles with 
the above vertices. We see that $\mathcal{C}(P)$ is a quotient of 
$\mathcal{C}''$ by a group of collineations conjugate to 
a group of diagonal matrices with positive eigenvalues. 

We will now analysis $\mathcal{C}''$ below using the above notations:
\subsection{(P2) The cone-point order $\ne 2$}
Now  assume $n \geq 3$.  If $n \geq 4$, then the angles of
the union hexagon are less than $\pi$.  Then 
\[b_1, c_1, a_2, c_2, a_3, b_3 >0.\] 
At each vertex of $\tri_0$, the cross-ratios of the
four lines containing the edges of the incident triangles, 
determine invariants
\[\rho_1 = b_3 c_2, \rho_2 = a_3 c_1, \rho_3 = a_2 b_1, \]
satisfying 
\begin{equation*}
\rho_1, \rho_2, \rho_3 > 1 
\end{equation*}
since $\tri_0\cup\tri_a \cup \tri_b\cup \tri_c$ is convex.
(See Figure \ref{fig:triangles}.)
If $n = 3$, then 
\begin{equation*}
 b_1, c_1, a_2, c_2, a_3, b_3 >0, \rho_1, \rho_2 > 1 
\end{equation*}
but $\rho_3$ may assume any positive real value. 
(When $\rho_3 = 1$, we have an angle $\pi$ at $[0,0,1]$.)

The group of diagonal matrices
\[\begin{bmatrix} 
\lambda & 0 & 0 \\
0 & \mu & 0 \\
0 & 0 & \nu \end{bmatrix}, 
\lambda\mu\nu = 1, \lambda, \mu, \nu > 0,
\]
acts on $\tri_0, \tri_a, \tri_b, \tri_c$, taking
\[ \begin{pmatrix} b_1 \\ c_1 \\ a_2 \\ c_2 \\ a_3 \\ b_3 \end{pmatrix} 
\mapsto \begin{pmatrix} (\mu/\lambda) b_1 \\ (\nu/\lambda) c_1 \\
(\lambda/\mu) a_2 \\ (\nu/\mu) c_2 \\ (\lambda/\nu) a_3 \\ (\mu/\nu)
b_3 \end{pmatrix}. \] 
Then $\sigma_1 = a_2b_3 c_1$ and
$\sigma_2 = a_3b_2c_2$ are invariants under the diagonal group action
and hence are invariants of the $7$-tuple under the action of
collineations. $\sigma_1, \sigma_2 > 0$ and
$\sigma_1\sigma_2 = \rho_1\rho_2 \rho_3$.

By applying such a diagonal matrix, we may assume 
$a_3 = 2$ and $b_3 = 2$, obtaining a slice for the $\PGL(3,\bR)$-action
on $7$-tuples:
\begin{equation}\label{eqn:coordinates}
a_2 = t, a_3 = 2, b_1 = \rho_3/t, b_3 =2,c_1 = \rho_2/2, c_2 = \rho_1/2 
\end{equation}
where $t = \sigma_1/\rho_2 > 0$ is an arbitrary positive number.
Elements of $\PGL(3, \bR)$ may be uniquely represented by 
elements of $\SL(3, \bR)$. 
The most general collineation sending 
$\tri_b$ to $\tri_c$ is given by the matrix
\begin{eqnarray}
A  &= & \alpha_1 \cdot \begin{bmatrix} 1 \\ 0 \\ 0 \end{bmatrix} 
\cdot [ 1, a_2, 0] + 
\beta_1 \cdot \begin{bmatrix} 0 \\ 1 \\ 0 \end{bmatrix} 
\cdot [0, -1, 0]  \nonumber \\
 & +& \gamma_1 \cdot \begin{bmatrix} a_3 \\ b_3 \\ -1 \end{bmatrix}
\cdot [0, c_2, 1] \nonumber \\
&=& \begin{bmatrix}
\alpha_1 & \alpha_1 a_2 + \gamma_1 a_3 c_2 & \gamma_1 a_3 \\
0 & -\beta_1 + \gamma_1 b_3 c_2 & \gamma_1 b_3 \\
0 & -\gamma_1 c_2 & -\gamma_1 \end{bmatrix} 
\end{eqnarray}
for $\alpha_1, \beta_1, \gamma_1 > 0$.
We obtain for $B$ and $C$ the following 
\begin{eqnarray}
B  &= & \alpha_2 \cdot \begin{bmatrix} -1 \\ b_1 \\ c_2 \end{bmatrix} 
\cdot [ 1, 0, a_3] + \beta_2 \cdot \begin{bmatrix} 0\\ 1\\ 0 \end{bmatrix} 
\cdot [0, 1, b_3] \nonumber \\
 & +& \gamma_2 \cdot \begin{bmatrix} 0 \\ 0\\ 1 \end{bmatrix}
\cdot [0, 0, -1]  \nonumber \\
&=& \begin{bmatrix}
-\alpha_2 & 0 & -\alpha_2 a_3 \\
\alpha_2 b_1 & \beta_2 & \beta_2 b_3 \alpha_2 a_3 b_1\\
\alpha_2 c_1 & 0 & -\gamma_2 + \alpha_2 a_3 c_1 \end{bmatrix} 
\end{eqnarray}
and 
\begin{eqnarray}
C  &= & \alpha_3 \cdot \begin{bmatrix} 1 \\ 0 \\ 0 \end{bmatrix} 
\cdot [-1, 0, 0] + \beta_3 \cdot \begin{bmatrix} a_2 \\ -1\\ c_2 \end{bmatrix} 
\cdot [b_1, 1, 0] \nonumber \\
 & +& \gamma_3 \cdot \begin{bmatrix} 0 \\ 0\\ 1 \end{bmatrix}
\cdot [c_1, 0, 1] \nonumber \\
&=& \begin{bmatrix}
-\alpha_3 + \beta_3 a_2 b_1 & \beta_3 a_2 & 0 \\
-\beta_3 b_1 &  -\beta_3 & 0\\
\gamma_3 c_1 + \beta_3 b_1 c_2 & \beta_3 c_2 & \gamma_3 \end{bmatrix} 
\end{eqnarray}
where 
\[\alpha_2, \beta_2, \gamma_2, \alpha_3, \beta_3, \gamma_3 > 0.\] 

$ABC = I$ and $\det(A) = \det(B) = \det(C) = 1$ imply
\begin{equation}\label{eqn:CBA}
\alpha_1 \alpha_2 \alpha_3 = \beta_1 \beta_2 \beta_3 = 
\gamma_1 \gamma_2 \gamma_3 = 1
\end{equation}
and 
\begin{equation}\label{eqn:det}
\alpha_1\beta_1\gamma_1 = \alpha_2\beta_2\gamma_2= 
\alpha_3 \beta_3\gamma_3 = 1.
\end{equation}
The invariants of $A$ and $B$ are given by:
\begin{equation}\label{eqn:lA}
\lambda(A) = \lambda_1 = \alpha_1, 
\end{equation}
\begin{equation}\label{eqn:tA}
\tau(A) = \tau_1 = -\beta_1 + \gamma_1(\rho_1 -1),
\end{equation}
\begin{equation}\label{eqn:lB}
\lambda(B) = \lambda_2 = \beta_2, 
\end{equation}
\begin{equation}\label{eqn:tB}
\tau(B) = \tau_2 = -\gamma_2 + \alpha_2(\rho_2 - 1).
\end{equation} 
For $C$ to be of order $n$ (for any $n \geq 2$),
\begin{equation}\label{eqn:lC}
\gamma_3 = 1, \mbox{ and } 
\end{equation}
\begin{equation}\label{eqn:tC}
-\alpha_3 + \beta_3 (\rho_3 -1) = 2\cos(2\pi/n).
\end{equation}
We now restrict our attention 
purely to the system of
equations from \eqref{eqn:CBA} to \eqref{eqn:tC}.
We wish to solve for $\alpha, \beta, \gamma, \rho$s given 
the invariant 
$(\lambda_1, \tau_1, \lambda_2, \tau_2) \in {\mathcal{R}}^2$
and the order $n \geq 3$. From equation \eqref{eqn:coordinates}, 
we obtain the coordinates of vertices.
It is clear that the set of such solutions correspond 
in one to one manner to ${\mathcal{R}}$. 
 
This system of equations appear in \cite{Gconv:90} where $\lambda_3$ is
replaced by $\gamma_3$ and $\tau_3$ by $2\cos(2\pi/n)$ where
$(\lambda_3, \tau_3)$ is the invariant of $C$ if $C$ were hyperbolic.

The way to solve the system is to realize that with $\alpha_1$ and
$\beta_2$ and $\gamma_3$ fixed by $\lambda_1, \lambda_2,$ and $1$ from
equations \eqref{eqn:lA}, \eqref{eqn:lB}, and \eqref{eqn:lC}, equations
\eqref{eqn:CBA} and \eqref{eqn:det} can be made into a system of linear
equations of rank five of six variables
\begin{equation*}
\log \alpha_2, \log \alpha_3, \log \beta_1, \log \beta_3,
\log \gamma_1, \log \gamma_2.
\end{equation*}
There is a space of solutions diffeomorphic to $\bR^+$, say
parameterized by a variable $s$. 
For each of the solutions of these equations, we plug into 
equations \eqref{eqn:tA}, \eqref{eqn:tB}, and \eqref{eqn:tC} to solve for 
$\rho_1, \rho_2, \rho_3$. Then we plug these to 
equation \eqref{eqn:coordinates} adding $\bR^+$-parameter 
from the variable $t$, obtaining a solution space diffeomorphic 
to $\bR^2$. Thus, we obtain for every fixed
\[(\lambda_1, \tau_1, \lambda_2, \tau_2) \in {\mathcal{R}}^2, 
(s, t) \in \bR^{+ 2},\] 
a unique solution. (Actually, 
equations 4-21 and 4-23 of \cite{Gconv:90} 
with setting $\lambda_3 = 1$ and $\tau_3 = 2\cos(2\pi/n)$
are the solutions for any choice of $s$ and $t$, $s, t >0$, 
there.) This completes the proof of Proposition \ref{prop:elem(P2)} 
when $n\geq 3$.

\subsection{(P2) the cone-point order $2$}
\begin{figure}[ht] %%
\centerline{\epsfxsize=3.8in \epsfbox{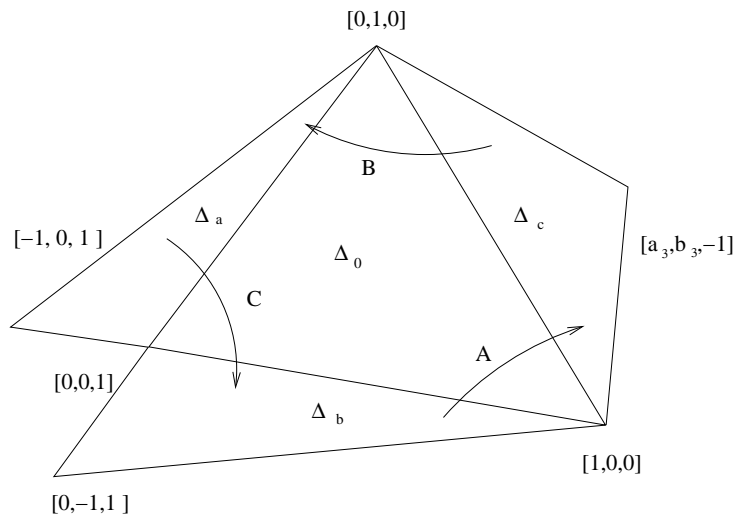}}
\caption{\label{fig:triangles2} 
The four triangles needed to understand the deformation
space of an annulus with one cone-point of order two.}
\end{figure}

Now consider $n = 2$.  Apply a unique collineation
(see Figure \ref{fig:triangles2}) so that 
$\clo(\tri_0)$ has  vertices 
\[ [1,0,0], [0,1,0], [0,0,1], \]
$\clo(\tri_a)$  is 
an adjacent triangle with vertices
\[ [0,1,0], [-1,0,1], [0,0,1],\] 
$\clo(\tri_b)$ is
an adjacent triangle with vertices
\[[0,-1,1], [1,0,0], [0,0,1], \] 
and $\clo(\tri_c)$ is an adjacent triangle with vertices 
\[[0,1,0], [1,0,0], [a_3, b_3, -1].\]
We can define cross-ratios $\rho_1, \rho_2$ as above.
We have 
\begin{equation}\label{eqn:coordinates2}
\rho_1 = b_3, \rho_2 = a_3, \rho_3 = 0, c_1 =1, c_2 = 1
\end{equation}
by equation \eqref{eqn:n=2}.
From the condition that $C^2 = \idnt$, 
and equation \eqref{eqn:n=2}, we see that 
$\alpha_3 = \beta_3 = \gamma_3 = 1$. 
The equations \eqref{eqn:CBA} and \eqref{eqn:det} and 
equations \eqref{eqn:lA}, \eqref{eqn:tA}, \eqref{eqn:lB}, 
\eqref{eqn:tB}, \eqref{eqn:lC}, and \eqref{eqn:tC}  
still apply here. Then with 
\[(\lambda_1, \tau_1, \lambda_2, \tau_2) \in {\mathcal{R}}^2\]
fixed, we obtain values of 
$\alpha_1, \beta_2$ by equations \eqref{eqn:lA}
and \eqref{eqn:lB} and hence
a unique solution for 
equations \eqref{eqn:CBA} and \eqref{eqn:det}.
Now, we plug the solution to equation \eqref{eqn:tA} 
and equation \eqref{eqn:tB} to obtain
$\rho_1$ and $\rho_2$. Equation \eqref{eqn:tC} is automatically
satisfied. This determines the seven-tuple.--(**)
\qed

\typeout{<>}

\subsection{A disk with two cone-points (P3)}
Now, we determine the topology of the deformation space of the
elementary $2$-orbifolds of type (P3).  Let $m$ and $n$ be the order
of the two cone-points.  Assume $m \geq 3$ without loss of generality.
By the same reasoning as (P2), 
the proof reduces again to determining the space $\mathcal{C}$ of
equivalence classes under the diagonal group action of $\mathcal{C}'$
of $7$-tuples $(\tri_0, \tri_a, \tri_b, \tri_c, A, B, C)$ where
$\tri_0$ is the standard triangle and $\tri_a, \tri_b, \tri_c$ are
adjacent triangles, $A$ is a collineation sending $\tri_b$
to $\tri_c$ and $B$ one sending $\tri_a$ to $\tri_c$ of order $m$ and
$C$ one sending $\tri_c$ to $\tri_a$ of order $n$. 
That is, we show
similarly to (P2) that $\mathcal{C}(P) = \mathcal{C}$ again by
openness and closedness.  

Assume for the moment $n\geq 3$ also.  We can introduce unknowns $b_1, c_1,
a_2, c_2, a_3, b_3$, $\rho_1, \rho_2, \rho_3$, $\sigma_1, \sigma_2$,
$\alpha_i, \beta_i, \gamma_i$, $i=1,2,3$, as above.  The only change
of the equations occur at the equations \eqref{eqn:lB} and \eqref{eqn:tB}
which change to
\begin{equation}\label{eqn:lB'}
\beta_2 = 1,
\end{equation}
\begin{equation}\label{eqn:tB'}
-\gamma_2 + \alpha_2(\rho_2 - 1) = 2\cos(2\pi/m). 
\end{equation}
The same procedure implies
${\mathcal{F}}: \mathcal{C}(P) \ra {\mathcal{R}}_A$ is 
an principal $\bR^2$-fibration where ${\mathcal{R}}_A$ 
is the space of invariants of $A$. Here $\mathcal{C}(P)$ has dimension four.

When $n=2$, an argument similar to (**)  $\mathcal{F}$ is a diffeomorphism
$\mathcal{C}(P) \ra {\mathcal{R}}_A$ 
so that $\mathcal{C}(P)$ is two-dimensional. (They confirm 
the conclusion of Proposition \ref{prop:elemdim}.)

\subsection{Spheres with three cone-points (P4)}
For elementary $2$-orbifolds of type (P4), first suppose that $r, m, n
\geq 3$, then the same considerations as above apply if we change
the equations \eqref{eqn:lA} and \eqref{eqn:tA} as well by
\begin{equation}\label{eqn:lA'}
\alpha_1 = 1, 
\end{equation}
\begin{equation}\label{eqn:tA'}
-\beta_1 + \gamma_1(\rho_1 -1)= 2\cos(2\pi/r).
\end{equation}
Thus $\mathcal{C}(P)$ is a $2$-cell.
% which equals $6\times 3 - 8\times 2$. There is no boundary
% consideration here. 

If $n=2$, then both $r\ge 3$ and $m\ge 3$ since $\chi(\Sigma)<0$.
So the second part of 
(**) applies, and  $\mathcal{C}(P)$ is a single point.

\subsection{Crowns with two prongs (A1)}
Again, we reduce an element of $\mathcal{C}(P)$ to 
a point of a certain configuration space.
(Here a method is slightly different.)
For an elementary $2$-orbifold $P$ of type (A1),
$\dev(\tilde P)$ for the universal cover $\tilde P$ of $P$ is a convex
disk invariant under $h(\vth)$ for the deck transformation $\vth$
corresponding to its principal closed geodesic boundary component
with the boundary orientation. 
A hyperbolic element $h(\vth)$ has three fixed points in $\rpt$ and
three lines through two of them, and four open triangles bounded by
subsegments of them are $h(\vth)$-invariant. The developing image
$\dev(\tilde P)$ is convex, and lies inside the closure of one of these open
triangles. Choose projective coordinates so that the open triangle 
is the standard coordinate triangle, that is, the repelling fixed point
of $h(\vth)$ is $[0,0,1]$, the attracting fixed point $[0,1,0]$ and
the saddle type fixed point $[1,0,0]$.  (Figures are similar to
elementary $2$-orbifolds of type (P2) here.)

\begin{figure}[ht] %%
\centerline{\epsfxsize=4.8in \epsfbox{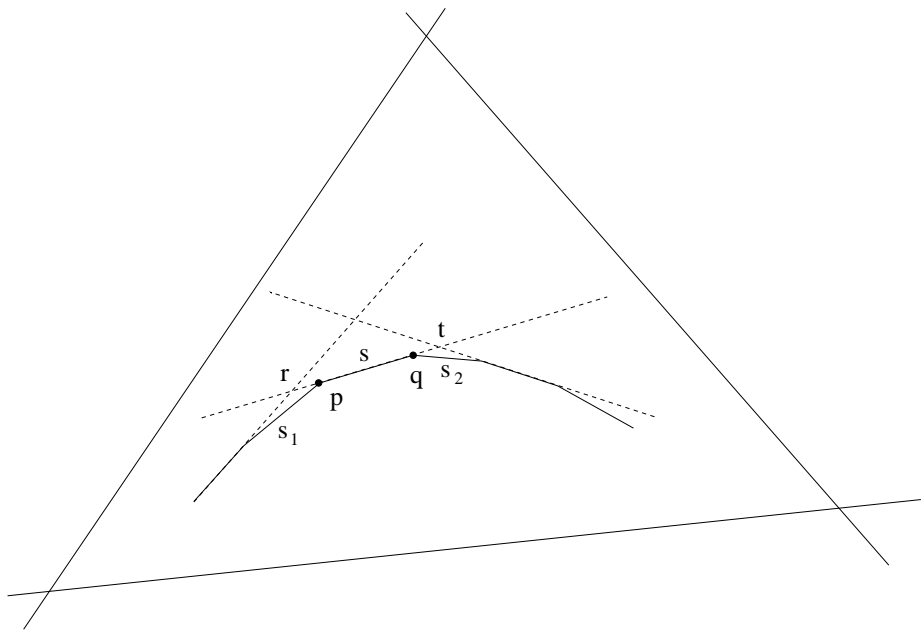}}
\caption{\label{fig:annulus} 
The annulus and reflections.}
\end{figure}
\typeout{<>}

%%5/03/11:00
Given a convex $\rpt$-structure on $P$, its universal covering space 
$\tilde P$ contains a disk $K'$ with two boundary components 
$l_1$ and $l_2$ where $l_1$ covers a principal
closed geodesic boundary of $P$ under the universal covering map and
$l_2$ covers the union of a segment in the singular locus
and a boundary $1$-orbifold.  
Let $p'$ and $q'$ be two endpoints of the boundary
$1$-orbifold in $P$, and let $p$ and
$q$ be two endpoints of a segment in $l_2$ corresponding to $p'$ and
$q'$ respectively; let $s$ be the segment connecting $p$ and $q$.
Let $s_1$ and $s_2$ be the other two segments in $l_2$ 
starting from $p$ and $q$ respectively. Then $\vth(s_1) = s_2$.

%5/2/11:45
The developing map $\dev$ embeds $K'$ to a convex domain 
in the standard triangle so that $l_1$ maps to an open segment 
connecting $[0,0,1]$ to $[0,1,0]$, and $l_2$ to an arc consisting 
of segments bent in one direction connecting $[0,0,1]$ to
$[0,1,0]$. Identify $\tilde P$ and associated objects 
with their developing images or holonomy images 
since $\dev$ is an embedding.
There is a unique $s_1$-invariant reflection $r_1$ 
whose fixed line contains $s_1$. Reflection in  $s_2$ equals
$r_2 = \vth r_1 \vth^{-1}$. For each segment 
of form $\vth^i(s_1)$, $i \in \bZ$, has an associated reflection 
$\vth^i r_1 \vth^{-i}$. Since $s$ corresponds to a $1$-orbifold, 
$r_1$ and $r_2$ act on the line $l(s)$ containing $s$. Therefore, 
the isolated fixed points $r$ and $t$ of $r_1$ and $r_2$ 
respectively lie on $l(s)$.
Since $r_1$ also acts on the line containing $\vth^{-1}(s)$,
$r$ must lie on $\vth^{-1}(l(s))$ the fixed line 
of $\vth^{-1} r_1 \vth$. Thus, $r$ is the unique intersection point of 
$l(s)$ and $\vth^{-1}(l(s))$. Similarly $t$ is the unique fixed point of 
$l(s)$ and $\vth(l(s))$. By post-composing $\dev$ 
by a diagonal matrix, we may assume that $r$ correspond to $[1,1,1]$ and 
$t$ to $\vth([1,1,1])$.
Thus $p$ and $q$ lie on $\ovl{[1,1,1]\vth([1,1,1])}$. Recall that $D$ is  
bounded by $l_1$ and $l_2$ the union of $\ovl{\vth^i(p)\vth^i(q)}$ 
and $\ovl{\vth^{i-1}(q)\vth^i(p)}$ for $i \in  \bZ$.

We showed that each element of $\mathcal{C}(P)$ corresponds 
to \[(\vth, p, q), p, q \in \ovl{[1,1,1]\vth([1,1,1])} \]
where $\vth$ is a positive projective transformation 
represented by a diagonal matrix so that $[0,0,1]$ 
correspond to the eigenvector with least eigenvalues,
and $[1,0,0]$ to one with the largest eigenvalue.

%5/2 9:00
%%%May need more work! I need to verify homeomorphism.
We can identify the space of such configurations $\mathcal{C}$ 
with $\mathcal{C}(P)$: Any element of $\mathcal{C}(P)$ 
correspond to a unique such element by the above standardization process. 

We can introduce 
the topology to the space $\mathcal{C}$ of above configurations 
as a subspace of $D \subset \PGL(3, \bR)$ and $A \times A$ where 
where $D$ is the space of diagonal matrices with 
strictly decreasing positive eigenvalues and
$A$ is the interior of the standard triangle.
Clearly, $\mathcal{C}$ is diffeomorphic to $\bR^2 \times \bR^2$.
Then the map from $\mathcal{C}(P)$ to $\mathcal{C}$ 
defined by the above work is clearly 
continuous since $\vth, p, q$ is determined by the holonomy
homomorphisms which depends continuously on the projective
structures. 

We now try to obtain the inverse $\iota$ of this map:
From an element of $\mathcal{C}$, we form the domain $K'$ 
in the standard triangle bounded by 
\[\vth^n(\ovl{[1,1,1]\vth([1,1,1])}), n \in \bZ, 
\ovl{[1,0,0][0,0,1]}.\]
More precisely, $K'$ is defined to be the 
union of the interior and the segments
\[\vth^n(\ovl{[1,1,1]\vth([1,1,1])}), n \in \bZ\]
and the interior of the segment 
$\ovl{[1,0,0][0,0,1]}$.

Let $r$ be the reflection with isolated fixed point 
$[1,1,1]$ and the segment of fixed points $\ovl{p\vth^{-1}(q)}$.
$K'$ and $r_1(K')$ meet in a side and adjacent 
sides extend each other as $r$ preserves their ambient 
lines. Thus, $K'\cup r(K')$ is convex
and so is $K'\cup r(K') \cup \vth r \vth^{-1}(K'\cup r(K'))$. 
Continuing in this manner, we see that 
a connected union of finitely many images 
of $K'$ under $\langle \vth, r \rangle$-action is convex. 
Therefore, $K^\infty$ equal to the union of all the images of
of $K'$ under $\langle \vth, r \rangle$-action is 
convex. Then \[K^\infty/\langle \vth, r \rangle\]
is an elementary orbifold of type (A1), 
which is convex.  

The map $\iota$ from $\mathcal{C}$ to $\mathcal{C}(P)$ defined 
in the above manner is also continuous: 
By definition of the topology, the composition 
\[\mathcal{H}\circ \iota: \mathcal{C} \ra U(P)_p/\PGL(3, \bR)\]
is continuous. Since $\mathcal{H}$ is a local homeomorphism, 
$\iota$ is continuous. 
(Here the proof is a little bit different from cases 
(P1),(P2),(P3), and (P4) since we can directly construct 
a convex projective structure from a point of 
the configuration space. However, there should be 
direct methods in (P1)-(P4) as well which we were 
unable to find.) 

We obtained that $\mathcal{C}(P)$ is diffeomorphic to 
a $4$-cell $\mathcal{C}$. 

The cross-ratio 
$[[1,1,1], \vth([1,1,1]); p, q] \in (0,1)$ is 
the invariant of the boundary $1$-orbifold. 
Fixing this cross-ratio, the space of choices 
of $p$ and $q$ on $\ovl{[1,1,1]\vth([1,1,1])}$ is 
diffeomorphic to $\bR$.
This becomes a fiber of the fiber map 
\[\mathcal{F}: \mathcal{C}(P) \ra {\mathcal{R}} \times (0,1)\]
given by sending the triple $(\vth, p, q)$ 
to the invariants of $\vth$ and the cross ratio. 

\subsection{One-pronged crowns (A2)}
We now go to elementary $2$-orbifolds of type (A2). 
Suppose that the order $p$ of the corner-reflector 
is greater than or equal to $3$. 
$\tilde P$ contains a disk 
$K'$ with two boundary arcs $l_1$ and $l_2$, 
where $l_1$ covers the principal closed geodesic boundary, and 
$l_2$ covers the other boundary component. Identify 
$K'$ with a convex domain in a standard triangle with 
$\vth$ acting on with properties as in (A1).
Here $\vth$ is the holonomy of the deck transformation
corresponding to the principal geodesic boundary component of $P$,
and has a diagonal matrix.

\begin{figure}[ht] %%
\centerline{\epsfxsize=4.8in \epsfbox{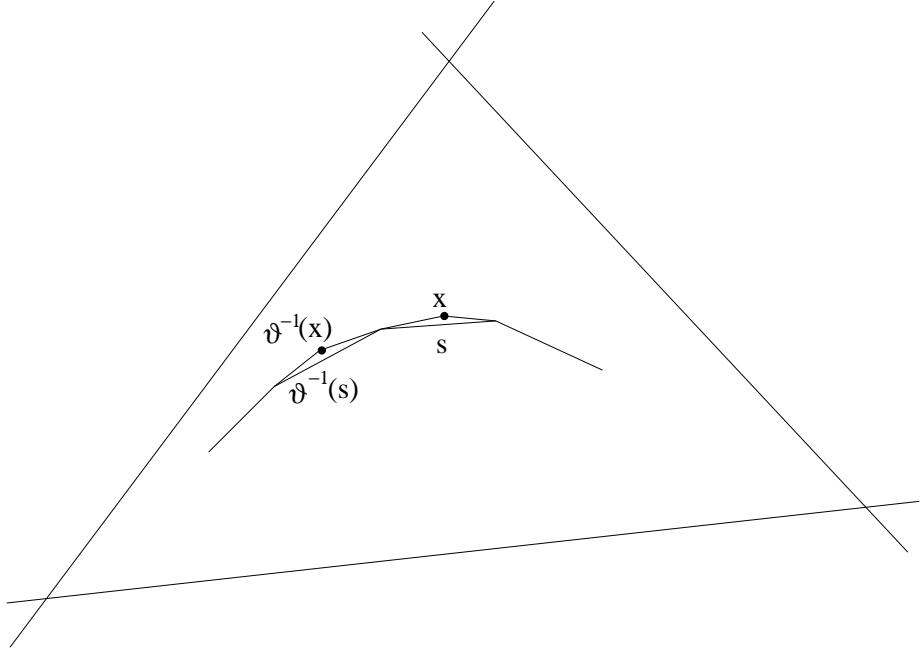}}
\caption{\label{fig:annulustwo} 
The annulus and reflections.}
\end{figure}
\typeout{<>}

We also identify associated objects.
$l_2$ is a union of segments of form $\vth^i(s)$, $i \in \bZ$, for
a segment $s$ with endpoints $p, \vth(p)$.
We post-compose $\dev$ with a diagonal matrix so that $p= [1,1,1]$ where 
$\vth$ does not change. The holonomy group contains a reflection $r$ 
with a line of fixed points containing $s$. $\vth^n r \vth^{-n}$ is 
a reflection with a line of fixed points containing the segment $\vth^n(s)$.
Let $x$ be the isolated fixed point of $r$. Then $\vth^{-1}(x)$ is
the isolated fixed point of $\vth^{-1} r \vth$. Let $l(x)$ and 
$l(\vth^{-1}(x))$ be the line from $[1,1,1]$ to points $x$ and 
$\vth^{-1}(x)$ respectively.  The corner-reflector 
has order $p \geq 3$, which is equivalent to the condition 
that the cross-ratio of lines satisfy
\begin{equation}\label{eqn:orderp}
[l(x), l(\vth^{-1}(x)); s, \vth^{-1}(s)]
= \frac{2}{\cos \frac{2\pi}{p} +1} 
\end{equation}
as shown in \cite{G:77}.
Thus, $x$ satisfies this condition. 

We compute the solution space of \eqref{eqn:orderp}:
For each line $l$ through $[1,1,1]$, we find a line $L(l)$ 
through $[1,1,1]$ so that 
\[[l, L(l); s, \vth^{-1}(s)] = 
\frac{2}{\cos \frac{2\pi}{p} +1}.\]
This is a projectivity $L$ from the pencil 
$F([1,1,1])$ of lines through $[1,1,1]$
to the same pencil. Using the map $\hat \vth$ from $F([1,1,1])$
to the pencil $F(\vth([1,1,1])$ of lines through $\vth([1,1,1])$
induced by $\vth$, we obtain a projectivity
\[\hat \vth\circ L: F([1,1,1]) \ra F(\vth([1,1,1])).\] 
By a classical result of Steiner 
(Chapters 6 and 7 of Coxeter \cite{Cox:60}),
the locus of 
$l \cap \hat \vth \circ L(l)$ for $l \in F([1,1,1])$ is a conic. 
Furthermore it passes through $[1,1,1]$ 
and $\vth([1,1,1])$. Since $x$ lies outside $D$, 
$x$ must be the unique connected arc in the conic outside $D$. 

For an element of $\mathcal{C}(P)$, we obtain a pair 
\[(\vth, x)\] 
where $x$ lies on the arc through $[1,1,1]$ and $\vth([1,1,1])$
given by above equations. 

As in case (A1), we can show that $\mathcal{C}(P)$ identifies with
the space of such configurations: We follow (P2) in this case
with a sketchy argument. We define a map 
\[ \iota: \mathcal{C} \ra \rpt(P)\]
first. We find a domain $K'$ bounded by 
\[\ovl{[1,0,0][0,0,1]},\vth^n(\ovl{[1,1,1]\vth([1,1,1])}),n\in \bZ;\]
more precisely, $K'$ is the union of the interior 
with the interior of $\ovl{[1,0,0][0,0,1]}$ 
and the segments above. 
There is a reflection $r$ with the isolated fixed point 
at $x$ and with the segment of fixed points 
$\ovl{[1,1,1]\vth([1,1,1])}$.
$K'/\langle \vth \rangle$ is an annulus
with principal boundary component $a$ and a once-broken geodesic 
circle boundary component $b$. We can silver the broken 
geodesic using $r$ and we obtain a point of $\rpt(P)$. 

Since $\mathcal{H}\circ \iota: \mathcal{C} \ra U(P)_p/\PGL(3, \bR)$ is 
a function sending a configuration $(\vth, x)$ to 
$(\vth, r)$ where $r$ is described above, it is 
continuous. Since $\mathcal{H}$ is a local homeomorphism 
$\iota$ is continuous. By the above description, 
the image of the composition is open. Thus, the image 
of $\iota$ in $\rpt(P)$ is open. The set $\mathcal{C}(P)$ 
is a subset of the image of $\iota$ since a convex $\rpt$-structure 
gives rise to a configuration by our construction. 
The closedness of $\mathcal{C}(P)$ in the image of 
$\iota$ follows since the image $\mathcal{C}(P)$ 
under $\mathcal{H}$ maps to a closed set in $U(P)_p/\PGL(3, \bR)$.
We can show as in (P2), that $\mathcal{C}(P)$ is an open 
and closed subset of the image of $\iota$ and hence 
equal to the image. 

Thus $\mathcal{C}(P)$ is diffeomorphic to 
a $3$-cell $\mathcal{C}$. The mapping
\[{\mathcal{F}}: \mathcal{C}(P) \ra {\mathcal{R}}_\vth \]
is a fibration with fibers the above arcs 
where ${\mathcal{R}}_\vth$ is the space of invariants of $\vth$. 

Suppose that the order of the corner-reflector is $2$. 
Then using the above notations, we see that the point $x$ lies on the 
line containing the segment $\vth^{-1}(s)$; and the point $\vth^{-1}(x)$ on
the line containing the segment $s$; $x$ is the unique intersection
point of lines $l(\vth^{-1}(s))$ and $\vth(l(s))$.
Thus, there is a unique choice of $x$
given $\vth$, and $\mathcal{C}(P)$ is a two-cell and 
\[{\mathcal{F}}: \mathcal{C}(P) \ra {\mathcal{R}}_\vth \] 
is a homeomorphism. 

%%5/03 2:35
\subsection{Disks with one boundary full $1$-orbifolds (A3)}
Again, we reduce an element of $\mathcal{C}(P)$ to a configuration:
First in case (A3), let $\vth$ be the deck-transformation 
associated with the cone-point of order $n$, $n\geq 3$. 
The universal cover $\tilde P$ of $P$ contains a disk where $\vth$ 
acts on so that its boundary map to the boundary of the underlying space 
$X_P$ of $P$. Let $s'_1$ be the boundary $1$-orbifold and $s'_2$
the segment in the singular locus of $P$.
Thus the disk is a convex polygon with $2n$ sides
\[s_1, s_2, \vth^1(s_1), \vth^1(s_2), \dots, \vth^{n-1}(s_1), \vth^{n-1}(s_2)\]
where $s_1$ map to $s'_1$ and $s_2$ to $s'_2$ in $P$.
The holonomy group contains a reflection $r$ with the line $l(s_2)$ 
of fixed points containing $s_2$. 
$\vth^i r \vth^{-1}$ has the line $\vth^i(l(s_2))$ of fixed points 
containing $\vth^i(s_2)$. Again, we find coordinates so that 
the isolated fixed point lie in $[0,0,1]$ and $\vth$ has a matrix form 
\begin{equation} 
\begin{bmatrix} 
\cos 2\pi/n & \sin 2\pi/n & 0 \\
-\sin 2\pi/n & \cos 2\pi/n & 0 \\
0 & 0 & 1
\end{bmatrix}.  
\end{equation} 
Since $s_1$ is a geodesic $1$-orbifold, $r$ and $\vth^{-1} r \vth$ act
on the line $l(s_1)$ containing $s_1$. Let $x$ be the isolated fixed point 
of $r$. Then $\vth^{-1}(x)$ is the isolated fixed point 
for $\vth^{-1} r \vth$. 
The fixed point $x$ lies on $l(s_1)$, as $r$ acts on $l(s_1)$, and 
on $l(\vth^{-1}(s_1))$, as $ r $ acts on $l(\vth^{-1}(s_1))$.
Thus, $x$ is the unique intersection point of $l(s_1)$ and 
$l(\vth^{-1}(s_1))$. The point $\vth(x)$ is the isolated fixed 
point of $\vth r \vth^{-1}$ with the line of fixed points $l(s_2)$. 
Since $\vth r \vth^{-1}$ acts on $l(s_1)$, $\vth(x)$ lies 
on $l(s_1)$. $\vth r \vth^{-1}$ acts on $l(\vth(s_1))$, 
$\vth(x)$ is the unique intersection point of $l(s_1)$ and 
$l(\vth(s_1))$. Post-compose $\dev$ with a matrix of form 
\begin{equation}
\begin{bmatrix} 
k \cos \theta & k\sin \theta & 0 \\
-k \sin \theta & k \cos \theta & 0 \\
0 & 0 & 1/k^2 
\end{bmatrix}, k > 0,
\end{equation} 
to send $x$ to $[1,1,1]$ and $\vth$ is unchanged.
Let $p$ and $q$ be endpoints of $s_1$ with $p$ separating
$x$ from $q$. The space of possible choices of $p$ and $q$ 
is diffeomorphic to $\bR^2$. The cross-ratio 
$[[1,1,1], \vth([1,1,1]), p, q]$ is the invariant of the boundary 
$1$-orbifold corresponding to $s_1$. 

Our configuration space $\mathcal{C}$ is the set of 
points $p, q \in \ovl{[1,1,1]\vth([1,1,1])}$. 
Given a point of $\mathcal{C}$, we can determine the reflection 
$r$, and as in case (A1), we show that a point 
of $\mathcal{C}$ gives rise to a convex orbifold of type (A3).
We identify $\mathcal{C}$ with $\mathcal{C}(P)$ by 
a homeomorphism. Thus $\mathcal{C}(P)$ is 
diffeomorphic to $\bR^2$ as $\vth$ is fixed now.

% unlike the case for (A1).  
There is an $\bR$-parameter family  of choices of
points to be labeled by 
$p$ and $q$ when the boundary invariant is fixed. 
This is a fiber of the fibration 
${\mathcal{F}}: \mathcal{C}(P) \ra (0,1)_{s'_1}$ where 
$(0,1)_{s'_1}$ is the space of invariants of the boundary 
$1$-orbifold $s'_1$.

\subsection{Disks with one corner-reflectors (A4)}
% We go over to the type-(A4)-orbifolds. 
Let $P$ be an orbifold of type (A4) with a corner-reflector of order $n$ 
and a cone-point 
of order $m$ where $1/n + 2/m < 1$. Then $ m \geq 3$. 
The proof is completely analogous to 
the two preceding cases.
If $n \geq 3$, then $\mathcal{C}(P)$ is diffeomorphic to $\bR$. 
If $n=2$, then $\mathcal{C}(P)$ is a single point.

\subsection{Pentagons (D2)}
% We go to pentagons (D2). ((D1) will be treated shortly.) 
Let $P$ be an orbifold of type (D2).
Suppose that $n \geq 3$. 
Let $e'_1, \dots, e'_5$ denote the edges and boundary orbifolds in 
the boundary of $X_P$ ordered in an appropriate orientation.
Let $e'_2, e'_3, e'_5$ be the edges and 
$e'_1$ and $e'_4$ be boundary $1$-orbifolds. 
Let $e'_2$ and $e'_3$ be the edges in the singular locus 
meeting in the corner-reflector of order $n$. 
$\tilde P$ contains a disk $D$ bounded by five geodesics 
$e_1, e_2, e_3, e_4,$ and $e_5$ so that $e_i$ maps to $e'_i$
for $i=1, \dots, 5$. Let $v_1, v_2, v_3, v_4, v_5$ denote 
the vertices of $D$ so that $e_i$ has vertices $v_i, v_{i+1}$ 
with cyclic indices. $\dev$ maps $D$ into a convex pentagon
in an affine patch of $\rpt$. Post-compose 
$\dev$ with a collineation so that 
$v_1, v_2, v_4, v_5$ map to 
\[[0,0,1], [1,0,1], [1,1,1], [0,1,1]\]
respectively. Then $v_3$ maps to a point $[s, t, 1]$ where $s > 1$ and 
$0 < t < 1$ by convexity of the pentagon. Identify $D$ and 
associated objects with their images in $\rpt$. 
Let $l_i$ denote the unique line containing $e_i$, $i=1, \dots, 5$. 
For $e_2, e_3, e_5$ there are associated reflections
$r_2, r_3$, and $r_5$ respectively. Their fixed lines 
are $l_2, l_3,$ and $l_5$ respectively. 
Let $x_2, x_3,$ and $x_5$ denote the respective 
isolated fixed points. $e_1$ and $e_4$ correspond to boundary
$1$-orbifolds. Since $r_2$ and $r_5$ act on $l_1$ 
and $r_4$ and $r_5$ on $l_2$, we have $x_2, x_5 \in l_1$ 
and $x_3, x_5 \in l_2$; that is, $x_5 = [1,0,0]$
and $x_2 = [t_2, 0, 1]$ and $x_3 = [t_3, 1, 1]$ 
with $t_2, t_3 > 1$. We compute the invariants of $e_1$ and $e_4$ 
\begin{eqnarray}\label{eqn:invt} 
\begin{bmatrix} [t_2, 0, 1], & [1,0,0]; & [1,0,1], & [0,0,1] 
\end{bmatrix} & = & \frac{t_2-1}{t_2} \nonumber \\
\begin{bmatrix} [t_3, 1, 1], & [1,0,0]; & [1,1,1], & [0,1,1] 
\end{bmatrix}  & = & \frac{t_3-1}{t_3}.
\end{eqnarray}
Since the corner-reflector has order $n$,
\begin{equation}\label{eqn:pentagon}
\begin{bmatrix} \ovl{v_3[t_2,0,1]}, \ovl{v_3[t_3,1,1]}; 
\ovl{v_3 [1,0,1]}, \ovl{v_3[1,1,1]} \end{bmatrix}
 = \frac{2}{\cos \frac{2\pi}{n} + 1}.
\end{equation}
We define our configuration space $\mathcal{C}$ to by 
$v_2, v_4$ satisfying equations \eqref{eqn:invt}
and $v_3$ satisfying equation \eqref{eqn:pentagon}.
Therefore, we found a map from $\mathcal{C}(P)$ to
the configuration space $\mathcal{C}$.

Let us explore the topology of the configuration
space $\mathcal{C}$:
Fix $t_2$ and $t_3$.
Steiner's theorem (Chapters 6 and 7 of 
Coxeter \cite{Cox:60}) implies the set of solutions $v_3$s
is a conic through $[t_2,0,1],[t_3,1,1], [1,1,1],$ and $[1,0,1]$. 
The subarc of the conic in the triangle with vertices 
$[t_2,0,1], [t_3,1,1], [1,0,0]$ outside the quadrilateral is 
the solution space for $v_3$. 

\begin{figure}[ht] %%
\centerline{\epsfxsize=3.8in \epsfbox{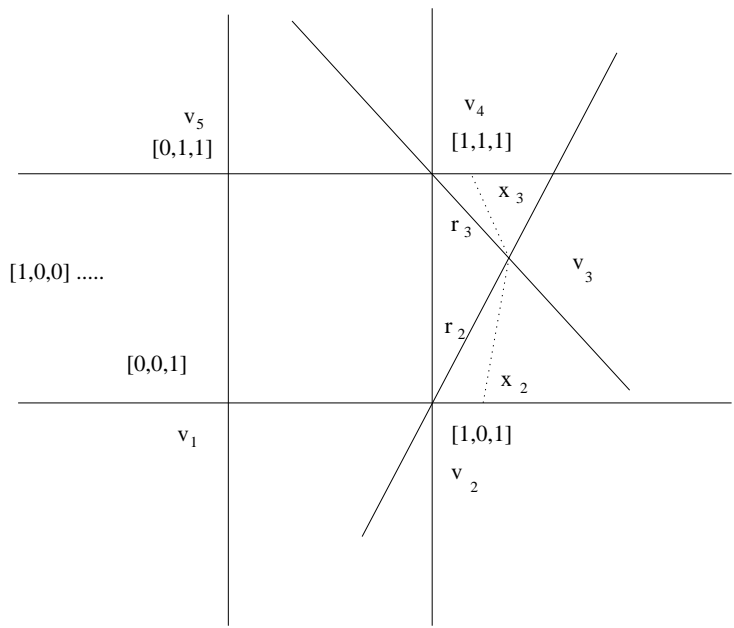}}
\caption{\label{fig:pentagon} 
The deformation spaces of pentagons.}
\end{figure}
\typeout{<>}

The above equation \eqref{eqn:pentagon} 
admits a solution since 
every cross-ratio between $0$ to $1$ in the equation 
is realized by a point on the triangle.
For each $t_2$ and $t_3$, there is 
an arc of solutions.
Therefore, $\mathcal{C}$ is diffeomorphic to $\bR^3$
and fibers over $(0,1)\times (0,1)$ with fibers 
equal to arcs given by the Steiner's theorem. 

The map from $\mathcal{C}$ to $U(P)_p$ given by
assigning to the configuration the reflections 
$r_1, r_2, r_3$ is obviously continuous. 
Also, given a point of $\mathcal{C}$, we can 
construct an $\rpt$-structure on $P$. 
Therefore, as in case (P2), we can identify
$\mathcal{C}(P)$ with $\mathcal{C}$, and 
the map ${\mathcal{F}}: \mathcal{C}(P) \ra (0,1)_{e_1} \times (0,1)_{e_4}$ 
is a principal  $\bR$-fibration where $(0,1)_{e_1}$ 
and $(0,1)_{e_4}$ are the space of invariants of $e_1$ and $e_4$ 
respectively. 

Next assume the cone-point order $n$ equals $2$.
Then $x_5=[0,0,1]$ and $x_2$ must be 
the intersection of the line through $v_3$ and $[1,0,1]$ 
with $l_1$. 
Furthermore
$x_3$ is the intersection of the line
through $v_3$ and $[1,1,1]$ with $l_2$. 
Let  $x_2 = [t_2,0,1]$ and $x_3 = [t_3,1,1]$.
Given $t_2, t_3$, let $v_3$ be 
the unique intersection point of the line 
through $x_2$ and $[1,1,1]$ with the line through
$x_3$ and $[1,0,1]$. Thus $\mathcal{C}(P)$ is a $2$-cell 
and 
\[{\mathcal{F}}:\mathcal{C}(P) \ra (0,1)_{e_1} \times (0,1)_{e_2}\] 
is a diffeomorphism.

%5/03/7:45
\subsection{Hexagons (D1)} Let $P$ be a hexagon.
The universal cover $\tilde P$ of $P$
contains a disk with boundary edges $e_1, e_2, \dots, e_6$ and
vertices $v_1, v_2, \dots, v_6$ so that $e_1, e_3,$ and $e_5$ are in
the interior of $\tilde P$ and correspond to singular edges of $P$ and
$e_2, e_4,$ and $e_6$ are geodesic segments mapping to the boundary
$1$-orbifolds of $P$.  Choose a developing map so that 
\[ \dev(v_1) = [0,0,1], \dev(v_2) = [1,0,1], \dev(v_3) = [1,1,1], 
\dev(v_6) = [0,1,1].\]  Then we have
\[\dev(v_4) = [s_4, t_4, 1], \dev(v_5) = [s_5, t_5, 1],
0< s_4, s_5 < 1, 1 < t_4, t_5.\]
As above identify objects
in $\tilde P$ with those in $\rpt$ by $\dev$. There are reflections
$r_1, r_3$, and $r_5$ whose lines of fixed points contain $e_1, e_3$,
and $e_5$ respectively. Let $x_1, x_3, x_5$ denote the respective
isolated fixed points. Let $l_i$ denote the line containing $e_i$ for
$i=1, \dots, 6$.  Since $l_2$ and $l_6$ are $r_1$-invariant, $x_1$ is
the intersection point of $l_2$ and $l_6$, that is, $x_1 = [0,1,0]$.
Since $l_2$ and $l_4$ are $r_2$-invariant, $x_3$ is of form $[1, t_3,
1]$ for $t_3 > 1$, and similarly $x_5 = [0, t_5, 1]$ for $t_5 > 1$.
The invariants for $e_2$ and $e_6$ are given by $(t_3-1)/t_3$ and
$(t_5-1)/t_5$ respectively.  Let $l_{35}$ be the line through $x_3$
and $x_5$. Then $l_{35}$ meet $l_3$ at $v_4$ and meet $l_5$ at
$v_5$. The invariant for $e_4$ is given by the cross-ratio $[x_3, x_5;
v_4, v_5]$.

The invariants for $e_2$ and $e_6$ determine  $t_3$ and $t_5$.
Given an invariant for $e_4$, choose $v_4$ and
$v_5$ on the segment $l_{35}$ connecting $x_3$ and $x_5$
giving us the correct invariant.  The solutions are parametrized
by $\bR$.

Thus 
\[\mathcal{F}: \mathcal{C}(P) \ra \mathcal{C}(\partial P) 
= (0,1)_{e_2} \times (0,1)_{e_4} \times (0,1)_{e_6}\] is 
a principal $\bR$-fibration, and $\mathcal{C}(P)$ is a $4$-cell.

\subsection{Quadrilaterals (D3)} 
Let $P$ be a quadrilateral, and let $e'_1$ be the boundary $1$-orbifold and 
$e'_2, e'_3,$ and $e'_4$ the edges of $P$ so that they are ordered 
with respect to an orientation. In the universal 
cover $\tilde P$, consider a disk $D$
bounded by geodesic segments $e_1, e_2, e_3,$ and $e_4$ mapping 
to $e'_1, e'_2, e'_3,$ and $e'_4$ respectively. 
Let $v_1, v_2, v_3,$ and $v_4$ 
denote the vertices so that $v_i$ and $v_{i+1}$ are the vertices 
of edges $e_i$ for cyclic indices. Post-composing $\dev$
with a projective automorphism if necessary, 
$v_1, v_2, v_3,$ and $v_4$ map to 
\[[0,0,1], [1,0,1], [1,1,1], [0,1,1]\]
respectively. 
Identify $D$ and associated objects 
by $\dev$ as before. 
%   this is repetitive, and by now doesn't need to be said.
%   Why not put this in the beginning once and for all?
% C. I have to do this unfortunately...
Let $l_i$ denote the lines containing $e_i$, $i=1, \dots, 4$.
Let $r_2, r_3,$ and $r_4$ denote the reflections in the holonomy 
group whose lines of fixed points are $l_2$, $l_3$ and $l_4$
respectively. Let $x_2, x_3$, and $x_4$ denote the respective 
isolated fixed points. Since $e_1$ correspond to a $1$-orbifold, 
$r_2$ and $r_4$ act on $l_1$. Thus, $x_2, x_4 \in l_1$. 

First, suppose that $n, m\ge 3$.
Then the cross-ratios are as follows:
\begin{equation} \label{eqn:square1}
[\ovl{[0,1,1]x_3}, \ovl{[0,1,1]x_4}; 
\ovl{[0,1,1][1,1,1]}, \ovl{[0,1,1][0,0,1]}]
= \frac{2}{\cos 2\pi/m + 1}, 
\end{equation}
\begin{equation}\label{eqn:square2}
[\ovl{[1,1,1]x_2}, \ovl{[1,1,1]x_3};  
\ovl{[1,1,1][1,0,1]}, \ovl{[1,1,1][0,1,1]}]
=  \frac{2}{\cos 2\pi/n + 1},
\end{equation}
and 
\[[x_2,x_4; [1,0,1],[0,0,1]] \] 
parametrizes  $e_1$. 
In this case, the pentagon with vertices $x_2, [1,0,1], x_3, [1,1,1], x_4$ 
is convex.
For any choice of $x_2, x_4\in l_1$, there is a line through $[1,0,1]$ 
containing $x_3$ satisfying \eqref{eqn:square1}.
Similarly  there is a line through $[1,1,1]$ 
containing  $x_3$ satisfying \eqref{eqn:square2}. 
Thus, any $x_2, x_4$ determines a unique  $x_3$.
The space of the choices of $x_2, x_4$ keeping
the boundary invariant constant is diffeomorphic to $\bR$. 
Thus, $\mathcal{C}(P)$ is a $2$-cell and the map 
\[{\mathcal{F}}: \mathcal{C}(P) \ra (0,1)_{e_1}\] 
is a principal $\bR$-fibration 
for the space $(0,1)_{e_1}$ of invariants of $e_1$.

Suppose that $n= 2$ and $m \geq 3$. While $x_3\in l_2$,
$x_2$ lies on the line containing $[0,1,1]$ and $[1,1,1]$. 
Since $x_2$ lies on $l_1$, it follows that $x_2 = [0,1,0]$. 
Given any $x_4$, there exists $x_3\in l_2$ satisfying
\[ [\ovl{[0,1,1]x_3}, \ovl{[0,1,1]x_4}; 
\ovl{[0,1,1][1,1,1]}, \ovl{[0,1,1][0,0,1]}] 
=  \frac{2}{\cos 2\pi/m + 1}. \]
Thus, $\mathcal{C}(P)$ is diffeomorphic to $\bR$ and 
$\mathcal{F}:\mathcal{C}(P) \ra (0,1)$ is a homeomorphism.

\subsection{Triangles (D4)}
In the last case (D4),  \cite{G:77} implies
that $\mathcal{C}(P)$ is diffeomorphic to $\bR^+$, and earlier
by Kac-Vinberg \cite{KacVin:67}. 
This completes the proof of Proposition \ref{prop:elemdim}. 
\qed

\subsection{Proof of Theorem A or \ref{thm:finaldim}} 
We assume the following inductive assumption:
Let $\Sigma$ be the $2$-orbifold whose components 
have negative Euler characteristic obtained from 
sewing a $2$-orbifold $\Sigma'$ whose components have 
negative Euler characteristic. We suppose that 
$\Sigma'$ satisfies the fibration property. 
If we show that $\Sigma$ has the fibration property,
the proof of Theorem \ref{thm:finaldim} follows from 
Proposition \ref{prop:elemdim}.

%02/11/14/3:05
First, we suppose that $\Sigma$ is obtained from 
$\Sigma'$ by pasting along simple closed curves.
Let $\Sigma'$ have two boundary principal closed geodesics
$l_1$ and $l_2$ corresponding to $l$ in $\Sigma$.
Take a convex $\rpt$-structure on $\Sigma'$ with $l_1$ 
and $l_2$ with matching invariants and find an  
isomorphism between neighborhoods of $l_1$ and $l_2$ 
in an appropriate ambient $2$-orbifold. 
Identifying and truncating these neighborhoods produces
a convex $\rpt$-structure on $\Sigma$. 
There is a principal fibration
\[{\mathcal{F}}': \mathcal{C}(\Sigma') \ra \mathcal{C}(\partial \Sigma')\]
with the action of a cell of dimension equal to
\[\dim \mathcal{C}(\Sigma') - \dim \mathcal{C}(\partial \Sigma').\]

Take the diagonal subset $\tri$ of 
$\mathcal{C}(\partial \Sigma')$ consisting of elements of 
where invariants of $l_1$ agree with that of $l_2$.
Since ${\mathcal{F}}'$ is a fibration, 
${\mathcal{F}}^{\prime -1}(\tri)$ is a cell of dimension two less
than that of $\mathcal{C}(\Sigma')$. 

%Then from ${\mathcal{F}}^{\prime -1}(\tri)$, 
%we construct all elements on $\Sigma$ by the $\bR^2$-space of
%choices of the gluing maps.

Proposition \ref{prop:principal1} states that there exists
an $\bR^2$-action $\Phi$ on $\mathcal{C}(\Sigma)$ 
and a $\Phi$-invariant fibration 
\[{\mathcal{SP}}: \mathcal{C}(\Sigma) \ra {\mathcal{F}}^{\prime, -1}(\tri)\]
where $\Phi$ acts transitively, freely, and properly on the fibers. 
We have the following commutative diagram:
\begin{equation} \label{eqn:fibrations}
\begin{array}{ccccc} 
{\mathcal{C}}(\Sigma)  & \stackrel{\mathcal{F}}{\longrightarrow} & 
{\mathcal{C}}(\partial \Sigma) & = 
& {\mathcal{C}}(\partial \Sigma) \nonumber \\
\downarrow \mathcal{SP} & & & & \downarrow {g} \nonumber \\
{\mathcal{F}^{\prime -1}}(\tri) \subset {\mathcal{C}}(\Sigma') & 
\stackrel{\mathcal{F}'}{\longrightarrow} & 
\tri \subset {\mathcal{C}}(\partial \Sigma') &
\stackrel{f}{\longrightarrow} & {\mathcal{C}}(\partial \Sigma' - l_1 - l_2)
\end{array}
\end{equation}
where $f$ is a function defined on $\tri$ forgetting 
about the values on $l_1$ and $l_2$, and $g$ a function 
sending invariants of $\partial \Sigma$ to corresponding 
invariants of $\partial \Sigma' - l_1 - l_2$. 
Since $f$ forgets the value of the invariants at $l_1$ and $l_2$, 
$f$ is a principal fibration with the $\bR^2$-action.
$g$ is the homeomorphisms since there is a one-to-one 
correspondence between the boundary components. 

From the diagram, $\mathcal{F}$ can be identified 
with $f \circ {\mathcal{F}'} \circ {\mathcal{SP}}$.
Therefore, $\mathcal{F}$ is a fibration
with the action of a cell of dimension 
\[ 2 + \dim {\mathcal{C}(\Sigma')} - 
\dim {\mathcal{C}(\partial \Sigma')} +2.\] 

Since ${\mathcal{F}^{\prime -1}}(\tri)$ is 
of codimension $2$ in the ambient space, 
$\mathcal{C}(\Sigma)$ is a cell of dimension $\dim \mathcal{C}(\Sigma')$. 
We verify the dimension formula:
The Euler characteristic of the underlying space of $\Sigma'$ 
equals $\chi(\Sigma)$. 
Furthermore $\Sigma'$ and $\Sigma$ have equal numbers of 
cone-points, corner-reflectors (of equal orders), and 
boundary full $1$-orbifolds respectively. 
Thus, the fibration property (i) holds for $\Sigma'$. 

Since we loose two boundary components $l_1$ and $l_2$, we obtain
$\dim {\mathcal{C}}(\partial \Sigma) = 
\dim {\mathcal{C}}(\partial \Sigma') - 4.$
Since the above cell is of dimension 
$\dim {\mathcal{C}(\Sigma)} - 
\dim {\mathcal{C}(\partial \Sigma)}$, the property (ii) holds 
for $\Sigma$.

Suppose that $\Sigma$ is obtained from $\Sigma'$ by 
cross-capping along a simple closed curve $l'$.
By Proposition \ref{prop:principal1}, $\mathcal{SP}$ maps 
$\mathcal{C}(\Sigma)$ diffeomorphically to $\mathcal{C}(\Sigma')$. 
The commutative diagram
\begin{equation}\label{eqn:fibrations2}
\begin{array}{ccccc} 
{\mathcal{C}}(\Sigma)  & \stackrel{\mathcal{F}}{\longrightarrow} & 
{\mathcal{C}}(\partial \Sigma) & = 
& {\mathcal{C}}(\partial \Sigma) \nonumber \\
\downarrow \mathcal{SP} & & & & \downarrow {g} \nonumber \\
{\mathcal{C}}(\Sigma') & 
\stackrel{\mathcal{F}'}{\longrightarrow} & 
{\mathcal{C}}(\partial \Sigma') &
\stackrel{f}{\longrightarrow} & {\mathcal{C}}(\partial \Sigma' - l'),
\end{array}
\end{equation}
where $f$ is the forgetting principal $\bR^2$-fibration and 
$g$ a homeomorphism, holds. Therefore (i) and (ii) easily follow
by verifying the dimension formula.

When $\Sigma$ is obtained from $\Sigma'$ be silvering a simple 
closed curve $l'$, the proof is exactly same as the above
if we change the meanings of the symbols correspondingly to
the silvering case.

Suppose that $\Sigma$ is obtained from $\Sigma'$ by 
folding along a simple closed curve $l'$.
Then $l'$ is purely hyperbolic.  Since the subspace 
of $\mathcal{R}$ consisting of purely hyperbolic elements is
diffeomorphic to an arc, the subspace $\tri$ of
foldable convex $\rpt$-structures is a cell of dimension 
$\dim \mathcal{C}(\Sigma') -1$ by the fibration property.  
The commutative diagram above \eqref{eqn:fibrations} 
holds where $f$ is a principal $\bR$-fibration and 
$g$ is a homeomorphism again. 
Since $\mathcal{SP}$ is a principal $\bR$-fibration, 
(i) and (ii) easily follow.

Suppose that $\Sigma$ is obtained from $\Sigma'$ by
pasting along two full $1$-orbifolds $l_1$ and $l_2$.
Again (i) and (ii) follow by the commutative 
diagram \eqref{eqn:fibrations} for this case.

Suppose that $\Sigma$ is obtained from
$\Sigma'$ by folding or silvering a $1$-orbifold $l'$. 
There is a unique way to do
this. Thus $\mathcal{C}(\Sigma)$ is diffeomorphic to
$\mathcal{C}(\Sigma')$, and (i) and (ii) follows
by diagram \ref{eqn:fibrations2} for this case.
\qed

\subsection{A theorem of Thurston}
% Recall that we denote by $\mathcal{T}(\Sigma)$ the subspace of
% $\mathcal{C}(\Sigma)$ consisting of equivalence classes 
% of hyperbolic $\rpt$-structures, and the subspace topology
% is equal to the usual topology of the Teichm\"uller space.
%
% I think this has already been said.
% C. delete

Since the holonomy of hyperbolic $\rpt$-surfaces are hyperbolic, 
the space of invariants of a boundary closed curve
is a one-dimensional subspace of $\mathcal{R}$ diffeomorphic 
to $\bR$: The lengths of the closed curves provide satisfactory invariants. 
Invariants of boundary full $1$-orbifolds are again simply lengths. 
Therefore, we define $\mathcal{T}(\partial \Sigma)$ to 
be the product of these lines of invariants. 

For the sake of completeness % , but not for the necessity of this paper,
we include the projective proof of the following theorem:
\begin{thm}[Thurston]\label{thm:Thur}
Let $\Sigma$ be a $2$-orbifold of negative Euler characteristic
and empty boundary. Then the deformation space
$\mathcal{T}(\Sigma)$ of hyperbolic $\rpt$-structures 
is diffeomorphic to a cell of dimension 
$-3\chi(X_\Sigma) + 2k + l$
where $X_\Sigma$ is the underlying space, and $k$
is the number of cone-points and $l$ is the number
of corner-reflectors.
\end{thm}
\begin{proof}
Since we repeat the proof of Theorem \ref{thm:finaldim},
we give a sketchy argument.

Let a $2$-orbifold $\Sigma$, each component of which
has negative Euler characteristic, 
be in a class $\mathcal{P}$ if 
the following hold:
\begin{itemize}
\item[(i)] The deformation space of hyperbolic 
$\rpt$-structures $\mathcal{T}(\Sigma)$ is diffeomorphic to
a cell of dimension 
\[-3\chi(X_\Sigma) + 2k + l + 2n\]
where $k$ is the number of cone-points,
$l$ the number of corner-reflectors,
$n$ is the number of boundary full $1$-orbifolds.
\item[(ii)] There exists a principal fibration
\[{\mathcal{F}}: \mathcal{T}(\Sigma) \ra \mathcal{T}(\partial \Sigma)\] 
with the action by a cell of dimension 
$\dim \mathcal{T}(\Sigma) - \dim \mathcal{T}(\partial \Sigma)$.
\end{itemize}

%%%Elementary case.....
%For elementary orbifolds, (i) and (ii) hold
%by the results in \S\ref{sec:teichmuller}.

Let $\Sigma$ be a $2$-orbifold whose components 
are orbifolds of negative Euler characteristic, 
and it splits into an orbifold $\Sigma'$ in $\mathcal{P}$.
We suppose that (i) and (ii)
hold for $\Sigma'$, and show that (i) and (ii) hold for 
$\Sigma$. Since $\Sigma$ eventually decomposes 
into a union of elementary $2$-orbifolds where (i) and (ii) hold,
we would have completed the proof by Proposition 
\ref{prop:hypele}.

Since we need to match lengths and find gluing maps
preserving lengths, the arguments are similar 
to the rest of the proof of Theorem \ref{thm:finaldim}.

We need the result of Proposition \ref{prop:principal1} 
for hyperbolic cases:
\begin{description} 
\item[(A)(I)(1)] 
Let the $2$-orbifold $\Sigma''$ be obtained 
from pasting along two closed curves $b, b'$ 
in a $2$-orbifold $\Sigma'$.
The map resulting from splitting
\[\mathcal{SP}: \mathcal{T}(\Sigma'') \ra 
\Delta \subset \mathcal{T}(\Sigma')\]
is a principal $\bR$-fibration, 
where $\Delta$ is the subset of $\mathcal{C}(\Sigma')$ 
where $b$ and $b'$ have equal invariants.
\item[(A)(I)(2)] 
Let $\Sigma''$ be obtained from $\Sigma'$ by cross-capping.
The resulting map
\[\mathcal{SP}: \mathcal{T}(\Sigma'') \ra \mathcal{T}(\Sigma')\]
is a diffeomorphism.
\item[(A)(II)(1)] Let $\Sigma''$ be obtained from $\Sigma'$ by 
silvering. The clarifying map
\[\mathcal{SP}: \mathcal{T}(\Sigma'') \ra \mathcal{T}(\Sigma')\]
is a diffeomorphism.
\item[(A)(II)(2)] Let $\Sigma''$ be obtained from $\Sigma'$ 
by folding a boundary closed curve $l'$. 
The unfolding map 
\[\mathcal{SP}: \mathcal{T}(\Sigma'') \ra \Delta \subset
\mathcal{T}(\Sigma')\]
is a principal $\bR$-fibration. 
\item[(B)(I)] Let $\Sigma''$ be obtained by pasting along 
two full $1$-orbifolds $b$ and $b'$ in $\Sigma'$.
The splitting map
\[\mathcal{SP}: \mathcal{T}(\Sigma'') \ra 
\Delta \subset \mathcal{T}(\Sigma')\] 
is a diffeomorphism where
$\Delta$ is a subset of $\mathcal{T}(\Sigma')$ where 
the invariants of $b$ and $b'$ are equal.
\item[(B)(II)] Let $\Sigma''$ be obtained by silvering 
or folding a full $1$-orbifold.
The clarifying or unfolding map
\[\mathcal{SP}: \mathcal{T}(\Sigma'') \ra \mathcal{T}(\Sigma')\]
is a diffeomorphism.
\end{description} 

Again, using diagrams \eqref{eqn:fibrations} 
and \eqref{eqn:fibrations2} adopted to each cases, 
we prove the theorem.
\end{proof}

\begin{appendix}
\section{Developing maps of elementary 
orbifolds} 

This appendix presents various developing maps of
elementary $2$-orbifolds which are generated by Maple. 
We will divide the orbifolds into various
one or two polygonal domains and develop them.  In the following, 
the {\em depth} means the maximum word length of the deck transformations
written using the standard generators. $(s, t)$ indicates certain
invariants analogous to the invariants for a pair-of-pants as given in
\cite{Gconv:90} (see \S \ref{sec:deformelement}).  
The pictures for (P1) are not given, and examples for
(D4) are given in \cite{Gconv:90}.

\vfill

\break

\begin{figure}[h]
\centerline{\epsfxsize=10.2cm \epsfbox{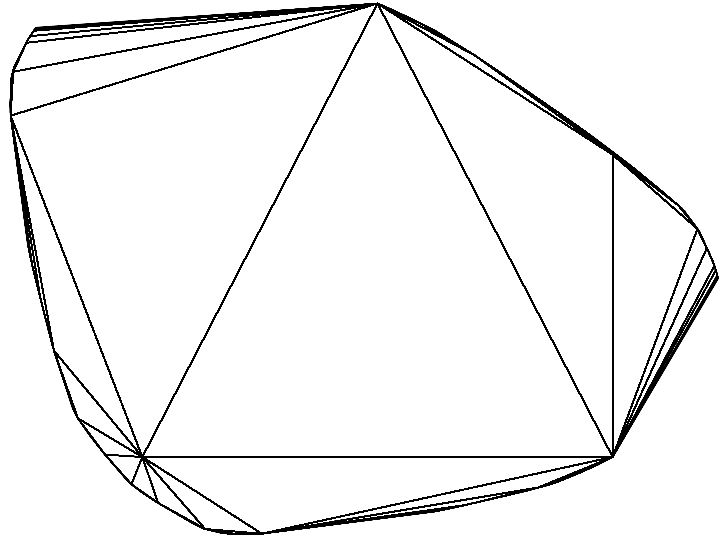}}
\caption{\label{fig:fig1} 
An annulus with a cone-point of order $5$, boundary invariants $(1/3.1, 4.1)$
and $(1/4.1, 5.1)$, $(s,t)=(2,1)$, depth $4$,
type (P2), and symbol $A(;5)$.}
\end{figure}
\typeout{<<fig1.eps>>}

\begin{figure}[h]
\centerline{\epsfxsize=10.2cm \epsfbox{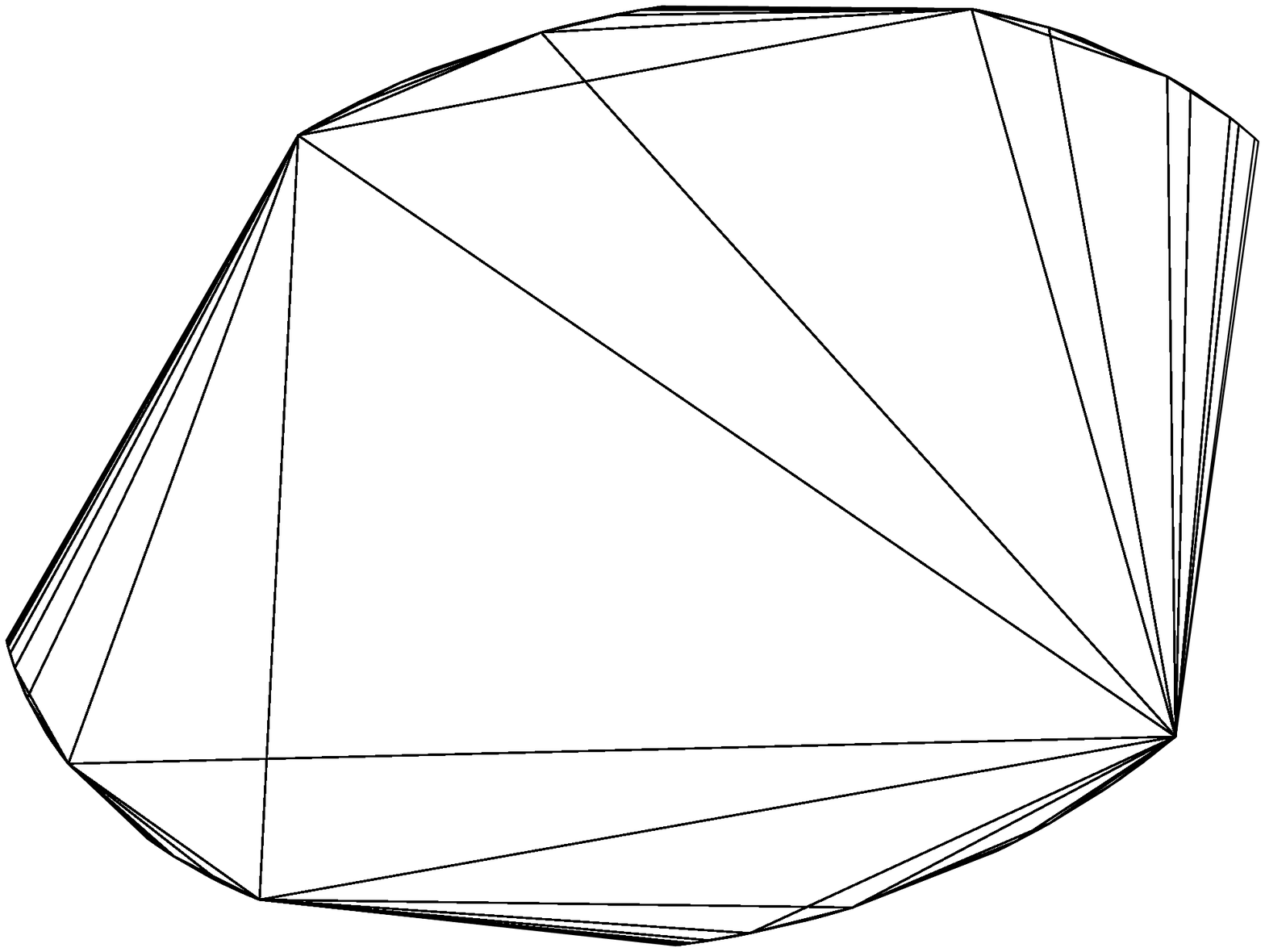}}
\caption{\label{fig:fig2} 
An annulus with a cone-point of order $2$, boundary
invariants $(1/2, 3)$, $(1/4, 5)$, $(s, t)=(2,1)$, 
depth $4$, type (P2), and symbol $A(;2)$.
%\marginpar{Check this $s,t$}
}
\end{figure}
\typeout{<<fig2.eps>>}

\vfill
\break

\begin{figure}[h]
\centerline{\epsfxsize=11cm \epsfbox{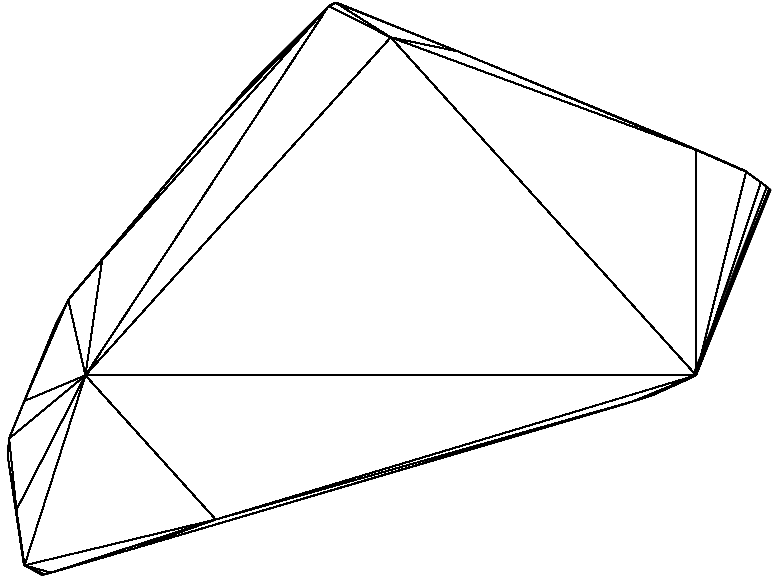}}
\caption{\label{fig:fig3} 
A disk with two cone-points of orders $3$ and $5$,
the boundary invariant $(1/5,6)$, $(s,t)=(1,1)$,
depth $4$, type (P3), and symbol $D(;3,5)$.
}
\end{figure}
\typeout{<<fig3.eps>>}

\begin{figure}[h]
\centerline{\epsfxsize=11cm \epsfbox{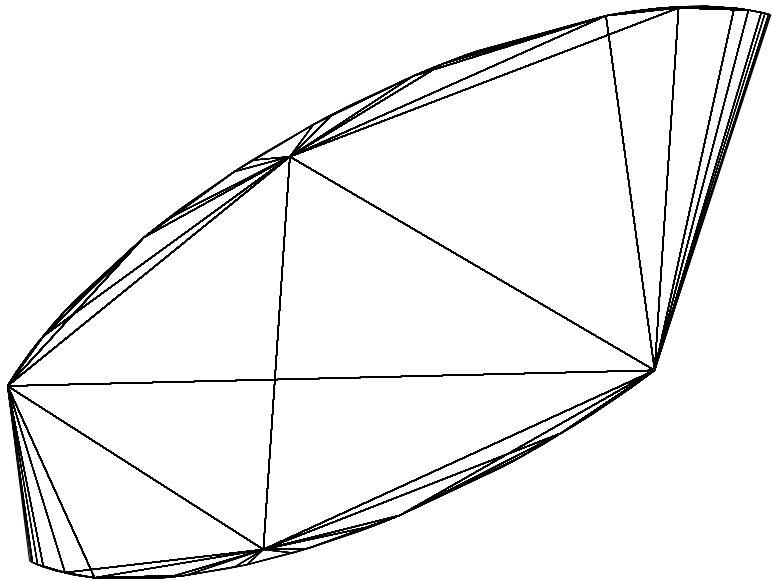}}
\caption{\label{fig:fig4} 
A disk with cone-points of orders $2$ and $7$, the boundary
invariant $(1/3, 4)$, depth $4$,
type (P3), and symbol $D(;2,7)$.
%\marginpar{$s,t$?}
}
\end{figure}
\typeout{<<fig4.eps>>}

\vfill
\break

\begin{figure}[h]
\centerline{\epsfxsize=11cm \epsfbox{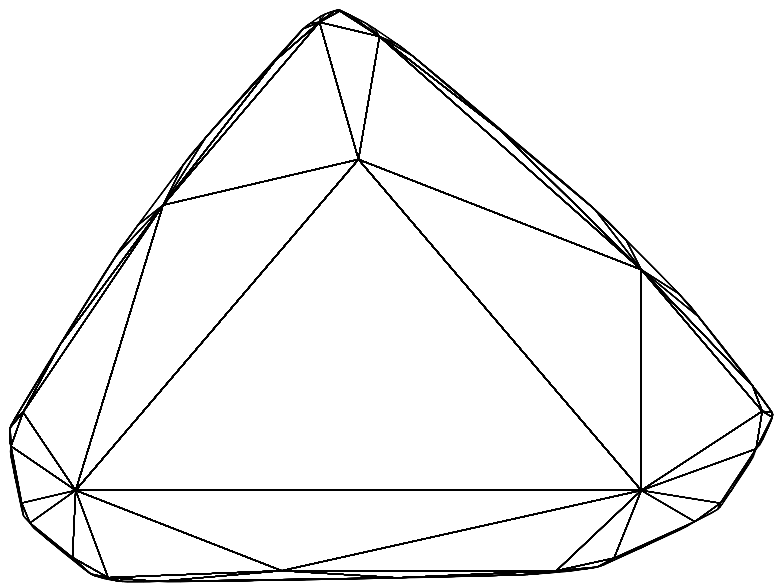}}
\caption{\label{fig:fig5} 
A sphere with three cone-points of orders $3,5,5$. $(s,t) = (2,2)$. 
depth $4$, type (P4), and symbol $\SI^2(;3,5,5)$.
}
\end{figure}
\typeout{<<fig5.eps>>}

\begin{figure}[h]
\centerline{\epsfxsize=11cm \epsfbox{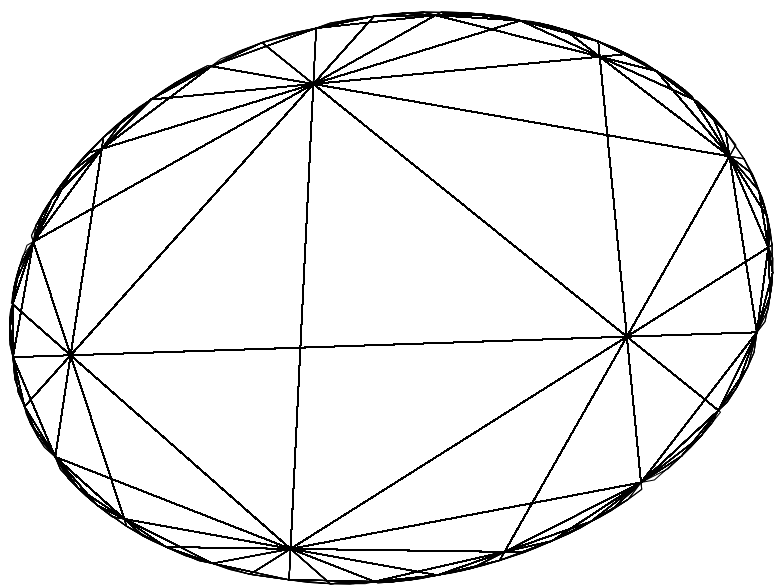}}
\caption{\label{fig:fig6} 
A sphere with order two cone-points orders $2,5,7$
depth $4$, type (P4), and symbol $\SI^2(;2,5,7)$.
}
\end{figure}
\typeout{<<fig6.eps>>}

\vfill
\break

\begin{figure}[h]
\centerline{\epsfxsize=11cm \epsfbox{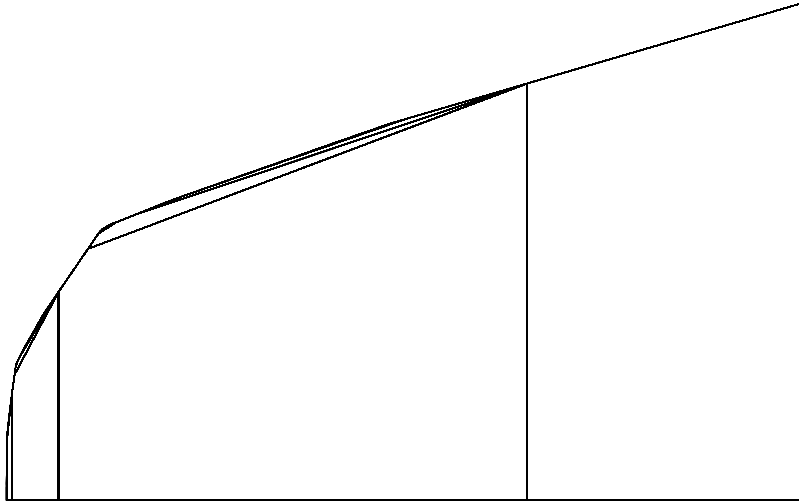}}
\caption{\label{fig:fig7}
An annulus with a boundary $1$-orbifold, boundary invariants $0.4$ 
and $(0.5, 5)$, depth $5$, type (A1), and symbol $A(2,2; )$. 
}
\end{figure}
\typeout{<<fig7.eps>>}

\begin{figure}[h]
\centerline{\epsfxsize=11cm \epsfbox{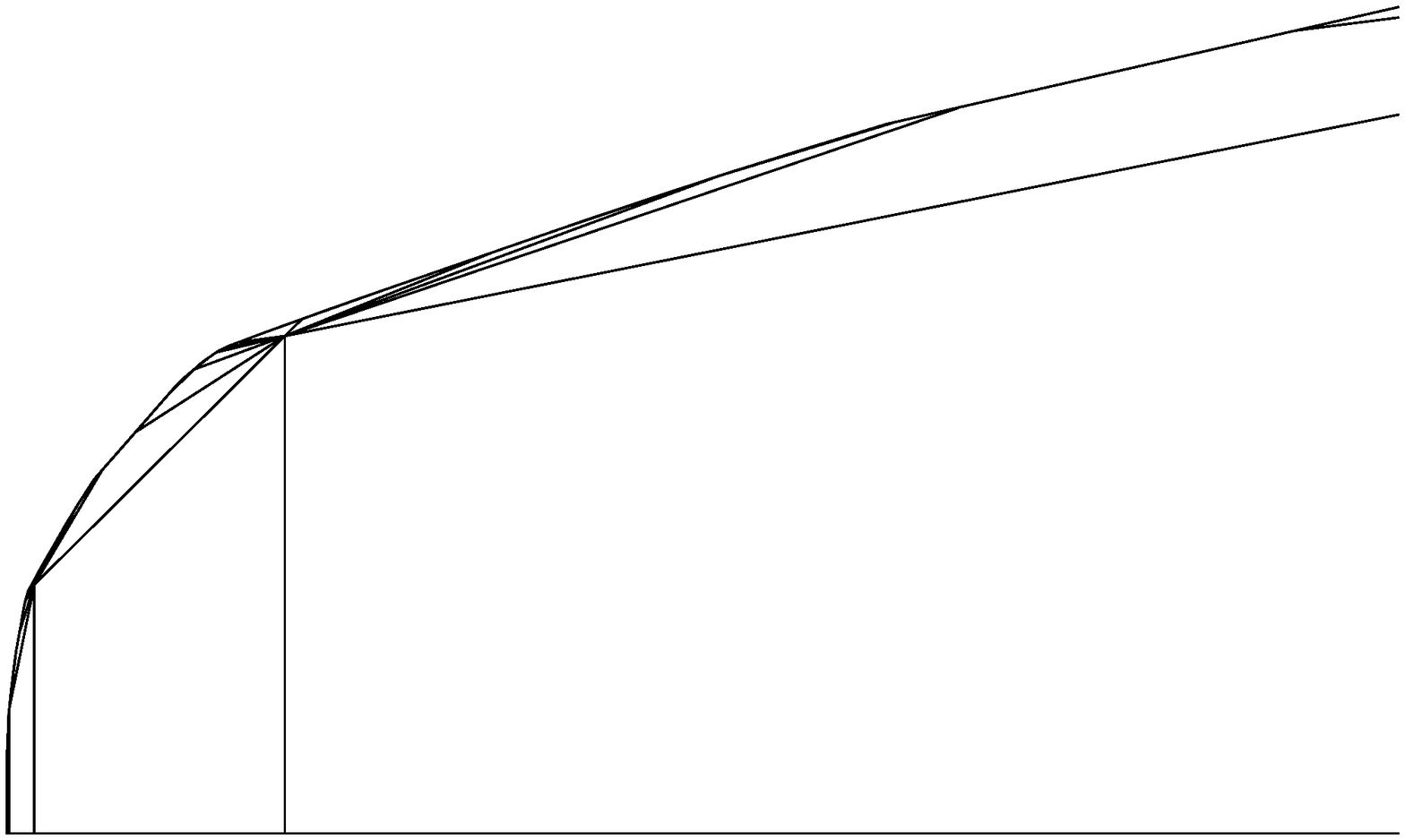}}
\caption{\label{fig:fig8}
An annulus with a corner-reflector of order $3$ , the boundary 
invariant  $(0.5, 5)$, depth $5$, type (A2), and
symbol $A(3; )$.
}
\end{figure}
\typeout{<<fig8.eps>>}

\vfill
\break

\begin{figure}[h]
\centerline{\epsfxsize=11cm \epsfbox{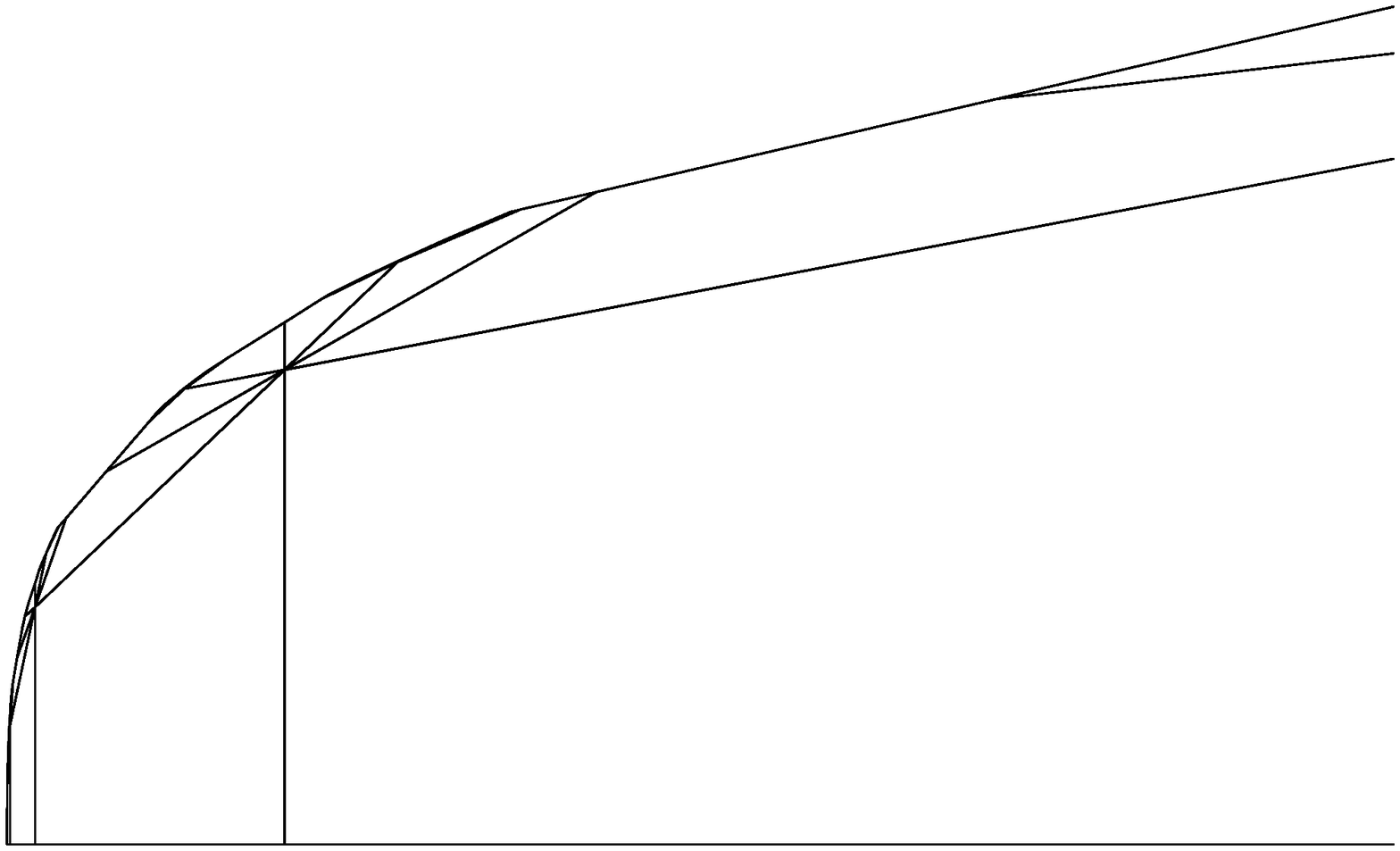}}
\caption{\label{fig:fig9}
An annulus with a corner-reflector of order  $2$, 
the boundary invariant  $(0.5, 5)$, depth $5$, type (A2), 
and symbol $A(2; )$
}
\end{figure}
\typeout{<<fig9.eps>>}

\begin{figure}[h]
\centerline{\epsfxsize=11cm \epsfbox{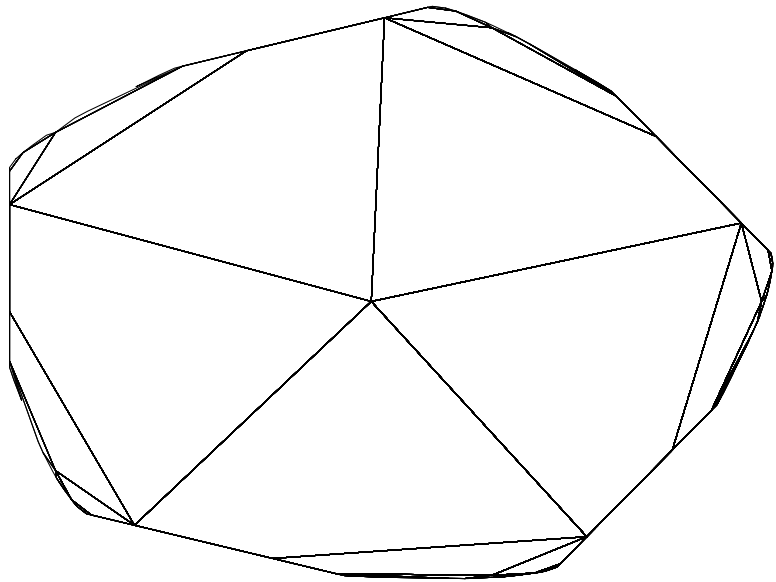}}
\caption{\label{fig:fig15}
A disk with one s-segment, and boundary 1-orbifold and a cone-point of
order $5$, boundary invariant $0.3$, depth $4$, type (A3),
and symbol $D^2(2,2;5)$.
%\marginpar{symbol?}
}
\end{figure}
\typeout{<<fig15.eps>>}

\vfill
\break

\begin{figure}[h]
\centerline{\epsfxsize=11cm \epsfbox{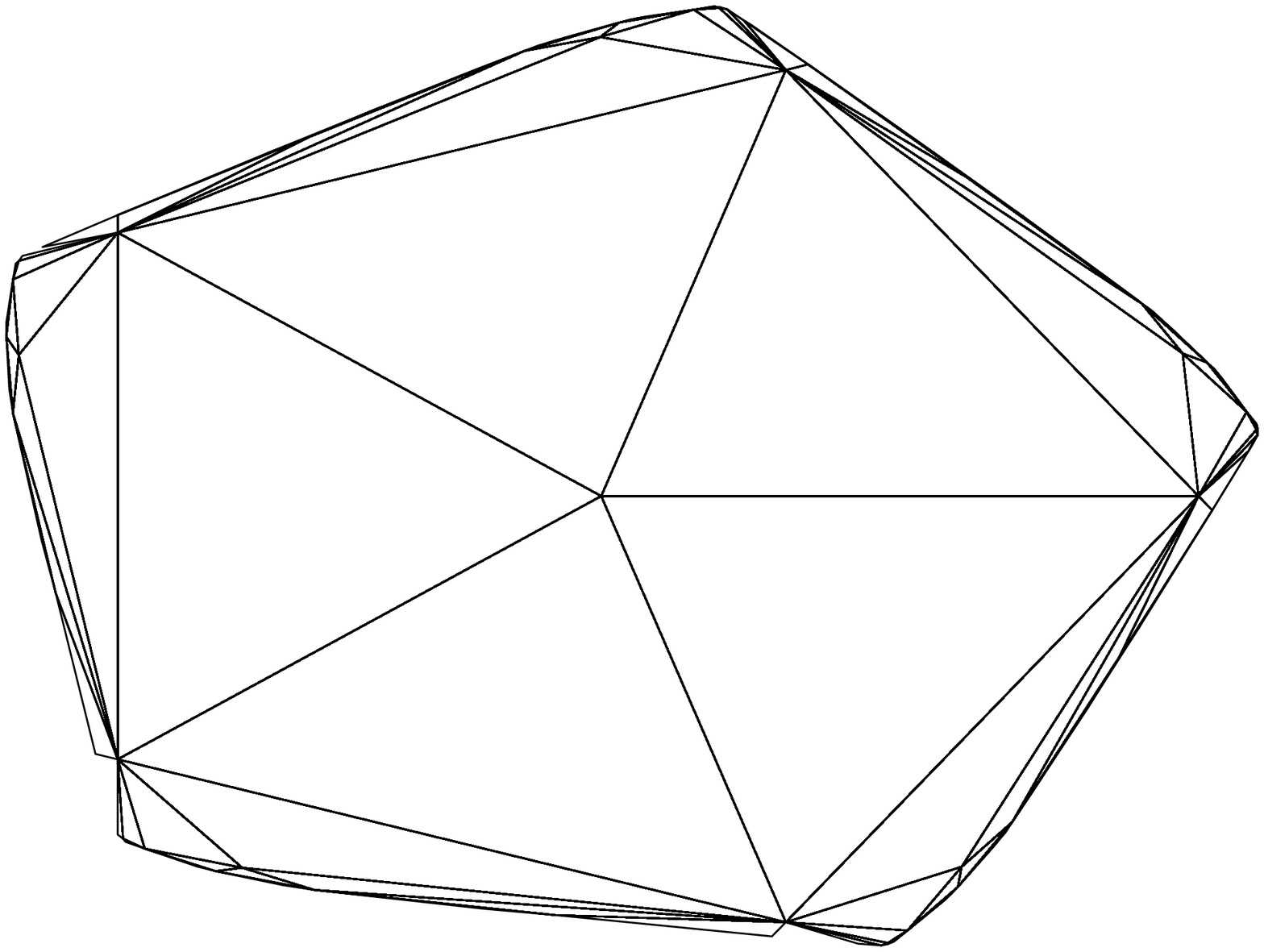}}
\caption{\label{fig:fig16}
A disk with one corner-reflector of order $3$ and one
cone-point of order $5$, depth $6$, type (A4), and
symbol $D^2(3;5)$.
%\marginpar{symbol?}
}
\end{figure}
\typeout{<<fig16.eps>>}

\begin{figure}[h]
\centerline{\epsfxsize=11cm \epsfbox{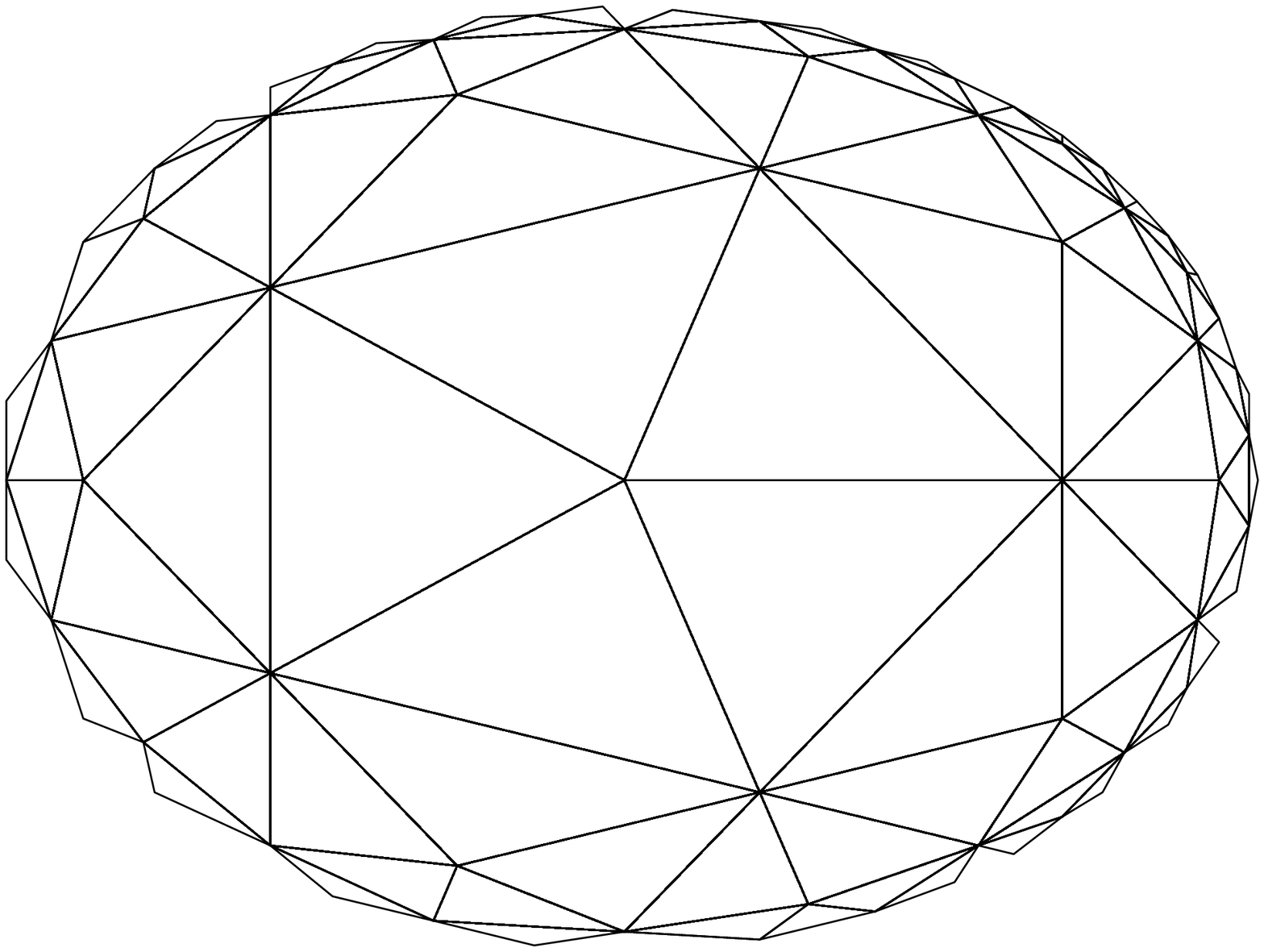}}
\caption{\label{fig:fig17}
A disk with one corner-reflector of order $2$ and one
cone-point of order $5$, depth $6$, type (A4),
and symbol $D^2(2;5)$.
%\marginpar{symbol?}
}
\end{figure}
\typeout{<<fig17.eps>>}

\vfill
\break

\begin{figure}[h]
\centerline{\epsfxsize=11cm \epsfbox{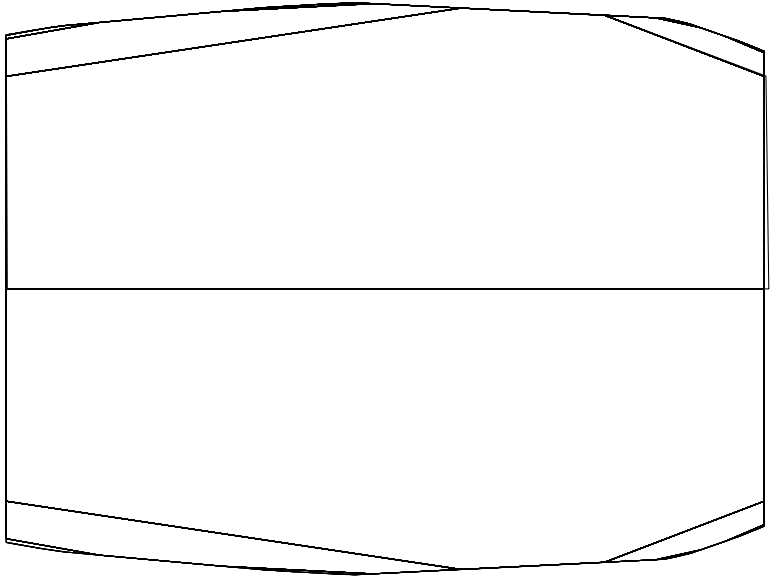}}
\caption{\label{fig:fig12}
A hexagon, boundary invariants $0.2, 0.4, 0.3$,
depth $5$, type (D1), and 
symbol $D^2(2,2,2,2,2,2;)$.
}
\end{figure}
\typeout{<<fig12.eps>>}

\begin{figure}[h]
\centerline{\epsfxsize=11cm \epsfbox{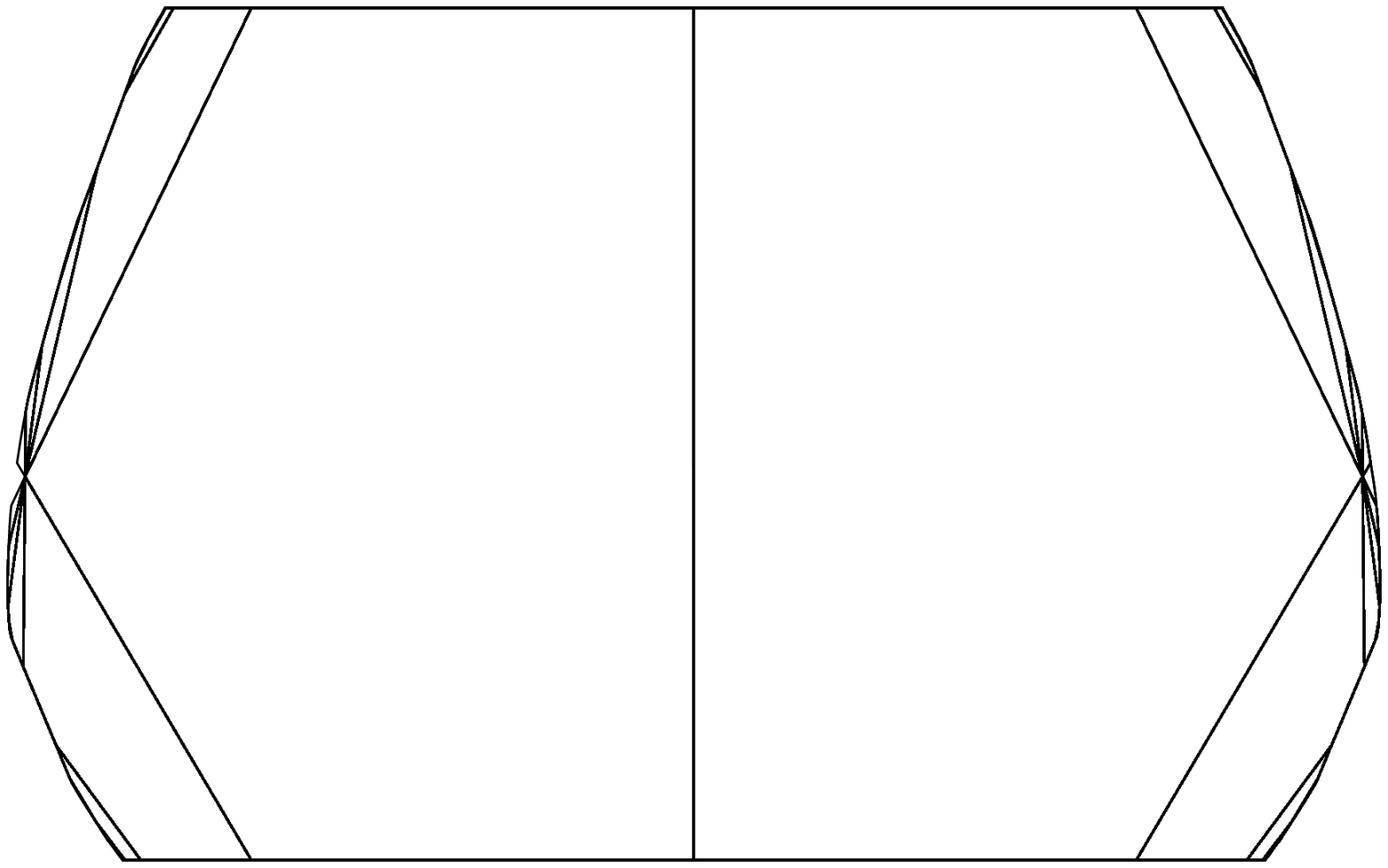}}
\caption{\label{fig:fig10}
A pentagon with a corner-reflector of order $5$,
boundary invariants $0.4, 0.3$, depth $5$, type (D2),
and symbol $D^2(2,2,2,2,5;)$.
}
\end{figure}
\typeout{<<fig10.eps>>}

\vfill
\break

\begin{figure}[h]
\centerline{\epsfxsize=11cm \epsfbox{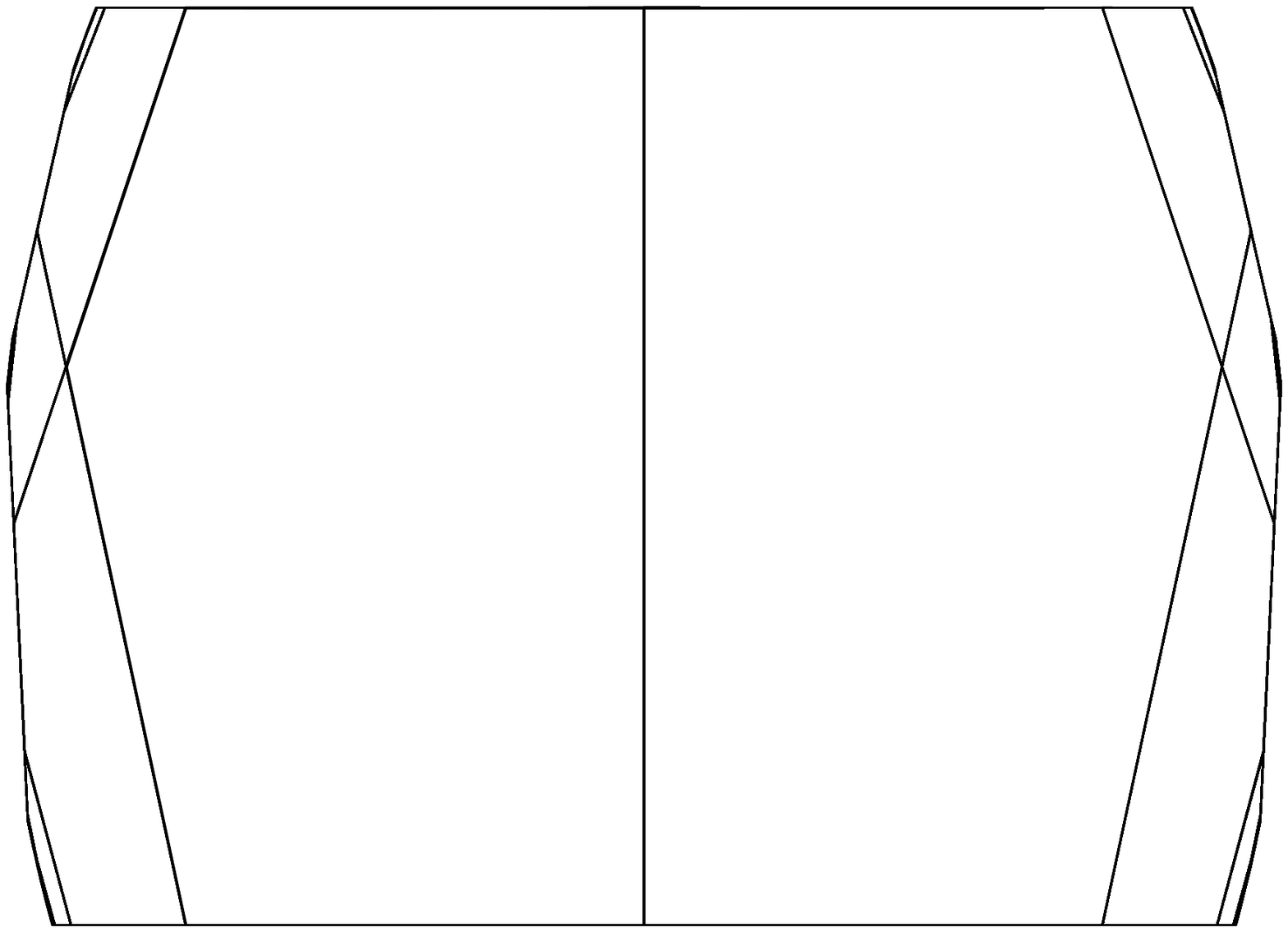}}
\caption{\label{fig:fig11}
A pentagon with a corner-reflector order $2$,
boundary invariants $0.4, 0.3$, depth $5$, type (D2),
and symbol $D^2(2,2,2,2,2;)$.
}
\end{figure}
\typeout{<<fig11.eps>>}

\begin{figure}[h]
\centerline{\epsfxsize=11cm \epsfbox{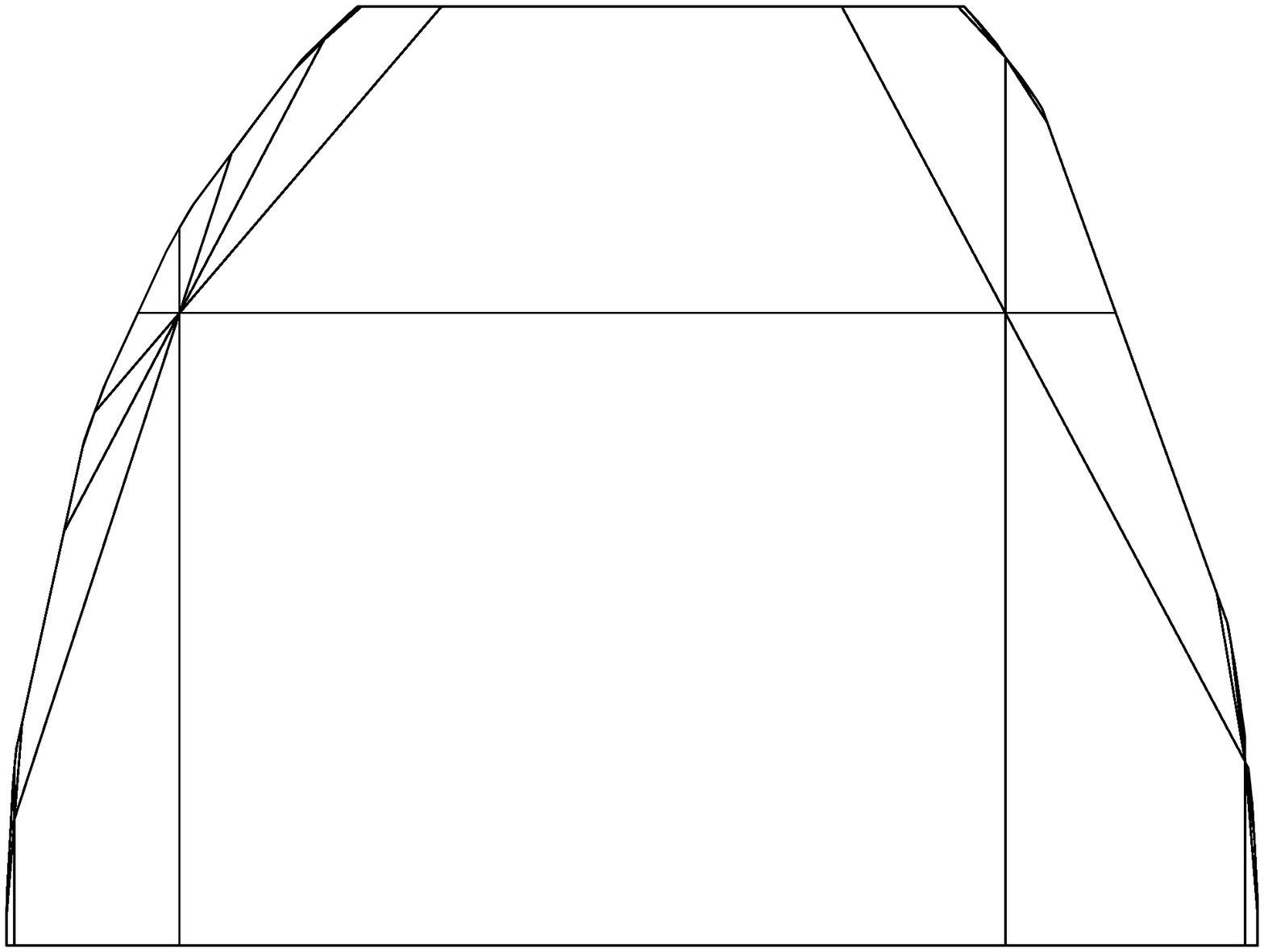}}
\caption{\label{fig:fig13}
A square with corner reflectors of order $3$ and $5$, 
the boundary invariant $0.15$, depth $5$, type (D3), 
and symbol $D^2(2,2,3,5;)$.
}
\end{figure}
\typeout{<<fig13.eps>>}

\vfill
\break

\begin{figure}[h]
\centerline{\epsfxsize=11cm \epsfbox{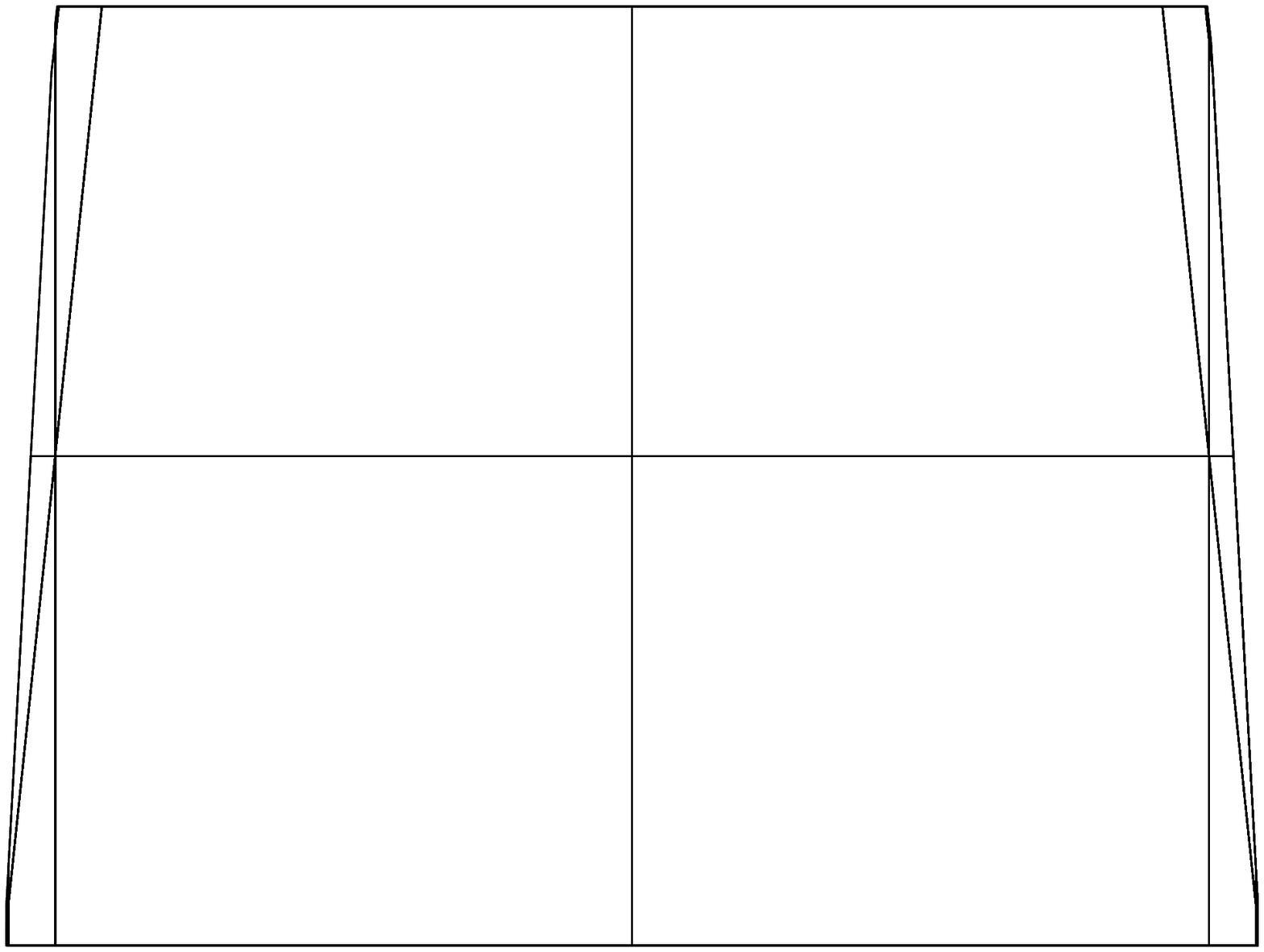}}
\caption{\label{fig:fig14}
A square with corner reflectors of order $2$ and $3$, 
the boundary invariant $0.15$, depth $4$, type (D3),
and symbol $D^2(2,2,2,3;)$.
}
\end{figure}
\typeout{<<fig14.eps>>}

\end{appendix}

\vfill

\break

\bibliographystyle{plain}

%Introduce duality involutions.
%Maybe in the introduction?

\end{document}